\documentclass[psamsfonts]{amsart}
\usepackage{amsfonts, mathrsfs}
\usepackage{ifthen}
\usepackage{amsthm}
\usepackage{amsmath}
\usepackage{graphicx}
\usepackage{amscd,amssymb,amsthm}
\usepackage{graphicx}
\usepackage{epstopdf}
\usepackage{hyperref}
\usepackage{algorithm2e}
\usepackage{mathtools}

\newcounter{minutes}
\divide\time by 60
\newcounter{hours}
\multiply\time by 60 \addtocounter{minutes}{-\time}

\newcommand{\real}{\operatorname{Re}}
\newcommand{\imag}{\operatorname{Im}}

\newtheorem{lemma}{Lemma}

\newtheorem{corollary}{Corollary}
\newtheorem{theorem}{Theorem}

\usepackage{xcolor}
\begin{document}

\title[Infinitely divisible modified Bessel distributions]{Infinitely divisible modified Bessel distributions}

\author[\'A. Baricz]{\'Arp\'ad Baricz}
\address{Department of Economics, Babe\c{s}-Bolyai University, 400591 Cluj-Napoca, Romania and Institute of Applied Mathematics, \'Obuda University, 1034 Budapest, Hungary}
\email{bariczocsi@yahoo.com}

\author[D.K. Prabhu]{Dhivya Prabhu K}
\address{Department of Mathematics,
Indian Institute of Technology Indore, Indore 453552, India}
\email{phd2001241004@iiti.ac.in}

\author[S. Singh]{Sanjeev Singh}
\address{Department of Mathematics,
Indian Institute of Technology Indore, Indore 453552, India}
\email{snjvsngh@iiti.ac.in}

\author[V.A. Vijesh]{Antony Vijesh V}
\address{Department of Mathematics,
Indian Institute of Technology Indore, Indore 453552, India}
\email{vijesh@iiti.ac.in}

\def\thefootnote{}
\footnotetext{ \texttt{File:~\jobname .tex,
         printed: \number\year-0\number\month-0\number\day,
          \thehours.\ifnum\theminutes<10{0}\fi\theminutes}
} \makeatletter\def\thefootnote{\@arabic\c@footnote}\makeatother

\dedicatory{\'A. Baricz dedicates this paper to Mourad E.H. Ismail on the occasion of his 80th birthday}

\maketitle

\begin{abstract}
In this paper we focus on continuous univariate probability distributions, like McKay distributions, $K$-distribution, generalized inverse Gaussian distribution and generalised McKay distributions, with support $[0,\infty),$ which are related to modified Bessel functions of the first and second kinds and in most cases we show that they belong to the class infinitely divisible distributions, self-decomposable distributions, generalized gamma convolutions and hyperbolically completely monotone densities. Some of the results are known, however the proofs are new and we use special functions technique. Integral representations of quotients of Tricomi hypergeometric functions as well as of quotients of Gaussian hypergeometric functions, or modified Bessel functions of the second kind play an important role in our study. In addition, by using a different approach we rediscover a Stieltjes transform representation due to Hermann Hankel for the product of modified Bessel functions of the first and second kinds and we also deduce a series of new Stieltjes transform representations for products, quotients and their reciprocals concerning modified Bessel functions of the first and second kinds. By using these results we obtain new infinitely divisible modified Bessel distributions with Laplace transforms related to modified Bessel functions of the first and second kind. Moreover, we show that the new Stieltjes transform representations have some interesting applications and we list some open problems which may be of interest for further research. In addition, we present a new proof via the Pick function characterization theorem for the infinite divisibility of the ratio of two gamma random variables and we present some new Stieltjes transform representations of quotients of Tricomi hypergeometric functions.
\end{abstract}

\section{\bf Introduction}

\subsection{Preliminaries on infinite divisibility} In probability theory, a {\em probability distribution is infinitely divisible} (see Steutel and Van Harn \cite{steutel}) if it can be expressed as the probability distribution of the sum of an arbitrary number of independent and identically distributed random variables. The concept of infinite divisibility of probability distributions was introduced in 1929 by Bruno de Finetti and the most basic results were discovered by Andrey Kolmogorov, Paul L\'evy and Aleksandr Khinchin. This type of decomposition of a distribution is used in probability and statistics to find families of probability distributions that might be natural choices for certain models or applications. Infinitely divisible distributions play also an important role in probability theory in the context of limit theorems. The characteristic function of any infinitely divisible distribution is then called an infinitely divisible characteristic function and such a function may be represented, for any value of $n$ as the $n$th power of some other characteristic function. More precisely, a probability distribution $\nu$ on the half-line $[0,\infty)$ is infinitely divisible if for any $n\in\mathbb{N}$ there exists a probability distribution $\nu_n$ on $(0,\infty)$ such that
$$\int_0^{\infty}e^{-xt}d\nu=\left[\int_0^{\infty}e^{-xt}d\nu_n\right]^n.$$
A function $f:(0,\infty)\to\mathbb{R}$ is {\em completely monotone} (or completely monotonic) if it has derivatives of all orders and $(-1)^nf^{(n)}(x)>0$ for all $x>0$ and $n\in\{0,1,2,\dots\}.$ The classes of completely
monotonic functions and infinitely divisible distributions are related by the following well-known result (see \cite[p. 425]{feller}).

\begin{lemma}\label{lemma1}
The function $\omega:(0,\infty)\to(0,\infty)$ is the Laplace transform of an infinitely divisible distribution if and only if $\omega(x)=e^{-\varphi(x)},$ where $\varphi(0^{+})=0$ and $\varphi:(0,\infty)\to(0,\infty)$ is a Bernstein function, that is, $\varphi'$ is completely monotonic.
\end{lemma}

It is important to mention here that every continuous-time L\'evy process has distributions that are necessarily infinitely divisible, and conversely every infinitely divisible distribution generates uniquely
a L\'evy process (see for example Steutel and Van Harn \cite{steutel}). Moreover, in various real life situations some concrete models require a random effect to be the sum of several independent random
components with the same distribution. In this kind of situations a very convenient way is to suppose the infinite divisibility of the distribution of these random effects. Similar situations may occur in biology, physics, economics and insurance.

The concept of {\em self-decomposability of probability measures} is due to Paul L\'evy and goes back to 1937. A random variable $X,$ distributed according to a law, is self-decomposable, if and only if for every $c\in(0,1),$ there exists a random variable $X_c,$ independent of $X,$ such that $X$ and $X_c+cX$ are equal in law. We know that a distribution is self-decomposable if and only if it is a weak limit of partial normed centered sums of simple sequence of independent random variables. More precisely, a probability distribution is self-decomposable if it is the limit of
$$\left(X_1+\dots+X_n-b_n\right)/a_n,$$
where the $X_i$'s are independent random variables and $\{a_n\}$ and $\{b_n\}$ are sequences of constants with $a_n\to\infty,$ $a_{n+1}/a_n\to1.$ It is also well known that every self-decomposable distribution is infinitely divisible, and the class of self-decomposable distributions is closed under the convolution and the weak convergence. Moreover, this class contains stable distributions and generalized gamma convolutions. The next auxiliary result (see \cite[p. 589]{feller}) is a characterization of self-decomposable distributions with support $[0,\infty).$

\begin{lemma}\label{lemma2}
A random variable $X$ with support $[0,\infty)$ having the Laplace transform $\omega:(0,\infty)\to(0,\infty)$ is self-decomposable if and only if for every $\alpha,$ where $\alpha\in(0,1),$ the function $\omega(x)/\omega(\alpha x)$ is a Laplace transform.
\end{lemma}

Now, let us focus on other two subclasses of infinitely distributions. By definition, a function $f:(0,\infty)\to(0,\infty)$ is {\em hyperbolically completely monotone} if for each $u>0$ the function $f(uv)f(u/v)$ is completely monotone as a function of $w=v+1/v,$ where $v>0.$ A distribution is said to be hyperbolically completely monotone if its probability density function is hyperbolically completely monotone. This class of lifetime distributions has been discovered by Lennart Bondesson (see \cite{bondesson} and the references therein) and it is strongly connected to the class of generalized gamma convolutions. The class of {\em generalized gamma convolutions} was introduced by Olof Thorin \cite{thorin} in 1977 and it is in fact the smallest class of distributions on $(0,\infty)$ that contains all gamma distributions and is closed under convolution and weak convergence. A positive continuous random variable $X$ belongs to the class of generalized gamma convolutions if its Laplace transform is of the form
$$L(s)=\exp\left(-as+\int_0^{\infty}\ln\left(\frac{t}{s+t}\right)d\mu(t)\right),$$
where $s\geq0,$ $a\geq0$ and $d\mu(t)$ is a non-negative measure. Since the gamma distribution is infinitely divisible and self-decomposable, so is every generalized gamma convolution. Moreover,
an interesting result of Bondesson (see \cite{bondesson}) states that a probability density, which is hyperbolically completely monotone, is the density of a generalized gamma convolution. The following result (see \cite[Theorem 3.1.2]{bondesson}) is a {\em Pick function characterization theorem} of generalized gamma convolutions.

\begin{lemma}\label{lemma3}
A probability distribution of a continuous random variable $X$ on $[0,\infty)$ is a generalized gamma convolution if and only if its moment generating function $\psi$ is analytic and zero-free in $\mathbb{C}\setminus[0,\infty)$ and satisfies $\imag\left[\psi'(s)/\psi(s)\right]\geq0$ for $\imag s>0.$
\end{lemma}

Another simple characterization of generalized gamma convolutions is due to Bondesson (see \cite[Theorem 6.1.1]{bondesson}) and is the following result.

\begin{lemma}\label{lemma4}
A function $\phi$ on $[0,\infty)$ is the Laplace transform of a generalized gamma convolution if and only if $\phi(0)=1$ and $\phi$ is hyperbolically completely monotone.
\end{lemma}

We note that Thorin's class of generalized gamma convolutions is closed with respect to change in scale, weak limits, and addition of independent random variables. Moreover, Lennart Bondesson \cite[Theorem 1]{bonde}
has shown that the generalized gamma convolution class also has the remarkable property of being closed with respect to multiplication of independent random variables.

It is also worth to mention here that to prove or disprove infinite divisibility of a certain distribution is sometimes a complicated task and it may need a specialized approach, see for example the papers of John Kent \cite{kent}, as well as of Pierre Bosch and Thomas Simon \cite{bosch0}, \cite{bosch} and the references therein. The Laplace transforms of probability measures are usually transcendental special functions, which led various authors to study the complete monotonicity of various quotients of special functions (like modified Bessel Bessel functions, Tricomi hypergeometric functions, parabolic cylinder functions) and of the logarithmic derivatives of solutions of differential and difference equations (see for example the papers of Philip Hartman \cite{hartman1} and \cite{hartman2}). This paper is motivated by the {\em special function technique} approach of Mourad Ismail and coauthors (see for example the papers \cite{ismail00}, \cite{ismail0}, \cite{kelker}, \cite{may}, \cite{miller} and \cite{ismail}) and our aim is to continue and complement the results from these papers by deducing a series of new Stieltjes transform representations for products and quotients of modified Bessel functions of the first and second kinds and by obtaining via these results new infinitely divisible modified Bessel distributions. Some of our main results are proved by using ideas from the above mentioned papers of Mourad Ismail, however we also use some ideas from the theory of generalized gamma convolutions and hyperbolically completely monotone densities.

The paper is organized as follows: in the remaining part of this section we recall some basic lemmas on Stieltjes transforms. In Section 2 we present a series of results on infinite divisibility of McKay distributions, $K$-distribution, generalized inverse Gaussian distribution and generalised McKay distributions. Moreover, in Section 3 we obtain a series of new Stieltjes transform representations for products and quotients of modified Bessel functions of the first and second kinds and by using these results we obtain new infinitely divisible modified Bessel distributions. These new distributions have Laplace transforms related to modified Bessel functions of the first and second kind. Section 4 is devoted to remarks and open problems related to the infinite divisibility of some distributions related to modified Bessel functions and Tricomi hypergeometric functions, while Section 5 contains all the proofs of the main results of this paper.

\subsection{Preliminaries on Stieltjes transforms} Before we present a series of results concerning distributions whose probability density function or Laplace transform involves modified Bessel functions we recall some basic results concerning Stieltjes transforms. The first two lemmas of this subsection are two variants of the so-called {\em representation theorem for Stieltjes transforms}. In some cases we use Lemma \ref{LemStrep1} and in other cases Lemma \ref{LemStrep2}. The third lemma in this subsection, that is, Lemma \ref{LemStinv} is the so-called {\em inversion theorem for Stieltjes transforms} or {\em Perron-Stieltjes inversion formula} and it is also a key ingredient in our proofs. For a  proof of Lemma \ref{LemStrep1} we refer to \cite[p. 235]{widder}, while for a proof of Lemma \ref{LemStrep2} we refer to \cite[p. 210]{widder} and Lemma \ref{LemStinv} can be found in \cite{stone}. Note also that Lemma \ref{LemStrep1} can be found in \cite[Lemma 2.1]{ismail00}, Lemma \ref{LemStrep2} is \cite[Theorem 1.2]{kelker} and Lemma \ref{LemStinv} can be found in \cite[Lemma 2.2]{ismail00}.

\begin{lemma}\label{LemStrep1}{\em (Representation theorem for Stieltjes transforms)}
A complex function $F(z)$ has the Stieltjes transform representation
\begin{equation}\label{repSt1}
F(z)=\int_{0}^{\infty}\frac{d\mu(t)}{z+t},\quad \mbox{with}\quad \int_{0}^{\infty}|d\mu(t)|<\infty
\end{equation}
if and only if the following conditions hold true
\begin{enumerate}
\item[\bf a.] $F(z)$ is analytic for $|\arg z|<\pi,$
\item[\bf b.] $F(z)=o(1)$ as $|z|\to \infty$ and $F(z)=o(|z|^{-1})$ as $|z|\to 0,$ uniformly in every sector $|\arg z|\leq\pi-\varepsilon$ for $\varepsilon>0$.
\end{enumerate}
\end{lemma}

\begin{lemma}\label{LemStrep2}{\em (Representation theorem for Stieltjes transforms)}
If
\begin{enumerate}
\item[\bf a.] $F(z)$ is analytic for $|\arg z|<{\pi}/{\theta}$ for some $\theta$ such that $0<\theta<1$,
\item[\bf b.] $F(z)=o(1)$ as $z\to\infty$ and $F(z)=o(|z|^{-1})$ as $z\to0,$ uniformly in every sector $|\arg z|\leq {\pi}/{\eta}$ with $\theta<\eta<1,$
\end{enumerate}
then the following Stieltjes transform representation is valid for all $x>0$
\begin{equation}\label{repSt2}
F(x)=\frac{1}{\pi}\int_{0}^{\infty}\frac{dt}{x+t}\frac{1}{2\pi \mathrm{i}}\int_{\mathcal{C}}\frac{ze^{\frac{z}{2}}F(te^{z})}{z^2+\pi^2}dz,
\end{equation}
where $\mathcal{C}$ is a rectifiable closed curve going around $[-\mathrm{i}\pi, \mathrm{i}\pi]$ in the positive direction and lying in the strip $\left|\imag z \right|<{\pi}/{\theta}$.
\end{lemma}

\begin{lemma}\label{LemStinv}{\em (Inversion theorem for Stieltjes transforms)}
If $F$ has the representation \eqref{repSt1}, then
\begin{equation}\label{inv1}
\mu(t_2)-\mu(t_1)=\lim_{\eta \to 0^{+}}\dfrac{1}{2\pi \mathrm{i}}\int_{t_1}^{t_2}\left[F(-t-\mathrm{i}\eta)-F(-t+\mathrm{i}\eta)\right]dt,
\end{equation}
where $\mu(t)$ is normalized by $\mu(0)=\mu(0^{+})=0$ and $2\mu(t)=\mu(t^{+})+\mu(t^{-})$ for $t>0.$
\end{lemma}

\section{\bf Distributions whose probability density function involves modified Bessel functions}

In this section our aim is to consider some known distributions whose probability density function contains the modified Bessel function of the first or second kind and to study whether these distributions belong to the class of infinitely divisible distributions or to one of its subclasses such as the self-decomposable distributions, generalized gamma convolutions and hyperbolically completely monotone densities.

\subsection{The McKay distribution of type I} First we consider the well-known McKay distribution of type I which involves the modified Bessel function of the first kind $I_{\mu}$. Its probability density function is given by
\begin{equation}\label{mcnoltypdf}
\varphi_{\mu,a,b}(x)=\frac{\sqrt{\pi}(b^2-a^2)^{\mu+\frac{1}{2}}}{(2a)^{\mu}\Gamma\left(\mu+\frac{1}{2}\right)}x^{\mu}e^{-bx}I_{\mu}(ax),
\end{equation}
where $b>a>0,$ $\mu>-\frac{1}{2}$ and the support of the distribution is $(0,\infty).$ In view of the asymptotic relation
$$I_{\mu}(x)\sim \frac{x^{\mu}}{2^{\mu}\Gamma(\mu+1)}$$
as $x\to 0$ and by using the Legendre duplication formula for the Euler gamma function
$$\Gamma(2x)=\frac{1}{\sqrt{\pi}}2^{2x-1}\Gamma(x)\Gamma\left(x+\frac{1}{2}\right)$$
we obtain that
$$\lim_{a\to0}\varphi_{\mu,a,b}(x)=\frac{b^{2\mu+1}}{\Gamma(2\mu+1)}e^{-bx}x^{2\mu},$$
which shows that the McKay distribution of type I for $a\to0$ reduces to a gamma distribution with shape parameter $2\mu+1$ and inverse scale parameter (rate parameter) $b.$ Since the gamma distribution
belongs to the class of infinitely divisible distributions, self-decomposable distributions and generalized gamma convolutions, it is natural to ask whether the McKay distribution of type I belongs to the above mentioned classes of infinitely divisible distributions. It is known (see \cite[p. 580]{bariczed}) that $x\mapsto e^{-x}x^{\mu}I_{\mu}(x)$ is completely monotonic on $(0,\infty)$ for all
$\mu\in\left[-\frac{1}{2},0\right],$ which in turn implies that $x\mapsto e^{-ax}(ax)^{\mu}I_{\mu}(ax)$ is completely monotonic on $(0,\infty)$ for all
$\mu\in\left[-\frac{1}{2},0\right]$ and $a>0.$ On the other hand, the function $x\mapsto e^{(a-b)x}$ is also completely monotonic on $(0,\infty)$ for all $b>a,$ and since the product of two completely monotonic functions is also completely monotonic, we arrive at the conclusion that {\em the probability density function $\varphi_{\mu,a,b}$ is completely monotonic on $(0,\infty)$ for all
$\mu\in\left(-\frac{1}{2},0\right]$ and $b>a>0$ and according to the Goldie-Steutel law the McKay distribution of type I is infinitely divisible for all $\mu\in\left(-\frac{1}{2},0\right]$ and $b>a>0.$} We note that with a more sophisticated analysis it is possible to show that the McKay distribution of type I is infinitely divisible for all $\mu>-\frac{1}{2}$ and $b>a>0.$

\begin{theorem}\label{th1}
If $\mu>-\frac{1}{2}$ and $b>a>0,$ then the McKay distribution, with probability density function defined by \eqref{mcnoltypdf},
belongs to the class of infinitely divisible distributions, self-decomposable distributions and generalized gamma convolutions.
\end{theorem}

\subsection{Another McKay type distribution} Now, we consider another distribution which is very similar to the McKay distribution of type I and involves also the modified Bessel function of the first kind $I_{\mu}$. Its probability density function is given by
\begin{equation}\label{newpdf}
\psi_{\mu,a,b}(x)=\frac{\sqrt{\pi}(b^2-a^2)^{\mu+\frac{3}{2}}}{2b(2a)^{\mu}\Gamma\left(\mu+\frac{3}{2}\right)}x^{\mu+1}e^{-bx}I_{\mu}(ax),
\end{equation}
where $b>a>0,$ $\mu>-1$ and the support of the distribution is $[0,\infty).$ The next result shows that this modified Bessel distribution belongs also to the class of infinitely divisible distributions.

\begin{theorem}\label{th2}
If $\mu>-1$ and $b>a>0,$ then the above distribution, with probability density function defined by \eqref{newpdf},
belongs to the class of infinitely divisible distributions.
\end{theorem}

\subsection{Generalization of the McKay distribution of type I} Next, we are going to study another McKay type distribution, which is in fact a generalization of the above McKay type distributions.
The support of this distribution is also $(0,\infty)$ and its probability density function is given by
\begin{equation}\label{pdfgener}
\varphi_{\mu,\nu,a,b}(x)=\frac{1}{c_{\mu,\nu,a,b}}x^{\nu-1}e^{-bx}I_{\mu}(ax),
\end{equation}
where $\mu+1>0,$ $\mu+\nu>0,$ $b>a>0,$
$$c_{\mu,\nu,a,b}=\frac{\left(a/2\right)^{\mu}}{b^{\mu+\nu}}\frac{\Gamma(\mu+\nu)}{\Gamma(\mu+1)}\cdot{}_2F_1\left(\frac{\mu+\nu}{2},\frac{\mu+\nu+1}{2},\mu+1,\frac{a^2}{b^2}\right)$$
and ${}_2F_1(a,b,c,x)$ stands for the Gaussian hypergeometric function. Observe that
$$c_{\mu,\mu+1,a,b}=\frac{(2a)^{\mu}\Gamma\left(\mu+\frac{1}{2}\right)}{\sqrt{\pi}b^{2\mu+1}}{}_1F_0\left(\mu+\frac{1}{2},\frac{a^2}{b^2}\right)=
\frac{(2a)^{\mu}\Gamma\left(\mu+\frac{1}{2}\right)}{\sqrt{\pi}\left(b^2-a^2\right)^{\mu+\frac{1}{2}}}$$
and
$$c_{\mu,\mu+2,a,b}=\frac{2(2a)^{\mu}\Gamma\left(\mu+\frac{3}{2}\right)}{\sqrt{\pi}b^{2\mu+2}}{}_1F_0\left(\mu+\frac{3}{2},\frac{a^2}{b^2}\right)=
\frac{(2b)(2a)^{\mu}\Gamma\left(\mu+\frac{3}{2}\right)}{\sqrt{\pi}\left(b^2-a^2\right)^{\mu+\frac{3}{2}}},$$
and consequently for all $x>0,$ $\mu>-\frac{1}{2}$ and $b>a>0$ we have that
$$\varphi_{\mu,\mu+1,a,b}(x)=\varphi_{\mu,a,b}(x)=\frac{\sqrt{\pi}(b^2-a^2)^{\mu+\frac{1}{2}}}{(2a)^{\mu}\Gamma\left(\mu+\frac{1}{2}\right)}x^{\mu}e^{-bx}I_{\mu}(ax)$$
and for all $x>0,$ $\mu>-1$ and $b>a>0$ we arrive at
$$\varphi_{\mu,\mu+2,a,b}(x)=\psi_{\mu,a,b}(x)=\frac{\sqrt{\pi}(b^2-a^2)^{\mu+\frac{3}{2}}}{2b(2a)^{\mu}\Gamma\left(\mu+\frac{3}{2}\right)}x^{\mu+1}e^{-bx}I_{\mu}(ax).$$

Now, recall that (see \cite[p. 578]{bariczed}) the function $x\mapsto e^{-x}x^{-\mu}I_{\mu}(x)$ is completely monotonic on $(0,\infty)$ for all
$\mu\geq-\frac{1}{2},$ which in turn implies that $x\mapsto e^{-ax}(ax)^{-\mu}I_{\mu}(ax)$ is completely monotonic on $(0,\infty)$ for all
$\mu\geq-\frac{1}{2}$ and $a>0.$ On the other hand, the function $x\mapsto e^{(a-b)x}$ is also completely monotonic on $(0,\infty)$ for all $b>a,$ and since the product of two completely monotonic functions is also completely monotonic, we arrive at the conclusion that the probability density function $\varphi_{\mu,\nu,a,b}$ is completely monotonic on $(0,\infty)$ for all
$\mu\geq-\frac{1}{2},$ $\mu+\nu\leq1$ and $b>a>0$ and according to the Goldie-Steutel law {\em the generalization of the McKay distribution of type I is infinitely divisible for all $\mu\geq-\frac{1}{2},$ $\mu+\nu\leq1$ and $b>a>0.$} The next theorem complements this result and shows that if $\mu$ and $\nu$ are in the domain bounded by the straight lines $\mu+1=\nu$ and $\mu+\nu=0$ in the $(\mu,\nu)$ right-half plane, then the generalization of the McKay distribution of type I belongs also to the class of infinitely divisible distributions. It would of interest to find the largest domain of $\mu$ and $\nu$ such that the above generalization of the McKay distribution is infinitely divisible.

\begin{theorem}\label{th3}
If $\mu+\nu>0,$ $\mu+1\geq\nu$ and $b>a>0,$ then the above distribution, with probability density function defined by \eqref{pdfgener},
belongs to the class of infinitely divisible distributions.
\end{theorem}

\subsection{One more McKay type distribution} Motivated by the result in the previous subsection, now we are going to study another McKay type distribution, which is in fact somehow similar to the above generalization of the McKay distribution of type I.
The support of this distribution is also $(0,\infty)$ and its probability density function is given by
\begin{equation}\label{pdfxi}
\xi_{\mu,a,b}(x)=\frac{1}{c_{\mu,a,b}}x^{2\mu}e^{-bx}\left[I_{\mu}(ax)\right]^2,
\end{equation}
where $\mu>-\frac{1}{4},$ $b>2a>0$ and
$$c_{\mu,a,b}=\frac{2^{4\mu}a^{2\mu}}{\pi b^{4\mu+1}}\frac{\Gamma\left(\mu+\frac{1}{2}\right)\Gamma\left(2\mu+\frac{1}{2}\right)}{\Gamma(\mu+1)}\cdot{}_2F_1\left(\mu+\frac{1}{2},2\mu+\frac{1}{2},\mu+1,\frac{4a^2}{b^2}\right).$$

Recall that (see \cite[p. 580]{bariczed}) the function $x\mapsto e^{-x}x^{\mu}I_{\mu}(x)$ is completely monotonic on $(0,\infty)$ for all
$\mu\in\left[-\frac{1}{2},0\right],$ which in turn implies that $x\mapsto e^{-ax}(ax)^{\mu}I_{\mu}(ax)$ and $x\mapsto e^{-2ax}(ax)^{2\mu}\left[I_{\mu}(ax)\right]^2$ are completely monotonic on $(0,\infty)$ for all
$\mu\in\left[-\frac{1}{2},0\right]$ and $a>0.$ On the other hand, the function $x\mapsto e^{(2a-b)x}$ is also completely monotonic on $(0,\infty)$ for all $b>2a,$ and since the product of two completely monotonic functions is also completely monotonic, we arrive at the conclusion that the probability density function $\xi_{\mu,a,b}$ is completely monotonic on $(0,\infty)$ for all
$\mu\in\left(-\frac{1}{4},0\right]$ and $b>2a>0$ and according to the Goldie-Steutel law {\em the above McKay type distribution is infinitely divisible for all $\mu\in\left(-\frac{1}{4},0\right]$ and $b>2a>0.$} We note that with a more sophisticated analysis it is possible to show that the above McKay type distribution is in fact infinitely divisible for all $\mu\in\left(-\frac{1}{4},\frac{1}{2}\right]$ and $b>2a>0.$

\begin{theorem}\label{th4}
If $\mu\in\left(-\frac{1}{4},\frac{1}{2}\right]$ and $b>2a>0,$ then the above distribution, with probability density function defined by \eqref{pdfxi},
belongs to the class of infinitely divisible distributions.
\end{theorem}

\subsection{The $K$-distribution or gamma-gamma distribution} Suppose that a random variable $X$ has gamma distribution with mean $\sigma$ and shape parameter $\alpha,$ with $\sigma$ being treated as a random variable having another gamma distribution, this time with mean $\mu$ and shape parameter $\beta.$ The result is that $X$ has the following probability density function (see \cite{jakeman})
\begin{equation}\label{pdfK}\omega_{\alpha,\beta,\mu}(x)=\frac{2}{\Gamma(\alpha)\Gamma(\beta)}\left(\frac{\alpha\beta}{\mu}\right)^{\frac{\alpha+\beta}{2}}x^{\frac{\alpha+\beta}{2}-1}K_{\alpha-\beta}\left(2\sqrt{\frac{\alpha\beta x}{\mu}}\right),\end{equation}
where $\alpha,\beta,\mu>0,$ and the support of the distribution is $(0,\infty).$  The $K$-distribution, which was used by Jakeman and Pusey (see \cite{jakeman}) to model microwave sea echo, is a compound (or mixture) probability distribution and it is also a product distribution: it is in fact the distribution of the product of two independent random variables, one having a gamma distribution with mean $1$ and shape parameter $\alpha>0$, the second having a gamma distribution with mean $\mu>0$ and shape parameter $\beta>0.$ It is well-known that gamma densities are hyperbolically completely monotone and the product of two hyperbolically completely monotone functions is also hyperbolically completely monotone (see \cite{bondesson} or \cite[Proposition 4]{bonde}). This in turn implies that the $K$-distribution belongs to the class of hyperbolically completely monotone distributions. Although the next theorem is obvious in view of the theory of hyperbolically completely monotone functions, we will give an alternative proof for this theorem by applying special functions technique.

\begin{theorem}\label{thK}
If $\alpha,\beta,\mu>0,$ then the $K$-distribution, with probability density function defined by \eqref{pdfK},
belongs to the class of infinitely divisible distributions, self-decomposable distributions, generalized gamma convolutions and hyperbolically completely monotone distributions.
\end{theorem}

It is worth to mention here that the $K$-distribution is also well-known in engineering sciences and it is called the gamma-gamma distribution. More precisely, the gamma-gamma
distribution is produced from the product of two independent gamma random variables and has been widely used in a variety of applications, for example in modeling various types of land and sea
radar clutters, in modeling the effects of the combined fading and shadowing phenomena, encountered in the mobile communications channels. Moreover, the gamma-gamma distribution is also applied in optical wireless systems, where transmission of optical signals through the atmosphere is involved. For more details we refer to the papers \cite{karag,karag2} and to the references therein. We also would like to point out that in \cite[Theorem 5]{baricz2} the authors proved that the probability density function of the gamma-gamma distribution $x\mapsto \omega_{\alpha,\beta,\mu}(x)$ is {\em geometrically concave} on $(0,\infty)$ for all $\alpha,\beta,\mu>0,$ that is the function $x\mapsto x\omega_{\alpha,\beta,\mu}'(x)/\omega_{\alpha,\beta,\mu}(x)$ is decreasing on $(0,\infty)$ for all $\alpha,\beta,\mu>0.$ It can be shown easily (see for example \cite[p. 102]{bondesson}) that this is equivalent to the fact that the probability density function of the $K$-distribution or gamma-gamma distribution is hyperbolically monotone (of order 1), that is, the function $w\mapsto \omega_{\alpha,\beta,\mu}(uv)\omega_{\alpha,\beta,\mu}(u/v)$ is decreasing on $(0,\infty)$ for all $u,v,\alpha,\beta,\mu>0,$ where $w=v+1/v.$ Clearly Theorem \ref{thK} states much more than this, the function $w\mapsto \omega_{\alpha,\beta,\mu}(uv)\omega_{\alpha,\beta,\mu}(u/v)$ is not only decreasing on $(0,\infty)$ for all $u,v,\alpha,\beta,\mu>0,$ it is even completely monotonic.

\subsection{The generalized inverse Gaussian distribution} The generalized inverse Gaussian distribution is a three-parameter family of continuous probability distributions with probability density function
\begin{equation}\label{newgigd}\pi_{\mu,a,b}(x)=\frac{(a/b)^{\mu/2}}{2K_{\mu}(\sqrt{ab})}x^{\mu-1}e^{-\frac{1}{2}(ax+b/x)}\end{equation}
and support $(0,\infty),$ where $a,b>0$ and $\mu$ is a real parameter. Note that Barndorff-Nielsen and Halgreen \cite{nielsen} have proved that the generalized Gaussian distribution is infinitely divisible. Moreover, Barndorff-Nielsen et al. \cite{nielsen2} have shown that the generalized inverse Gaussian distribution is a first hitting time for certain time-homogeneous diffusion processes provided the parameter $\mu$ is negative, and in this case the infinite divisibility of this distribution then follows from the general central limit theorem. In addition, Halgreen \cite{halgreen} and Bondesson \cite{bondes} have shown that the generalized inverse Gaussian distribution is a generalized gamma convolution, and according to Bondesson \cite[p. 74]{bondesson} we also know that the generalized inverse Gaussian distribution belongs to the class of hyperbolically completely monotone densities and hence to the class of generalized gamma convolutions and self-decomposable distributions. In the proof of the next result we would like to show there is a special function approach to show that the generalized inverse Gaussian distribution is a generalized gamma convolution.

\begin{theorem}\label{Thnewgigd}
If $\mu\in\mathbb{R}$ and $a,b>0,$ then the generalized inverse Gaussian distribution, with probability density function defined by \eqref{newgigd},
belongs to the class of generalized gamma convolutions and hence to the class of self-decomposable distributions and infinitely divisible distributions.
\end{theorem}

\section{\bf Distributions whose Laplace transform involves modified Bessel functions}

At the end of his book Lenart Bondesson \cite{bondesson} wrote the followings: ''Since the class of infinitely divisible distributions is extremely large, it is not a very interesting
class. In fact, as the research during the last two decades has shown, infinitely divisibility seems to be more a rule than an exception. This is not surprising if one considers that the class of
(univariate) distributions which are infinitely divisible with respect to the maximum operation contains all distributions. On the other hand, to investigate whether or not infinitely divisibility holds for a
particular distribution may lead to a deep insight into the structure of that and other distributions and also to a lot of by-products. (Cf. Riemann's hypothesis in
mathematics, for example.) This work would certainly not have been written had not Steutel (1973) asked whether the lognormal distribution is infinitely divisible and had not Thorin
(1977) attacked and solved this problem.''

In this section our aim is to consider some new lifetime distributions whose Laplace transform contains the modified Bessel function of the first and/or second kind and to study whether these distributions belong to the class of infinitely divisible distributions or to one of its subclasses such as the self-decomposable distributions and generalized gamma convolutions. We also think that it is interesting to find out whether or not particular distributions are infinitely divisible and the first part of this section was in fact motivated by the conjecture of Ismail and Miller \cite[p. 234]{miller} which is the following:
{\em if $\nu>\mu\geq0$ and $b>a>0,$ then the function
$$x\mapsto \left(\frac{b}{a}\right)^{\mu-\nu}\frac{I_{\mu}(a\sqrt{x})I_{\nu}(b\sqrt{x})}{I_{\mu}(b\sqrt{x})I_{\nu}(a\sqrt{x})}$$
is the Laplace transform of an infinitely divisible probability distribution.} This conjecture is still open and Theorems \ref{theolap1} and \ref{theolap2} show that it is possible to generate with some slight modifications from the above quotient of modified Bessel functions of the first kind a Laplace transform of an infinitely divisible probability distribution.

\subsection{Quotients of modified Bessel functions of the first kind} Due to Ismail and Kelker \cite[Theorem 1.10]{kelker}, we know that the function
\begin{equation}\label{iskellap1}x \mapsto \rho_{\mu,a}(x)=\frac{1}{2^{\mu}\Gamma(\mu+1)}\frac{(a\sqrt{x})^{\mu}}{I_{\mu}(a\sqrt{x})}\end{equation} is a Laplace transform of an infinitely divisible distribution on $[0, \infty)$, when $\mu>0$ and $a>0$. It is worth to mention here that the above mentioned distribution with the Laplace transform in \eqref{iskellap1} belongs also to the class of self-decomposable distribution and the proof follows naturally according to Lemma \ref{lemma2}. Moreover, the next result shows that in fact the above function is a Laplace transform of a generalized gamma convolution.

\begin{theorem}\label{thiskellap1}
If $\mu>-1$ and $a>0$, then $x \mapsto \rho_{\mu,a}(x)$ is the Laplace transform of a generalized gamma convolution and therefore it is also the Laplace transform of a self-decomposable distribution.
\end{theorem}

\begin{theorem}\label{theolap1}
If $\mu>-1$, $\nu>\sigma>-1$ and $b>a>0,$ then
\begin{equation}\label{omegamunu1}
\Omega_{\mu,\nu,\sigma,a,b}(x)=\left(\frac{b}{a}\right)^{\mu-\nu}\frac{I_{\mu}(a\sqrt{x})I_{\nu}(b\sqrt{x})}{I_{\mu}(b\sqrt{x})I_{\nu}(a\sqrt{x})}\cdot \rho_{\sigma,b}(x)
\end{equation}
is the Laplace transform of an infinitely divisible distribution with support $[0,\infty).$
\end{theorem}

\begin{theorem}\label{theolap2}
If $\mu>-1,$ $\nu>\frac{1}{2}$ and $b>a>0,$ then the function
\begin{equation}\label{omegamunu2}
\Omega_{\mu,\nu,a,b}(x)=\left(\frac{b}{a}\right)^{\mu-\nu}\dfrac{I_{\mu}(a\sqrt{x})I_{\nu}(b\sqrt{x})}{I_{\mu}(b\sqrt{x})I_{\nu}(a\sqrt{x})}\cdot e^{-b\sqrt{x}}
\end{equation}
is the Laplace transform of an infinitely divisible distribution with support $[0,\infty).$
\end{theorem}

\subsection{Stieltjes transform representations and infinite divisibility}

Now, we are going to present a series of results related to Stieltjes transforms of modified Bessel functions of the first and second kinds, their products and quotients and our aim is to show that these new Stieltjes transform representations are strongly related to some new infinitely divisible modified Bessel distributions. The proofs of this series of results are based on representation and inversion theorems for Stieltjes transforms.
Some of the proofs related to infinitely divisibility were inspired from the papers of Mourad Ismail, however in each case we will also show the integral representation of the probability density functions in question.

The first result related to Stieltjes transforms is the following one.

\begin{theorem}\label{thIfirst}
If $a>0$, $\mu>-\frac{1}{2}$ and $|\arg z|<\pi,$ then the following Stieltjes transform representation is valid
$$e^{-a\sqrt{z}}z^{-\frac{\mu}{2}}I_{\mu}(a\sqrt{z})=\dfrac{1}{\pi}\int_{0}^{\infty}\frac{t^{-\frac{\mu}{2}}}{z+t}J_{\mu}(a\sqrt{t})\sin(a\sqrt{t})dt,$$
which can be rewritten as the two-fold Laplace transform
$$e^{-a\sqrt{z}}z^{-\frac{\mu}{2}}I_{\mu}(a\sqrt{z})=\dfrac{1}{\pi}\int_{0}^{\infty}e^{-zs}\left[\int_{0}^{\infty}e^{-st}t^{-\frac{\mu}{2}}J_{\mu}(a\sqrt{t})\sin(a\sqrt{t})dt\right]ds.$$
\end{theorem}

It is important to mention here that due to Ismail \cite[Theorem 2]{ismail} we know that if $\mu>\frac{1}{2},$ then the function $x\mapsto 2^{\mu}\Gamma(\mu+1)x^{-\mu/2}I_{\mu}(\sqrt{x})e^{-\sqrt{x}}$ is the Laplace transform of an infinitely divisible distribution, which is not a generalized gamma convolution. This in turn implies that the function $$x\mapsto \frac{2^{\mu}\Gamma(\mu+1)}{a^{\mu}}x^{-\frac{\mu}{2}}e^{-a\sqrt{x}}I_{\mu}(a\sqrt{x})$$
is also the Laplace transform of an infinitely divisible distribution, which is not a generalized gamma convolution. In view of Theorem \ref{thIfirst} we can clearly state that the function
$$s\mapsto \frac{1}{\pi}\frac{2^{\mu}\Gamma(\mu+1)}{a^{\mu}}\int_{0}^{\infty}e^{-st}t^{-\frac{\mu}{2}}J_{\mu}(a\sqrt{t})\sin(a\sqrt{t})dt$$ is the corresponding probability density function for the above distribution in the case when $\mu>\frac{1}{2}.$ This complements the above mentioned result of Ismail \cite[Theorem 2]{ismail}.

It is worth mentioning that the following formula $$z^{-\frac{\mu}{2}}e^{-a\sqrt{z}}I_{\mu}(b\sqrt{z})=\frac{1}{\pi}\int_{0}^{\infty}\frac{t^{-\frac{\mu}{2}}}{z+t}J_{\mu}(b\sqrt{t})\sin(a\sqrt{t})dt,$$ where $\real \mu>-\frac{1}{2},$ $a>b>0$ is available in \cite[eq. (19), p. 226]{erdelyi2}. Theorem \ref{thIfirst} can be extended to the case when $a>0,$ $\real \mu>-\frac{1}{2}$ and $|\arg z|<\pi,$ which clearly complements the above formula for the case $a=b$.

The next result is about a similar result for the product of modified Bessel functions.

\begin{theorem}\label{thmprod00}
Let $a\leq b$, $\mu, \nu>-1,$ $\nu-\mu<1$ and $\left|\arg z\right| < \pi$. Then the following Stieltjes transform representation is valid
\begin{equation}\label{eqproddifpar}
z^{\frac{\nu-\mu}{2}}I_{\mu}(a\sqrt{z})K_{\nu}(b\sqrt{z})=\frac{1}{2}\int_{0}^{\infty} \frac{t^{\frac{\nu-\mu}{2}}}{z+t}J_{\mu}(a\sqrt{t})J_{\nu}(b\sqrt{t})dt,
\end{equation}
which in fact can be rewritten as follows
$$z^{\frac{\nu-\mu}{2}}I_{\mu}(a\sqrt{z})K_{\nu}(b\sqrt{z})=\frac{1}{2}\int_{0}^{\infty}e^{-zs}\left[\int_{0}^{\infty}e^{-st}
t^{\frac{\nu-\mu}{2}}J_{\mu}(a\sqrt{t})J_{\nu}(b\sqrt{t})dt\right]ds.$$
\end{theorem}

It is worth mentioning that in Theorem \ref{thmprod00} we can also assume that the order of the modified Bessel functions $\mu$ and $\nu$ are complex numbers and Theorem \ref{thmprod00} is also valid with the following conditions: $a\leq b,$ $\real \mu>-1,$ $\real \nu>-1,$ $\real (\nu-\mu)<1$ and $\left|\arg z\right| < \pi.$ Moreover, during the process of our work we found out that the integral representation in \eqref{eqproddifpar} is available in the literature (see \cite[p. 96]{MOS66} and \cite[(24), p. 227]{erdelyi2}) in the form
\begin{equation}\label{eqhankel}z^{\nu-\mu}I_{\mu}(az)K_{\nu}(bz)=\int_{0}^{\infty}\frac{t^{\nu-\mu+1}}{z^2+t^2}J_{\mu}(at)J_{\nu}(bt)dt\end{equation}
with the better conditions  $a\leq b,$ $\real \mu>-1,$ $\real \nu>-1,$ $\real (\nu-\mu)<2$ and $\real z>0.$ By replacing $z$ by $\sqrt{z}$ and introducing the suitable transformation, the above formula becomes \eqref{eqproddifpar} with the conditions $a\leq b,$ $\real \mu>-1,$ $\real \nu>-1,$ $\real (\nu-\mu)<2$ and $\left|\arg z\right|<\pi$. This means that the function $z\mapsto z^{\frac{\nu-\mu}{2}}I_{\mu}(a\sqrt{z})K_{\nu}(b\sqrt{z})$ has the corresponding Stieltjes transform representation with above mentioned conditions. But we observe that for the case $a=b$, using the asymptotic relations \cite[eq. 10.30.4]{nist} and \cite[eq. 10.25.3]{nist}, we have that as $z \to\infty$ $$z^{\frac{\nu-\mu}{2}}I_{\mu}(a\sqrt{z})K_{\nu}(b\sqrt{z})\sim z^{\frac{\nu-\mu}{2}}\cdot\left(\frac{e^{a\sqrt{z}}}{\sqrt{2\pi a\sqrt{z}}}\right)\cdot\left(\sqrt{\frac{2\pi}{a\sqrt{z}}}e^{-a\sqrt{z}}\right)=\frac{1}{a}z^{\frac{\nu-\mu-1}{2}},$$ which does not tends to zero as $z \to \infty$ in the case when $\real (\nu-\mu)\geq1.$ This violates one of the conditions in Lemma \ref{LemStrep1} and thus, the condition $\real (\nu-\mu)<2$ mentioned in \cite[p. 96]{MOS66} and \cite[eq. (24), p. 227]{erdelyi2} is not correct and we believe that the condition $\real (\nu-\mu)<1$ that we mentioned in Theorem \ref{thmprod00} is accurate.

We also note that for $\mu=\nu=n$ integer the equation \eqref{eqhankel} appears in \cite{mitra} with a reference to \cite{hankel}, while a more general form of \eqref{eqhankel} for integer orders has been considered by von Hermann Hankel \cite{hankel}. Namely, the extension of \eqref{eqhankel} in the case when $k\in\mathbb{N}$ and $n+m+q$ is an even integer appears in the form
$$-\frac{1}{2\Gamma(k+1)}\left(\frac{\partial}{\partial r^2}\right)^k r^qJ_n(ar)\left[Y_n(br)-\pi \mathrm{i}J_n(br)\right]=\int_0^{\infty}t^{q+1}\frac{J_n(at)J_m(bt)}{(t^2-r^2)^{k+1}}dt.$$

The following result provides an immediate application of \eqref{eqproddifpar}.
	
\begin{theorem}\label{thprod1}
If $\mu>-1$ and $\left|\arg z\right|<\pi,$ then the following representations hold true
\begin{equation}\label{eqprod1}
2I_{\mu}(\sqrt{z})K_{\mu}(\sqrt{z})=\int_{0}^{\infty}\frac{J_{\mu}^2(\sqrt{t})}{z+t}dt=\int_{0}^{\infty} e^{-zs} \left[ \int_{0}^{\infty} e^{-st}J_{\mu}^2(\sqrt{t})dt\right]ds.
\end{equation}
and consequently for $\mu>0$ the function $x\mapsto 2\mu I_{\mu}(\sqrt{x})K_{\mu}(\sqrt{x})$ is the Laplace transform of an infinitely divisible probability distribution with support $[0,\infty),$ and its corresponding probability density function is as follows $$\varsigma_{\mu}(x)=\mu\int_{0}^{\infty} e^{-tx}J_{\mu}^2(\sqrt{t})dt=\frac{\mu}{x}e^{-\frac{1}{2x}}I_{\mu}\left(\frac{1}{2x}\right).$$
\end{theorem}

It is worth to mention here that in fact a Stieltjes transform, at least formally, can be viewed as a two-fold Laplace transform, namely
$$\int_{0}^{\infty}\frac{d\mu(t)}{z+t}=\int_{0}^{\infty}e^{-zs}\int_{0}^{\infty}e^{-st}d\mu(t)ds,$$
and this relation immediately implies the second equation in Theorems \ref{thIfirst} and \ref{thmprod00} as well as the second equality in \eqref{eqprod1}. It is also interesting to mention that the formula \eqref{eqprod1} is already available in literature (see \cite[eq. (22), p.226]{erdelyi2}) with condition $\mu \in \mathbb{C}$. But for the case $\real \mu \leq-1$, $\mu \neq -1$ by using the asymptotic relations \cite[eq. 10.30.1]{nist}, \cite[eq. 10.30.2]{nist} and the formula \cite[eq. 10.27.3]{nist}, we observe that as $z \to 0$ $$I_{\mu}(\sqrt{z})K_{\mu}(\sqrt{z})=I_{\mu}(\sqrt{z})K_{-\mu}(\sqrt{z})\sim\dfrac{\Gamma(-\mu)}{2^{2\mu+1}\Gamma(\mu+1)}z^{\mu}.$$ It is clear that $\left|z\right|I_{\mu}(\sqrt{z})K_{\mu}(\sqrt{z})$ does not tends to zero as $z \to 0.$ This violates one of the conditions in Lemma \ref{LemStrep1} and we can obtain the same conclusion for the case $\mu=-1$ by using the additional formula \cite[eq. 10.27.1]{nist}. Thus, the condition $\mu \in \mathbb{C}$ mentioned in \cite[eq. (22), p.226]{erdelyi2} is incorrect and the condition $\real\mu-1$ or more simply $\mu>-1$ that we have assumed in Theorem \ref{thprod1} is accurate.

Moreover, the next corollary contains some immediate applications of the formula in \eqref{eqprod1} concerning the product of modified Bessel functions of the first and second kind. The first part of the next corollary is well-known and it was proved by using non-trivial arguments in 2007 by Penfold et al. \cite{penfold} for $\mu\geq0,$ in 2009 by Baricz \cite{baricz} for $\mu\geq-\frac{1}{2}$ and in 2021 by Segura \cite{segura} for $\mu\geq-1.$

\begin{corollary}\label{corprodcm}
The following assertions are valid.
\begin{enumerate}
\item[\bf a.] The function $x\mapsto I_{\mu}(x)K_{\mu}(x)$ is decreasing on $(0,\infty)$ for all $\mu>-1.$
\item[\bf b.] The function $x\mapsto I_{\mu}(\sqrt{x})K_{\mu}(\sqrt{x})$ is completely monotonic on $(0,\infty)$ for all $\mu>-1.$
\item[\bf c.] The function $x\mapsto xI_{\mu}(\sqrt{x})K_{\mu}(\sqrt{x})$ is a Bernstein function on $(0,\infty)$ for all $\mu>-1.$
\item[\bf d.] The function $x\mapsto \left[xI_{\mu}(\sqrt{x})K_{\mu}(\sqrt{x})\right]^{-1}$ is completely monotonic on $(0,\infty)$ for all $\mu\geq0.$
%\item[\bf e.] The function $x\mapsto I_{\mu}(\sqrt{x})K_{\mu}(\sqrt{x})$ is geometrically concave on $(0,\infty)$ for all $\mu>0.$
%\item[\bf f.] The function $x\mapsto I_{\mu}(x)K_{\mu}(x)$ is geometrically concave on $(0,\infty)$ for all $\mu>0.$
\item[\bf e.] If $\mu>0$ and $x>0,$ then the following inequality is valid
$$I_{\mu}(x)K_{\mu}(x)\leq \frac{\pi c_L^2}{\sqrt{3}x^{\frac{1}{3}}},$$
where $c_{L}=sup_{t\in\mathbb{R_{+}}}\sqrt[3]{t}J_{0}(t)\simeq 0.7857468704{\ldots}.$
\end{enumerate}
\end{corollary}

It is interesting to note that recently in \cite[Theorem 1]{BMPS16} the authors obtained a more general bound for $\mu\geq\nu$ and $x>0$ as follows
$$I_{\mu}(x)K_{\nu}(x)\leq\frac{2\pi^{\frac{3}{2}}c_{L}}{\sqrt{3}\Gamma\left(\frac{2}{3}\right)\Gamma\left(\frac{5}{6}\right)(2x)^{\frac{1}{3}}}.$$
We note that since $$c_L<\frac{2^{\frac{2}{3}}\sqrt{\pi}}{\Gamma\left(\frac{2}{3}\right)\Gamma\left(\frac{5}{6}\right)}\simeq 1.8407427466{\ldots},$$ the upper bound in part {\bf e} of the above corollary is clearly sharper than the above bound from \cite{BMPS16} in the case when $\mu=\nu.$

Another infinitely divisible distribution related to the product of $I_{\mu}(a\sqrt{x})K_{\nu}(b\sqrt{x})$ can be found in the next theorems.

\begin{theorem}\label{theprodIKexprepr2}
If $\mu,\nu>-1,$ $\nu-\mu<1$, $a,b>0$ and $|\arg z|<\pi,$ then the following Stieltjes transform representation holds true
$$z^{\frac{\nu-\mu}{2}}e^{-a\sqrt{z}}I_{\mu}(a\sqrt{z})K_{\nu}(b\sqrt{z})=\frac{1}{2}\int_{0}^{\infty}\frac{t^{\frac{\nu-\mu}{2}}J_{\mu}(a\sqrt{t})}{z+t}
\left[J_{\nu}(b\sqrt{t})\cos(a\sqrt{t})-Y_{\nu}(b\sqrt{t})\sin(a\sqrt{t})\right]dt,$$
which can be rewritten as $$\frac{I_{\mu}(a\sqrt{z})K_{\nu}(b\sqrt{z})}{z^{\frac{\mu-\nu}{2}}e^{a\sqrt{z}}}=\frac{1}{2}\int_{0}^{\infty}e^{-zs}\left[\int_{0}^{\infty}e^{-st}t^{\frac{\nu-\mu}{2}}J_{\mu}(a\sqrt{t})
\left[J_{\nu}(b\sqrt{t})\cos(a\sqrt{t})-Y_{\nu}(b\sqrt{t})\sin(a\sqrt{t})\right]dt \right]ds.$$
\end{theorem}

The next theorem shows that the function
$$x\mapsto \frac{2^{\mu-\nu+1}\Gamma(\mu+1)b^{\nu}}{a^{\mu}\Gamma(\nu)}e^{-a\sqrt{x}}x^{\frac{\nu-\mu}{2}}I_{\mu}(a\sqrt{x})K_{\nu}(b\sqrt{x})$$
is the Laplace transform of an infinitely divisible distribution on $(0, \infty)$ with the conditions $a,b,\nu>0$ and $\mu>\frac{1}{2}$. Based on Theorem \ref{theprodIKexprepr2} we can express the probability density function of the corresponding distribution as follows
$$s\mapsto \frac{2^{\mu-\nu}\Gamma(\mu+1)b^{\nu}}{a^{\mu}\Gamma(\nu)}\int_{0}^{\infty}e^{-st}t^{\frac{\nu-\mu}{2}}J_{\mu}(a\sqrt{t})\left[J_{\nu}(b\sqrt{t})\cos(a\sqrt{t})-Y_{\nu}(b\sqrt{t})\sin(a\sqrt{t})\right]dt$$
whenever $a,b,\nu>0$, $\mu>\frac{1}{2}$ and $\nu-\mu<1.$

\begin{theorem}\label{theprodIKexp}
If $a,b,\nu>0$ and $\mu>\frac{1}{2},$ then the function
$$\chi_{\mu,\nu,a,b}(x)=\frac{2^{\mu-\nu+1}\Gamma(\mu+1)b^{\nu}}{a^{\mu}\Gamma(\nu)}e^{-a\sqrt{x}}x^{\frac{\nu-\mu}{2}}I_{\mu}(a\sqrt{x})K_{\nu}(b\sqrt{x})$$
is the Laplace transform of an infinitely divisible distribution on $(0,\infty)$.
\end{theorem}

The next result is also related to an infinitely divisible distribution, but involves only modified Bessel functions of the second kind.

\begin{theorem}\label{thprodK1}
If $a,b,\mu,\nu>0,$ $\nu+\mu<1$ and $\left|\arg z\right|<\pi,$ then the following representation holds true
\begin{equation}\label{eqprodK1}
z^{\frac{\mu+\nu}{2}}K_{\mu}(a\sqrt{z})K_{\nu}(b\sqrt{z})=-\frac{\pi}{4}\int_{0}^{\infty}\frac{t^{\frac{\mu+\nu}{2}}}{z+t}
\left[J_{\mu}(a\sqrt{t})Y_{\nu}(b\sqrt{t})+J_{\nu}(b\sqrt{t})Y_{\mu}(a\sqrt{t})\right] dt,
\end{equation}
which can be rewritten as
$$z^{\frac{\mu+\nu}{2}}K_{\mu}(a\sqrt{z})K_{\nu}(b\sqrt{z})= \int_{0}^{\infty}e^{-zs}\left[-\frac{\pi}{4}\int_{0}^{\infty}e^{-st}t^{\frac{\mu+\nu}{2}}\left[J_{\mu}(a\sqrt{t})Y_{\nu}(b\sqrt{t})+J_{\nu}(b\sqrt{t})Y_{\mu}(a\sqrt{t})\right] dt\right] ds.$$
\end{theorem}	

The following result shows that the above product of modified Bessel functions generates naturally a generalized gamma convolution. On the other hand, if we let $a, b, \mu, \nu$ and $s$ be strictly positive real numbers, and $\nu+\mu<1$, then the following function $$s\mapsto \frac{-\pi a^{\mu}b^{\nu}}{2^{\mu+\nu}\Gamma(\mu)\Gamma(\nu)}\int_{0}^{\infty}e^{-st}t^{\frac{\mu+\nu}{2}}\left[J_{\mu}(a\sqrt{t})Y_{\nu}(b\sqrt{t})+J_{\nu}(b\sqrt{t})Y_{\mu}(a\sqrt{t})\right] dt$$ is in fact a probability density function and we clearly have that
$$\int_{0}^{\infty}\left[\int_{0}^{\infty}e^{-st}t^{\frac{\mu+\nu}{2}}\left[J_{\mu}(a\sqrt{t})Y_{\nu}(b\sqrt{t})+J_{\nu}(b\sqrt{t})Y_{\mu}(a\sqrt{t})\right] dt\right]ds=\dfrac{-2^{\mu+\nu}\Gamma(\mu)\Gamma(\nu)}{\pi a^{\mu}b^{\nu}}.$$ This follows from Theorem \ref{thprodK1} and the fact that function in \eqref{07prod} is the Laplace transform of an infinitely divisible distribution. Moreover, in the case when $a=b=1$ in view of the
integral representation \cite[eq. (65), p. 97]{erdelyi}
$$J_{\mu}(x)Y_{\nu}(x)+J_{\nu}(x)Y_{\mu}(x)=-\frac{4}{\pi}\int_0^{\infty}J_{\mu+\nu}(2x\cosh s)\cosh((\mu-\nu)s)ds$$
we obtain that \eqref{eqprodK1} reduces to
$$z^{\frac{\mu+\nu}{2}}K_{\mu}(\sqrt{z})K_{\nu}(\sqrt{z})=\int_{0}^{\infty}\left[\int_{0}^{\infty}\frac{t^{\frac{\mu+\nu}{2}}}{z+t}J_{\mu+\nu}(2\sqrt{t}\cosh s)\cosh((\mu-\nu)s)ds\right]dt,$$
where $\mu,\nu>0,$ $\nu+\mu<1$ and $\left|\arg z\right|<\pi,$ as before.

It is worth also to mention that in Theorem \ref{thprodK1} it is possible to assume that the orders $\mu$ and $\nu$ are complex numbers with the conditions such that $\real\mu>0$, $\real\nu>0$ and $\real(\nu+\mu)<1$ During the process of our work we found out that the integral representation in \eqref{eqprodK1} is available in the literature (see \cite[p. 96]{MOS66}) in the following forms
$$(-1)^{l+1}\frac{2}{\pi}z^{\mu+\nu+2l}K_{\nu}(az)K_{\mu}(bz)=\int_{0}^{\infty}\frac{t^{\mu+\nu+2l+1}}{z^2+t^2}
\left[J_{\nu}(at)Y_{\mu}(bt)+J_{\mu}(bt)Y_{\nu}(at)\right] dt,$$
where $l\in\{0,\pm1,\pm2,\ldots\},$ $\real(\nu+l)>-1,$ $\real(\mu+l)>-1,$ $l-1<\real(\mu+\nu+2l)<l,$ $\real z>0$ and
$$(-1)^{l+1}\frac{2}{\pi}z^{\mu+\nu+2l-1}K_{\nu}(az)K_{\mu}(bz)=\int_{0}^{\infty}\frac{t^{\mu+\nu+2l}}{z^2+t^2}
\left[J_{\nu}(at)J_{\mu}(bt)-Y_{\nu}(at)Y_{\mu}(bt)\right] dt,$$
where $l\in\{0,\pm1,\pm2,\ldots\},$ $\real(\nu+l)>-\frac{1}{2},$ $\real(\mu+l)>-\frac{1}{2},$ $l-\frac{1}{2}<\real(\mu+\nu+2l)<l,$ $\real z>0.$ Clearly our product representation \eqref{eqprodK1} corresponds to case when
$z^{\mu+\nu+2l}$ reduces to $z^{\mu+\nu}$ and in this case the conditions will be as follows $\real\nu>-1,$ $\real\mu>-1,$ $-1<\real(\nu+\mu)<0$ and $\real z>0,$ which complements our conditions $\real\mu>0,$ $\real\nu>0$, $\real(\nu+\mu)<1$ and $\left|\arg z\right|<\pi.$

\begin{theorem}\label{thinfdivprodK}
If $a, b,$ $\mu$, and $\nu$ are strictly positive real numbers, then the function
\begin{equation}\label{07prod}
\vartheta_{\mu,\nu,a,b}(x)=\frac{a^{\mu}b^{\nu}}{2^{\mu+\nu-2}\Gamma(\mu)\Gamma(\nu)}x^{\frac{\mu+\nu}{2}}K_{\mu}(a\sqrt{x})K_{\nu}(b\sqrt{x})
\end{equation}
is the Laplace transform of an infinitely divisible distribution on $(0, \infty),$ which is a self-decomposable distribution and a generalized gamma convolution.
\end{theorem}

It is worth mentioning that Ismail and Kelker \cite[Theorem 1.8]{kelker} proved that if $\mu>\nu>-1,$ then the function
$x\mapsto (\sqrt{x})^{\nu-\mu}K_{\nu}(\sqrt{x})/K_{\mu}(\sqrt{x})$ is the Laplace transform of an infinitely divisible distribution with support $(0,\infty)$
and consequently taking into account that $K_{\mu}(x)\sim 2^{\mu-1}x^{-\mu}\Gamma(\mu)$ as $\mu\to\infty$ for $x>0$ fixed, then the function $x\mapsto (\sqrt{x})^{\nu}K_{\nu}(\sqrt{x})/\left[2^{\nu-1}\Gamma(\nu)\right]$ is also the Laplace transform of an infinitely divisible distribution with support $(0,\infty)$ whenever $\nu>-1.$ However, after verifying the proof of \cite[Theorem 1.8]{kelker} we arrived to the conclusion that the above results are only true in the case when $\mu>\nu>0,$ and $\nu>0,$ respectively.

The next result is analogous to Theorem \ref{thprodK1}.

\begin{theorem}\label{thprodeqI}
If $a,b>0$, $\mu,\nu>-1$, $\mu+\nu>-1$ and $\left|\arg z\right|<\pi,$ then the following Stieltjes transform representation holds true
\begin{equation}\label{prodeqI}
e^{-(a+b)\sqrt{z}}z^{-\frac{\mu+\nu}{2}}I_{\mu}(a\sqrt{z})I_{\nu}(b\sqrt{z})=
\frac{1}{\pi}\int_{0}^{\infty}\dfrac{t^{-\frac{\mu+\nu}{2}}}{z+t}\left[J_{\mu}(a\sqrt{t})J_{\nu}(b\sqrt{t})\sin((a+b)\sqrt{t})\right]dt,
\end{equation}
which can be rewritten as
$$e^{-(a+b)\sqrt{z}}z^{-\frac{\mu+\nu}{2}}I_{\mu}(a\sqrt{z})I_{\nu}(b\sqrt{z})=
\frac{1}{\pi}\int_{0}^{\infty}e^{-zs}\left[\int_{0}^{\infty}e^{-st}t^{-\frac{\mu+\nu}{2}}\left[J_{\mu}(a\sqrt{t})J_{\nu}(b\sqrt{t})\sin((a+b)\sqrt{t})\right]dt\right]ds.$$
\end{theorem}

The product of two modified Bessel functions of the first kind also generates an infinitely divisible distribution, which however will not be a generalized gamma convolution.

\begin{theorem}\label{thprodeqIinfdiv}
If $a,b>0$ and $\mu,\nu>\frac{1}{2},$ then the function
$$\zeta_{\mu,\nu,a,b}(x)=\frac{2^{\mu+\nu}\Gamma(\mu+1)\Gamma(\nu+1)}{a^{\mu}b^{\nu}}e^{-(a+b)\sqrt{x}}x^{-\frac{\mu+\nu}{2}}I_{\mu}(a\sqrt{x})I_{\nu}(b\sqrt{x})$$
is a Laplace transform of an infinitely divisible distribution on $(0, \infty),$ which is not a generalized gamma convolution.
\end{theorem}

From Theorems \ref{thprodeqI} and \ref{thprodeqIinfdiv} for $a,b>0$ and $\mu,\nu>\frac{1}{2}$ we immediately conclude that  $$s\mapsto\frac{2^{\mu+\nu}\Gamma(\mu+1)\Gamma(\nu+1)}{a^{\mu}b^{\nu}\pi}\int_{0}^{\infty}e^{-st}t^{-\frac{\mu+\nu}{2}}
\left[J_{\mu}(a\sqrt{t})J_{\nu}(b\sqrt{t})\sin(a+b)\sqrt{t}\right]dt$$
is a probability density function and
$$\int_{0}^{\infty}\left[\int_{0}^{\infty}e^{-st}t^{-\frac{\mu+\nu}{2}}\left[J_{\mu}(a\sqrt{t})J_{\nu}(b\sqrt{t})\sin(a+b)\sqrt{t}\right]dt\right]ds=\frac{\pi a^{\mu}b^{\nu}}{2^{\mu+\nu}\Gamma(\mu+1)\Gamma(\nu+1)}.$$

The next result is the counterpart of Theorem \ref{thprodK1}.

\begin{theorem}\label{recprodKrepr}
If $a, b> 0,$ $\mu, \nu \in \mathbb{R},$  $\mu+\nu > 1$ and $\left|\arg z\right|<\pi,$ then the following representation holds true
\begin{equation}
\frac{e^{-(a+b)\sqrt{z}}}{z^{\frac{\mu+\nu}{2}}K_{\mu}(a\sqrt{z})K_{\nu}(b\sqrt{z})}
=\frac{4}{\pi^3}\int_{0}^{\infty}\frac{t^{-\frac{\mu+\nu}{2}}}{z+t}\Gamma_{\mu,\nu,a,b}(t)dt
\end{equation}
or equivalently
$$\frac{e^{-(a+b)\sqrt{z}}}{z^{\frac{\mu+\nu}{2}}K_{\mu}(a\sqrt{z})K_{\nu}(b\sqrt{z})}
=\frac{4}{\pi^3}\int_{0}^{\infty}e^{-zs}\left[\int_{0}^{\infty}{e^{-st}}{t^{-\frac{\mu+\nu}{2}}}
\Gamma_{\mu,\nu,a,b}(t) dt\right] ds,$$
where
$$\Gamma_{\mu,\nu,a,b}(t)=\frac{{}_1T_{\mu,\nu,a,b}(t)\cos((a+b)\sqrt{t})-{}_2T_{\mu,\nu,a,b}(t)\sin((a+b)\sqrt{t})}
{\left[J^2_{\mu}(a\sqrt{t})+Y^2_{\mu}(a\sqrt{t})\right]\left[J^2_{\nu}(b\sqrt{t})+Y^2_{\nu}(b\sqrt{t})\right]}$$
with
$${}_1T_{\mu,\nu,a,b}(t)=J_{\mu}(a\sqrt{t})Y_{\nu}(b\sqrt{t})+J_{\nu}(b\sqrt{t})Y_{\mu}(a\sqrt{t})$$
and
$${}_2T_{\mu,\nu,a,b}(t)=J_{\mu}(a\sqrt{t})J_{\nu}(b\sqrt{t})-Y_{\mu}(a\sqrt{t})Y_{\nu}(b\sqrt{t}).$$
\end{theorem}

The next theorem shows that the reciprocal of the product of two modified Bessel functions of the second kind also generates an infinitely divisible distribution.

\begin{theorem}\label{recprodKinfdiv}
If $a,b>0,$ $\mu,\nu>\frac{1}{2},$ then the function
$$\kappa_{\mu,\nu,a,b}(x)=\frac{2^{\mu+\nu-2}\Gamma(\mu)\Gamma(\nu)}{a^{\mu}b^{\nu}}\frac{e^{-(a+b)\sqrt{x}}}{x^{\frac{\mu+\nu}{2}}K_{\mu}(a\sqrt{x})K_{\nu}(b\sqrt{x})}$$ is the Laplace transform of an infinitely divisible distribution on $(0, \infty),$ which is a generalized gamma convolution.
\end{theorem}

It is worth to mention here that in view of Theorems \ref{recprodKrepr} and \ref{recprodKinfdiv} we conclude that the function
$$s\mapsto \frac{2^{\mu+\nu}\Gamma(\mu)\Gamma(\nu)}{a^{\mu}b^{\nu}\pi^3}\int_{0}^{\infty}e^{-st}t^{-\frac{\mu+\nu}{2}}
\Gamma_{\mu,\nu,a,b}(t)dt$$
is a probability density function on $(0,\infty)$ whenever $a,b>0,$ $\mu,\nu>\frac{1}{2}$ and we also have that
$$\int_{0}^{\infty}\left[\int_{0}^{\infty}e^{-st}t^{-\frac{\mu+\nu}{2}}
\Gamma_{\mu,\nu,a,b}(t)dt\right]ds=
\frac{a^{\mu}b^{\nu}\pi^3}{2^{\mu+\nu}\Gamma(\mu)\Gamma(\nu)}.$$

Now, we focus on the quotient of modified Bessel functions of the first and second kind.

\begin{theorem}\label{theoquotIK}
If $a, b>0,$ $\mu>-1$, $\mu+\nu>0$ and $\left|\arg z\right|<\pi,$ then the following representation holds
\begin{equation*}		
\frac{z^{-\frac{\mu+\nu}{2}}}{e^{(a+b)\sqrt{z}}}\frac{I_{\mu}(a\sqrt{z})}{K_{\nu}(b\sqrt{z})}
=-\frac{2}{\pi^2}\int_{0}^{\infty}\frac{t^{-\frac{\mu+\nu}{2}}}{z+t}J_{\mu}(a\sqrt{t})\gamma_{\nu,a,b}(t)dt,
\end{equation*}
which can be rewritten as
$$\frac{z^{-\frac{\mu+\nu}{2}}}{e^{(a+b)\sqrt{z}}}\frac{I_{\mu}(a\sqrt{z})}{K_{\nu}(b\sqrt{z})}
=-\frac{2}{\pi^2}\int_{0}^{\infty}e^{-zs}\left[\int_{0}^{\infty}\frac{e^{-st}J_{\mu}(a\sqrt{t})}{t^{\frac{\mu+\nu}{2}}}
\gamma_{\nu,a,b}(t)dt \right]ds,$$
where
$$\gamma_{\nu,a,b}(t)=\frac{J_{\nu}(b\sqrt{t})\cos((a+b)\sqrt{t})+Y_{\nu}(b\sqrt{t})\sin((a+b)\sqrt{t})}{J^2_{\nu}(b\sqrt{t})+Y^2_{\nu}(b\sqrt{t})}.$$
\end{theorem}

In particular, when $a\to0,$ we obtain the following result.

\begin{corollary}\label{theoquotIKcoro}
If $b>0$, $\nu>\frac{1}{2}$ and $\left|\arg z\right|<\pi,$ then the following representation holds true
$$z^{-\frac{\nu}{2}}e^{-b\sqrt{z}}\dfrac{1}{K_{\nu}(b\sqrt{z})}=-\frac{2}{\pi^2}\int_{0}^{\infty}\dfrac{t^{-\frac{\nu}{2}}}{z+t}\gamma_{\nu,0,b}(t)dt.$$
\end{corollary}

Finally, we show that the above quotient of modified Bessel functions of the first and second kind also generates an infinitely divisible distribution.

\begin{theorem}\label{theoquotIKinfdiv}
If $a,b>0$, $\nu,\mu>\frac{1}{2},$ then the function
\begin{equation*}
\varepsilon_{\mu,\nu,a,b}(x)=\frac{2^{\mu+\nu-1}}{a^{\mu}b^{\nu}}\Gamma(\nu)\Gamma(\mu+1)e^{-(a+b)\sqrt{x}}x^{-\frac{\mu+\nu}{2}}\dfrac{I_{\mu}(a\sqrt{x})}{K_{\nu}(b\sqrt{x})}
\end{equation*}
is the Laplace transform of an infinitely divisible distribution on $(0, \infty)$. Moreover, if $a,b>0$, $\nu>0$ and $\mu>-1,$ then the reciprocal of $e^{(a+b)\sqrt{x}}\varepsilon_{\mu,\nu,a,b}(x),$ that is
$$\frac{e^{-(a+b)\sqrt{x}}}{\varepsilon_{\mu,\nu,a,b}(x)}=\frac{a^{\mu}b^{\nu}}{2^{\mu+\nu-1}\Gamma(\nu)\Gamma(\mu+1)}x^{\frac{\mu+\nu}{2}}\frac{K_{\nu}(b\sqrt{x})}{I_{\mu}(a\sqrt{x})}$$
is also a Laplace transform of an infinitely divisible distribution on $(0,\infty).$
\end{theorem}

Observe that if $a,b>0$ and $\mu, \nu>\frac{1}{2},$ then from Theorems \ref{theoquotIK} and \ref{theoquotIKinfdiv} we obtain that the function
\begin{equation*} s\mapsto\frac{-2^{\mu+\nu}\Gamma(\nu)\Gamma(\mu+1)}{a^{\mu}b^{\nu}\pi^2}\int_{0}^{\infty}e^{-st}t^{-\frac{\mu+\nu}{2}}J_{\mu}(a\sqrt{t})
\gamma_{\nu,a,b}(t)dt
\end{equation*}
is a probability density function and
\begin{equation*}
\int_{0}^{\infty}\left[\int_{0}^{\infty}e^{-st}t^{-\frac{\mu+\nu}{2}}J_{\mu}(a\sqrt{t})
\gamma_{\nu,a,b}(t)dt\right]ds=-\frac{a^{\mu}b^{\nu}\pi^2}{2^{\mu+\nu}\Gamma(\nu)\Gamma(\mu+1)}.
\end{equation*}

\section{\bf Remarks, open problems and challenges for future direction}

\subsection{\bf Remarks on quotients of Tricomi hypergeometric functions}

In this subsection our aim is to point out some facts related to quotients of Tricomi hypergeometric functions. For this recall that Goovaerts et al. \cite{goova} proved that the distribution of the ratio of two independent gamma-distributed random variables is infinitely divisible and this result was in fact a solution to an unsolved problem given by Steutel \cite{steutel0} in a survey about the theory of infinite divisibility. More precisely, Goovaerts et al. \cite{goova} based on a relation between the Laplace transform of the density of the quotient of two gamma random variables and the Tricomi confluent hypergeometric function, proved that the distribution of the ratio of two independent gamma-distributed random variables is a generalized gamma convolution. In their proof the Laplace transform, or rather its logarithmic derivative, is obtained by the argument principle and contour integration. Recall also that Ismail and Kelker \cite[Theorem 1.5]{kelker} proved with the Stieltjes transform technique that the distribution of the two gamma random variables is self-decomposable, hence is infinitely divisible. In this subsection we would like to point out that there is another way to show that the distribution of the ratio of two gamma random variables is a generalized gamma convolution.

\begin{theorem}\label{Thquogamma}
The distribution of the quotient of two gamma random variables is a generalized gamma convolution.
\end{theorem}

Note that Ismail and Kelker's proof of \cite[Theorem 1.5]{kelker} was based on the following integral representation \cite[p. 885]{kelker}
\begin{equation}\label{intfor}
\frac{\psi(a+1,c+1,z)}{\psi(a,c,z)}=\int_0^{\infty}\frac{t^{-c}e^{-t}
\left|\psi(a,c,te^{\mathrm{i}\pi})\right|^{-2}}
{(z+t)\Gamma(a+1)\Gamma(a-c+1)}dt,
\end{equation}
where $a>0$, $c<1$ and $|\arg z|<\pi$. Our proof of Theorem \ref{Thquogamma} is also based on \eqref{intfor}, however our approach is based on the Pick function characterization theorem, that is, Lemma \ref{lemma3}. Moreover, here we would like to show that it is possible to obtain similar Stieltjes transform representations for quotients of Tricomi hypergeometric functions. These results complement the earlier results of Ismail and Kelker \cite[Theorem 1.4]{kelker}.

\begin{theorem}\label{tricrepr}
If $a>0$, $c<1$ and $|\arg z|<\pi,$ then the following representation are valid
\begin{equation}\label{triceqco}\frac{\psi(a,c-1,z)}{\psi(a,c,z)}=\frac{1-c}{a-c+1}+\int_{0}^{\infty}\frac{z}{z+t}\frac{t^{-c}e^{-t}\left|\psi(a,c,e^{\mathrm{i}\pi}t)\right|^{-2}}{\Gamma(a)\Gamma(a-c+2)}dt,\end{equation}
$$\frac{\psi(a+1,c,z)}{\psi(a,c,z)}=\frac{1}{a-c+1}-\int_{0}^{\infty}\frac{z}{z+t}\frac{(a-c+1)t^{-c}e^{-t}}{\Gamma(a+1)\Gamma(a-c+2)}\left|\psi(a,c,e^{\mathrm{i}\pi}t)\right|^{-2}dt,$$
\begin{equation}\label{triceqnew}\frac{\psi(a,c+1,z)}{\psi(a,c,z)}=1+\int_{0}^{\infty}\frac{t^{-c}e^{-t}}{z+t}\frac{\left|\psi(a,c,e^{\mathrm{i}\pi}t)\right|^{-2}}{\Gamma(a)\Gamma(a-c+1)}dt,\end{equation}
$$\frac{\psi(a-1,c,z)}{\psi(a,c,z)}=z-c+a+\int_{0}^{\infty}\frac{z}{z+t}\frac{t^{-c}e^{-t}\left|\psi(a,c,e^{\mathrm{i}\pi}t)\right|^{-2}}{\Gamma(a)\Gamma(a-c+1)}dt.$$	
\end{theorem}

The above representations follow naturally from \eqref{intfor}, however for \eqref{triceqnew} we give a more detailed proof via the Stieltjes representation and inversion theorems. We also mention here that by using the Stieltjes transform technique Ismail and Kelker \cite[eq. (1.5)]{kelker} proved the following result
$$\frac{\psi(a,c-1,z)}{\psi(a,c,z)}=\int_{0}^{\infty}\frac{z}{z+t}\frac{t^{-c}e^{-t}\left|\psi(a,c,e^{\mathrm{i}\pi}t)\right|^{-2}}{\Gamma(a)\Gamma(a-c+2)}dt$$
for $a>0$, $1<c<a+1$ and $|\arg z|<\pi,$ and our result \eqref{triceqco} complements this result naturally.

We also mention that according to Lennart Bondesson (see \cite[Example 4.3.1]{bondesson}) the distribution of the quotient of two gamma random variables belongs also to the class of hyperbolically completely monotone densities. More precisely, if $X$ and $Y$ are independent gamma distributed random variables with parameters $(\alpha, \beta)$ and $(\alpha_0, \beta_0),$ then the probability density function of the quotient of these two random variables $Z={X}/{Y}$ is given by \cite[p. 889]{kelker}
$$f(x)=\frac{\Gamma(\alpha+\alpha_0)}{\Gamma(\alpha)\Gamma(\alpha_0)}\left(\frac{\beta_0}{\beta}\right)^{\alpha}x^{\alpha-1}\left(1+\frac{\beta_0}{\beta}x\right)^{-(\alpha+\alpha_0)},$$
where $x>0.$ Now, if $u, v>0$ and $w=v+{1}/{v},$ then we have that (see \cite[Example 4.3.1]{bondesson})
\begin{align*}
f(uv)f\left(\frac{u}{v}\right)&=\left[\dfrac{\Gamma(\alpha+\alpha_0)}{\Gamma(\alpha)\Gamma(\alpha_0)}\left(\dfrac{\beta_0}{\beta}\right)^{\alpha}\right]^2u^{2\alpha-2}
\left[\left(1+\frac{\beta_0}{\beta}uv\right)\left(1+\frac{\beta_0}{\beta}\frac{u}{v}\right)\right]^{-(\alpha+\alpha_0)}\\
&=\left[\dfrac{\Gamma(\alpha+\alpha_0)}{\Gamma(\alpha)\Gamma(\alpha_0)}\left(\dfrac{\beta_0}{\beta}\right)^{\alpha}\right]^2u^{2\alpha-2}
\left[1+\left(\frac{u\beta_0}{\beta}\right)^2+\frac{\beta_0}{\beta}uw\right]^{-(\alpha+\alpha_0)},
\end{align*}
and thus clearly $w\mapsto f(uv)f\left({u}/{v}\right)$ is completely monotonic on $(0, \infty)$ for all $\alpha,\alpha_0,\beta,\beta_0>0.$

Finally, we note that the integral representation \eqref{intfor} was also quite important in the study of the infinite divisibility of the Whittaker distribution, see \cite{assefa} for more details. Moreover, we also note that it is possible to obtain some similar results on ratios of Tricomi hypergeometric functions where the difference between the parameters is not necessarily an integer, see \cite[Section 2]{ferreira} for more details.

\subsection{The non-central chi square distribution} The non-central chi square distribution is a noncentral generalization of the chi-squared distribution and it has the next probability density function
$$\chi_{\mu, \lambda}(x)=\frac{1}{2}e^{-\frac{(x+\lambda)}{2}}\left(\frac{x}{\lambda}\right)^{\frac{\mu}{4}-\frac{1}{2}}I_{\frac{\mu}{2}-1}(\sqrt{\lambda x}),$$
where $\lambda>0$ is the non-central parameter and $\mu>0$ is the degree of freedom. According to Ismail and Kelker \cite[Theorem 1.6]{kelker} we know that the non-central chi square distribution is infinitely divisible
for all degrees of freedom (including fractional degrees of freedom). Moreover, due to Bondesson \cite[Example 9.2.2]{bondesson} we also know that the non-central chi square distribution belongs to the so-called class of generalized convolutions of mixtures of exponential distributions, introduced by Bondesson \cite[p. 43]{bondesson2} and which is in fact the smallest class of distributions which is closed under convolution and weak convergence and which contains all mixtures of exponential distributions. We know that the class of generalized convolutions of mixtures of exponential distributions is a subclass of the infinitely divisible distributions, however contains all generalized gamma convolutions and hence hyperbolically completely monotone densities. Thus, it is very natural to ask whether the non-central chi square distribution belongs to the class of generalized gamma convolutions or to the class of hyperbolically completely monotone densities. The next result suggest that under some conditions on the parameters $\mu$ and $\lambda$ the non-central chi square distribution belongs to the class of hyperbolically completely monotone densities, and hence to the class of generalized gamma convolutions. {\bf The first open problem} is to find the optimal range for the parameters $\mu$ and $\lambda$ such that the non-central chi square distribution belongs to the class of hyperbolically completely monotone densities.

\begin{theorem}\label{noncentralchihcm}
Let $\mu>1$, $u,v>0$ and $w=v+{1}/{v}.$ If $\lambda\leq\mu$, then $w\mapsto \chi_{\mu, \lambda}(uv)\chi_{\mu, \lambda}\left({u}/{v}\right)$ is strictly decreasing on $(2, \infty)$ and if $2\lambda\leq\mu,$ then $w\mapsto \chi_{\mu, \lambda}(uv)\chi_{\mu, \lambda}\left({u}/{v}\right)$ is strictly convex on $(2, \infty)$.
\end{theorem}

\subsection{Hyperbolically complete monotonicity of McKay distributions} In the previous section, namely in Theorems \ref{th1}, \ref{th2}, \ref{th3} and \ref{th4} we studied the infinite divisibility and self-decomposability
of the McKay type distributions and we also verified whether these distributions belong to the class of generalized gamma convolutions. However, we were not able to check whether these distributions belong to the class of hyperbolically completely monotone densities. In the case of the $K$-distribution in Theorem \ref{thK} the McDonald's integral representation was crucial, however, as far as we know there is no analogous result for modified Bessel functions of first kind. Moreover, the following result shows that $w \mapsto I_{\mu}(uv)I_{\mu}({u}/{v})$ is absolutely monotonic on $(2, \infty)$ and thus a more sophisticated analysis is needed to verify whether the McKay type distributions belong to the class of hyperbolically completely monotone densities. Recall that a function $f:(0,\infty)\to\mathbb{R}$ is {\em absolutely monotone} (or absolutely monotonic) if it has derivatives of all orders and $f^{(n)}(x)>0$ for all $x>0$ and $n\in\{0,1,2,\dots\}.$

\begin{theorem}\label{thprodIabsmon}
If $\mu>-\frac{1}{2}$, $u,v >0$ and $w=v+{1}/{v}$, then the function $w \to I_{\mu}(uv)I_{\mu}\left({u}/{v}\right)$ is absolutely monotonic on $(2, \infty)$.
\end{theorem}

{\bf The second open problem} is to verify under which conditions the distributions considered in Theorems \ref{th1}, \ref{th2}, \ref{th3} and \ref{th4} belong to the class of hyperbolically completely monotone densities, and
under which conditions the distributions considered in Theorems \ref{th2}, \ref{th3} and \ref{th4} belong to the class of self-decomposable distributions as well as of generalized gamma convolutions. Moreover, it would be also a good idea to see under which conditions the probability density function in Theorem \ref{thprod1} is hyperbolically completely monotonic.

\subsection{Self-decomposability of modified Bessel distributions} Now, let $w(x)$ be the Laplace transform of a probability distribution on $[0,\infty).$ It is well-known that if $-dw(x)/dx$ is the Stieltjes transform of a positive measure, then the original distribution is self-decomposable, and hence is infinitely divisible. In other words, a non-negative random variable is self-decomposable if its Laplace transform satisfies the conditions
\begin{equation}\label{selfcond}\int_0^{\infty}e^{-xt}d\varpi(t)=e^{-h(x)},\quad h'(x)=\int_0^\infty\frac{d\varpi(t)}{x+t},\quad d\varpi(t)\geq0.\end{equation}
Moreover, according to Lemma \ref{lemma1} a probability measure $d\omega$ supported on $[0,\infty)$ is infinitely divisible if and only if its Laplace transform satisfies the conditions
\begin{equation}\label{infdivcond}\int_0^{\infty}e^{-xt}d\omega(t)=e^{-h(x)},\quad h(0)=0,\quad \mbox{and}\ h'(x)\ \mbox{is completely monotonic}. \end{equation}
Recall that self-decomposable functions are infinitely divisible and a probability distribution satisfying \eqref{selfcond} and \eqref{infdivcond} is called a generalized gamma convolution.
In Theorems \ref{theprodIKexp}, \ref{thprodeqIinfdiv} and \ref{theoquotIKinfdiv} we have some infinitely divisible modified Bessel distributions, whose Laplace transform can be written as Stieltjes transforms, however these Stieltjes transforms does not have positive kernels. Thus, {\bf the third open problem} is to verify under which conditions the distributions considered in Theorems \ref{theprodIKexp}, \ref{thprodeqIinfdiv} and \ref{theoquotIKinfdiv} are self-decomposable.

\section{\bf Proofs of the main results}

\begin{proof}[\bf Proof of Theorem \ref{th1}]
According to Prudnikov et al. \cite[eq. 2.15.3.2]{prudnikov} we have that
$$\int_0^{\infty}x^{\mu}e^{-bx}I_{\mu}(ax)dx=\frac{(2a)^{\mu}\Gamma\left(\mu+\frac{1}{2}\right)}{\sqrt{\pi}(b^2-a^2)^{\mu+\frac{1}{2}}},$$
which in turn implies that the Laplace transform
$$L[\varphi_{\mu,a,b}(x)]=\int_0^{\infty}e^{-xt}\varphi_{\mu,a,b}(t)dt$$
of the probability density function $\varphi_{\mu,a,b}$ is given by
$$L[\varphi_{\mu,a,b}(x)]=\frac{\sqrt{\pi}(b^2-a^2)^{\mu+\frac{1}{2}}}{(2a)^{\mu}\Gamma\left(\mu+\frac{1}{2}\right)}
\int_0^{\infty}t^{\mu}e^{-(b+x)t}I_{\mu}(at)dt=\left[\frac{b^2-a^2}{(x+b)^2-a^2}\right]^{\mu+\frac{1}{2}}.$$

First we show that the McKay distribution, with probability density function in \eqref{mcnoltypdf},
is infinitely divisible. In view of Lemma \ref{lemma1}, the infinite divisibility is equivalent to the complete monotonicity on $(0,\infty)$ of
$$x\mapsto -\frac{d \ln L[\varphi_{\mu,a,b}(x)]}{dx}=\left(\mu+\frac{1}{2}\right)\left(\frac{1}{x+b-a}+\frac{1}{x+b+a}\right)$$
when $\mu>-\frac{1}{2}$ and $b>a>0,$ which is certainly true.

Now, we show that the McKay distribution, with probability density function in \eqref{mcnoltypdf}, belongs
to the self-decomposable subclass of infinitely divisible distributions. In view of Lemma \ref{lemma2} and the notation
$\omega_{\mu,a,b}(x)=L[\varphi_{\mu,a,b}(x)]/L[\varphi_{\mu,a,b}(\alpha x)]$, the self-decomposability is equivalent to the
complete monotonicity on $(0,\infty)$ of
$$x\mapsto-\frac{d \ln \omega_{\mu,a,b}(x)}{dx}=\left(\mu+\frac{1}{2}\right)\left[\frac{(1-\alpha)(b-a)}{(x+b-a)(\alpha x+b-a)}+\frac{(1-\alpha)(b+a)}{(x+b+a)(\alpha x+b+a)}\right]$$
when $\alpha\in(0,1),$ $\mu>-\frac{1}{2}$ and $b>a>0,$ which is certainly true.

Finally, we show that the McKay distribution in question belongs also to the class of generalized gamma convolutions. The moment generating function of this McKay distribution is
$$M_X(z)=\int_0^{\infty}e^{tz}\varphi_{\mu,a,b}(t)dt=L[\varphi_{\mu,a,b}(-z)]=\left[\frac{b^2-a^2}{(z-b)^2-a^2}\right]^{\mu+\frac{1}{2}}.$$
Of course this moment generating function $\psi=M_X$ is analytic and zero-free in $\mathbb{C}\setminus[0,\infty)$ and in order to apply Lemma \ref{lemma3} we just need to verify the inequality $\imag\left[\psi'(s)/\psi(s)\right]\geq0$ for $\imag s>0.$ In this case we have that
$$\frac{\psi'(s)}{\psi(s)}=-\left(\mu+\frac{1}{2}\right)\left(\frac{1}{s+a-b}+\frac{1}{s-a-b}\right),$$
which in turn implies that for $s=x+\mathrm{i}y$ and $y=\imag s>0$ we certainly have
$$\imag \frac{\psi'(s)}{\psi(s)}=\left(\mu+\frac{1}{2}\right)\left(\frac{y}{(x+a-b)^2+y^2}+\frac{y}{(x-a-b)^2+y^2}\right)>0.$$
This shows that the logarithmic derivative of the moment generating function is a Pick function and thus Lemma \ref{lemma3} shows that the McKay distribution with probability density function
in \eqref{mcnoltypdf} belongs to the class of generalized gamma convolutions.

We mention that there is another proof for the last affirmation. Namely, if we apply Lemma \ref{lemma4}, then we just need to observe that if $\phi(x)=L[\varphi_{\mu,a,b}(x)],$ then $\phi(uv)\phi(u/v)$ can be written as
$$\phi(uv)\phi(u/v)=\frac{(b^2-a^2)^{\mu+\frac{1}{2}}u^{-2\mu-1}}{\left[w+\frac{u^2+(b-a)^2}{(b-a)u}\right]^{\mu+\frac{1}{2}}\left[w+\frac{u^2+(b+a)^2}{(b+a)u}\right]^{\mu+\frac{1}{2}}}$$
and this is clearly completely monotonic in $w=v+1/v>0$ for all $\mu>-\frac{1}{2}$ and $b>a>0.$
\end{proof}

\begin{proof}[\bf Proof of Theorem \ref{th2}]
According to Prudnikov et al. \cite[eq. 2.15.3.2]{prudnikov} we have that
$$\int_0^{\infty}x^{\mu+1}e^{-bx}I_{\mu}(ax)dx=\frac{2b(2a)^{\mu}\Gamma\left(\mu+\frac{3}{2}\right)}{\sqrt{\pi}(b^2-a^2)^{\mu+\frac{3}{2}}},$$
which in turn implies that the Laplace transform
$$L[\psi_{\mu,a,b}(x)]=\int_0^{\infty}e^{-xt}\psi_{\mu,a,b}(t)dt$$
of the probability density function $\psi_{\mu,a,b}$ is given by
\begin{align*}L[\psi_{\mu,a,b}(x)]=\frac{\sqrt{\pi}(b^2-a^2)^{\mu+\frac{3}{2}}}{2b(2a)^{\mu}\Gamma\left(\mu+\frac{3}{2}\right)}
\int_0^{\infty}t^{\mu+1}e^{-(b+x)t}I_{\mu}(at)dt=\left(1+\frac{x}{b}\right)\left[\frac{b^2-a^2}{(x+b)^2-a^2}\right]^{\mu+\frac{3}{2}}.\end{align*}

Now, we show that the distribution, with probability density function in \eqref{newpdf},
is infinitely divisible. In view of Lemma \ref{lemma1}, the infinite divisibility is equivalent to the complete monotonicity on $(0,\infty)$ of
$$x\mapsto-\frac{d \ln L[\psi_{\mu,a,b}(x)]}{dx}=\frac{\mu+1}{x+b-a}+\frac{\mu+1}{x+b+a}+\frac{a^2}{(x+b)(x+b-a)(x+b+a)}$$
when $\mu>-1$ and $b>a>0,$ which is certainly true.
\end{proof}

\begin{proof}[\bf Proof of Theorem \ref{th3}]
By using the integral formula (see for example the book of Prudnikov et al. \cite[eq. 2.15.3.2]{prudnikov})
$$\int_0^{\infty}x^{\nu-1}e^{-bx}I_{\mu}(ax)dx=c_{\mu,\nu,a,b}$$
we obtain that the Laplace transform
$$L[\varphi_{\mu,\nu,a,b}(x)]=\int_0^{\infty}e^{-xt}\varphi_{\mu,\nu,a,b}(t)dt$$
of the probability density function $\varphi_{\mu,\nu,a,b}$ is given by
\begin{align*}L[\varphi_{\mu,\nu,a,b}(x)]=\left(\frac{b}{x+b}\right)^{\mu+\nu}\cdot
\frac{{}_2F_1\left(\frac{\mu+\nu}{2},\frac{\mu+\nu+1}{2},\mu+1,\frac{a^2}{(x+b)^2}\right)}
{{}_2F_1\left(\frac{\mu+\nu}{2},\frac{\mu+\nu+1}{2},\mu+1,\frac{a^2}{b^2}\right)}.\end{align*}
Recall that the Gaussian hypergeometric function satisfies
$$\frac{d}{dx}\left[{}_2F_1(a,b,c,x)\right]=\frac{ab}{c}\cdot{}_2F_1(a+1,b+1,c+1,x)$$
and for the parameters $a,$ $b$ and $c$ such that $-1\leq a \leq c$ and $0<b\leq c$ the next integral representation of K\"ustner \cite{kustner}
\begin{equation}\label{intkust}
\frac{z\cdot {}_2F_1(a+1,b+1,c+1,z)}{{}_2F_1(a,b,c,z)}=\int_0^1\frac{z}{1-tz}dq_{a,b,c}(t)
\end{equation}
is valid, where $z\in\mathbb{C}\setminus[1,\infty),$ $q_{a,b,c}(1)-q_{a,b,c}(0)=1$ and $q_{a,b,c}$ is non-decreasing self-mapping of $[0,1].$ To prove that the generalization
of the McKay distribution is infinitely divisible, in view of Lemma \ref{lemma1}, we just need to observe that the expression
\begin{align*}x\mapsto-\frac{d\ln L\left[\varphi_{\mu,\nu,a,b}(x)\right]}{dx}&=\frac{\mu+\nu}{x+b}+\frac{a^2(\mu+\nu)(\mu+\nu+1)}{2(\mu+1)(x+b)^3}
\frac{{}_2F_1\left(\frac{\mu+\nu+2}{2},\frac{\mu+\nu+3}{2},\mu+2,\frac{a^2}{(x+b)^2}\right)}
{{}_2F_1\left(\frac{\mu+\nu}{2},\frac{\mu+\nu+1}{2},\mu+1,\frac{a^2}{(x+b)^2}\right)}\\
&=\frac{\mu+\nu}{x+b}\left[1+\frac{a^2(\mu+\nu+1)}{2(\mu+1)}\int_0^1\frac{dq_{\mu,\nu}(t)}{(x+b-a\sqrt{t})(x+b+a\sqrt{t})}\right]\end{align*}
is completely monotonic on $(0,\infty)$ for all $\mu+\nu>0,$ $\nu\leq\mu+1$ and $b>a>0,$ where in view of \eqref{intkust} we have that $q_{\mu,\nu}(1)-q_{\mu,\nu}(0)=1$ and $q_{\mu,\nu}$ is a non-decreasing self-mapping of $[0,1].$
\end{proof}

\begin{proof}[\bf Proof of Theorem \ref{th4}]
By using the integral formula (see for example the book of Prudnikov et al. \cite[eq. 2.15.20.5]{prudnikov})
$$\int_0^{\infty}x^{2\mu}e^{-bx}\left[I_{\mu}(ax)dx\right]^2=c_{\mu,a,b}$$
we obtain that the Laplace transform
$$L[\xi_{\mu,a,b}(x)]=\int_0^{\infty}e^{-xt}\xi_{\mu,a,b}(t)dt$$
of the probability density function $\xi_{\mu,a,b}$ is given by
\begin{align*}L[\xi_{\mu,a,b}(x)]=\left(\frac{b}{x+b}\right)^{4\mu+1}\cdot
\frac{{}_2F_1\left(\mu+\frac{1}{2},2\mu+\frac{1}{2},\mu+1,\frac{4a^2}{(x+b)^2}\right)}
{{}_2F_1\left(\mu+\frac{1}{2},2\mu+\frac{1}{2},\mu+1,\frac{4a^2}{b^2}\right)}.\end{align*}
To prove that the McKay type distribution with probability density function in \eqref{pdfxi} is infinitely divisible, in view of Lemma \ref{lemma1}, we just need to observe that the expression
\begin{align*}x\mapsto-\frac{d\ln L\left[\xi_{\mu,a,b}(x)\right]}{dx}&=\frac{4\mu+1}{x+b}+\frac{8a^2\left(\mu+\frac{1}{2}\right)\left(2\mu+\frac{1}{2}\right)}{(\mu+1)(x+b)^3}
\frac{{}_2F_1\left(\mu+\frac{3}{2},2\mu+\frac{3}{2},\mu+2,\frac{4a^2}{(x+b)^2}\right)}
{{}_2F_1\left(\mu+\frac{1}{2},2\mu+\frac{1}{2},\mu+1,\frac{4a^2}{(x+b)^2}\right)}\\
&=\frac{4\mu+1}{x+b}\left[1+\frac{2a^2(2\mu+1)}{\mu+1}\int_0^1\frac{dq_{\mu}(t)}{(x+b-2a\sqrt{t})(x+b+2a\sqrt{t})}\right]\end{align*}
is completely monotonic on $(0,\infty)$ for all $-1\leq\mu+\frac{1}{2}\leq\mu+1,$ $0<2\mu+\frac{1}{2}\leq\mu+1$ and $b>2a>0,$ where in view of \eqref{intkust} we have that $q_{\mu}(1)-q_{\mu}(0)=1$ and $q_{\mu}$ is a non-decreasing self-mapping of $[0,1].$
\end{proof}

\begin{proof}[\bf Proof of Theorem \ref{thK}]
We will consider only the case when $\alpha \neq \beta$ for the proof of the infinite divisibility, self-decomposability as well as when we prove that the $K$-distribution belongs to the class of generalized gamma convolutions. However, in the case of hyperbolically complete monotonicity we will consider also the case when $\alpha=\beta.$

By using Lemma \ref{lemma1}, first we show that the $K$-distribution is infinitely divisible. In view of the integral formula (see for example the book of Prudnikov et al. \cite[eq. 2.16.8.4]{prudnikov})
$$\int_0^{\infty}t^{q-1/2}e^{-rt}K_{2\nu}(2s\sqrt{t})dt=\frac{\Gamma\left(q+\nu+\frac{1}{2}\right)\Gamma\left(q-\nu+\frac{1}{2}\right)}{2sr^qe^{-\frac{s^2}{2r}}}W_{-q,\nu}\left(\frac{s^2}{r}\right),$$
where $W_{\kappa,\mu}$ stands for the Whittaker function of the second kind, we arrive at
$$L\left[\omega_{\alpha,\beta,\mu}(x)\right]=\int_0^{\infty}e^{-xt}\omega_{\alpha,\beta,\mu}(t)dt=\left(\frac{\alpha\beta}{\mu x}\right)^{\frac{\alpha+\beta-1}{2}}e^{\frac{\alpha\beta}{2\mu x}}W_{-\frac{\alpha+\beta-1}{2},\frac{\alpha-\beta}{2}}\left(\frac{\alpha\beta}{\mu x}\right).$$
Observe that in view of the connection between the Whittaker function of the second kind and Tricomi hypergeometric function
$$\label{wtrel}W_{\kappa,\mu}(x)=e^{-\frac{x}{2}}x^{\mu+\frac{1}{2}}\cdot\psi\left(\mu-\kappa+\frac{1}{2},1+2\mu,x\right),$$
we obtain that the Laplace transform of the $K$-distribution is
$$L\left[\omega_{\alpha,\beta,\mu}(x)\right]=\left(\frac{\alpha\beta}{\mu x}\right)^{\alpha}\psi\left(\alpha,1+\alpha-\beta,\frac{\alpha\beta}{\mu x}\right).$$
Now, recall the recurrence relation \cite[eq. 13.3.22]{nist}
$$\psi'(a,c,x)=-a\psi(a+1,c+1,x)$$
and the integral representation \eqref{intfor}, which is valid for $|\arg z|<\pi,$ $a>0$ and $c<1.$ In view of Lemma \ref{lemma1} and the above relations, to show that the $K$-distribution is infinitely divisible we
just need to show that for $\theta(x)=L\left[\omega_{\alpha,\beta,\mu}(x)\right]$ we have that
$$-\frac{d}{dx}\left[\ln\theta(x) \right]=\frac{\alpha}{x}+\frac{\alpha\beta}{\mu x^2}\frac{\psi'\left(\alpha,1+\alpha-\beta,\frac{\alpha\beta}{\mu x}\right)}{\psi\left(\alpha,1+\alpha-\beta,\frac{\alpha\beta}{\mu x}\right)}=\frac{\alpha}{x}-\frac{\alpha^2\beta}{\mu x^2}\frac{\psi\left(\alpha+1,2+\alpha-\beta,\frac{\alpha\beta}{\mu x}\right)}{\psi\left(\alpha,1+\alpha-\beta,\frac{\alpha\beta}{\mu x}\right)},$$
that is,
$$-\frac{d}{dx}\left[\ln\theta(x) \right]=\frac{\alpha}{x}\left[1-\frac{\alpha\beta}{\mu x}\int_0^{\infty}\frac{\omega_{\alpha,\beta}(t)}{\frac{\alpha\beta}{\mu x}+t}dt\right],$$
where
$$\omega_{\alpha,\beta}(t)=\frac{t^{\beta-\alpha-1}e^{-t}}{\Gamma(\alpha+1)\Gamma(\beta)}\left|\psi\left(\alpha,1+\alpha-\beta,te^{\mathrm{i}\pi}\right)\right|^{-2},$$
is completely monotonic on $(0,\infty).$ To show this observe that if in
\begin{equation}\label{ltKdist}\frac{\alpha\beta}{\mu x}\frac{\psi\left(\alpha+1,2+\alpha-\beta,\frac{\alpha\beta}{\mu x}\right)}{\psi\left(\alpha,1+\alpha-\beta,\frac{\alpha\beta}{\mu x}\right)}=\int_0^{\infty}\frac{\frac{\alpha\beta}{\mu x}\omega_{\alpha,\beta}(t)}{\frac{\alpha\beta}{\mu x}+t}dt\end{equation}
we make the change of variable $(\alpha\beta)/(\mu x)=s$ and $s$ tends to infinity, then in view of the asymptotic expansion \cite[eq. 13.7.3]{nist}
$$\psi(a,c,x)\sim x^{-a}\left(1+a(c-a-1)\frac{1}{x}+\frac{1}{2}a(a+1)(a+1-c)(a+2-c)\frac{1}{x^2}+\dots\right),$$
which is valid for large real $x$ and fixed $a$ and $c,$ we clearly have that
\begin{equation}\label{pdfome}\int_0^{\infty}\omega_{\alpha,\beta}(t)dt=1.\end{equation}
This implies that
$$-\frac{d}{dx}\left[\ln\theta(x) \right]=\frac{\alpha}{x}\left[\int_0^{\infty}\omega_{\alpha,\beta}(t)dt-\frac{\alpha\beta}{\mu x}\int_0^{\infty}\frac{\omega_{\alpha,\beta}(t)}{\frac{\alpha\beta}{\mu x}+t}dt\right]=
\int_0^{\infty}\frac{\alpha\omega_{\alpha,\beta}(t)}{x+\frac{\alpha\beta}{\mu t}}dt,$$
and this is completely monotonic in $x$ on $(0,\infty)$ for all $\alpha,\beta,\mu>0$ such that $\alpha<\beta.$ All we need is to observe that the
Whittaker function of the first kind is symmetric in the second parameter, that is, $W_{\kappa,\mu}=W_{\kappa,-\mu}$ and this shows that for every $\alpha,\beta,\mu>0$ the
Laplace transform of the $K$-distribution $L\left[\omega_{\alpha,\beta,\mu}(x)\right]$ is the Laplace transform of an infinitely divisible distribution. Another way to show the complete monotonicity of
$-d\left[\ln \theta(x)\right]/dx$ is to use the Kummer transformation \cite[eq. 13.2.40]{nist}
$$\psi(a,c,x)=x^{1-c}\psi(a-c+1,2-c,x)$$
and then to apply the previous approach for the next version of the Laplace transform
\begin{equation}\label{e1K}L\left[\omega_{\alpha,\beta,\mu}(x)\right]=\left(\frac{\alpha\beta}{\mu x}\right)^{\beta}\psi\left(\beta,1-\alpha+\beta,\frac{\alpha\beta}{\mu x}\right).\end{equation}
More precisely, by using the above approach it is possible to verify that the above function is the Laplace transform of infinite divisible distribution for all $\alpha, \beta, \mu>0$ such that $\beta<\alpha.$

Next, we show the self-decomposability of the $K$-distribution. For this we show that for $\theta(x)=L\left[\omega_{\alpha,\beta,\mu}(x)\right]$ we have that for every $a,$ $a\in(0,1),$ the function $\theta(x)/\theta(ax)$ is a Laplace transform. In view of Lemma \ref{lemma1} we have that
$\eta(x)=\theta(x)/\theta(ax)$ is a Laplace transform if and only if $-{d}\left[\ln\eta(x)\right]/dx$ is completely monotonic on $(0,\infty).$ By using again the integral formula \eqref{ltKdist} we arrive at
\begin{align*}
-\frac{d}{dx}\left[\ln\eta(x)\right]&=-\frac{\alpha\beta}{a\mu x^2}\frac{\psi'\left(\alpha,1+\alpha-\beta,\frac{\alpha\beta}{a\mu x}\right)}{\psi\left(\alpha,1+\alpha-\beta,\frac{\alpha\beta}{a\mu x}\right)}+\frac{\alpha\beta}{\mu x^2}\frac{\psi'\left(\alpha,1+\alpha-\beta,\frac{\alpha\beta}{\mu x}\right)}{\psi\left(\alpha,1+\alpha-\beta,\frac{\alpha\beta}{\mu x}\right)}\\&=
\frac{\alpha^2\beta}{a\mu x^2}\frac{\psi\left(\alpha+1,2+\alpha-\beta,\frac{\alpha\beta}{a\mu x}\right)}{\psi\left(\alpha,1+\alpha-\beta,\frac{\alpha\beta}{a\mu x}\right)}-\frac{\alpha^2\beta}{\mu x^2}\frac{\psi\left(\alpha+1,2+\alpha-\beta,\frac{\alpha\beta}{\mu x}\right)}{\psi\left(\alpha,1+\alpha-\beta,\frac{\alpha\beta}{\mu x}\right)}\\&=
\frac{\alpha}{x}\left[\int_0^{\infty}\frac{\frac{\alpha\beta}{a\mu x}\omega_{\alpha,\beta}(t)}{\frac{\alpha\beta}{a\mu x}+t}dt-\int_0^{\infty}\frac{\frac{\alpha\beta}{\mu x}\omega_{\alpha,\beta}(t)}{\frac{\alpha\beta}{\mu x}+t}dt\right]\\&=\int_0^{\infty}\frac{\alpha^2\beta(1-a)\omega_{\alpha,\beta}(t)}{a\mu t\left(x+\frac{\alpha\beta}{t\mu}\right)\left(x+\frac{\alpha\beta}{at\mu}\right)}dt,
\end{align*}
which is clearly completely monotonic as a function of $x$ on $(0,\infty)$, $a\in(0,1)$, for all $\alpha,\beta,\mu>0$ and $\alpha<\beta.$ This in turn implies that $\eta(x)=\theta(x)/\theta(ax)$ is a Laplace transform and in view of Lemma \ref{lemma2} we conclude that if $\alpha<\beta$ the $K$-distribution is self-decomposable. Moreover, we can verify that if $\beta<\alpha$ then the $K$-distribution is also self-decomposable, and this conclusion follows by repeating the above process for the other version of the Laplace transform \eqref{e1K}.

Now, we prove that the $K$-distribution belongs to the class of generalized gamma convolutions and to do this we shall use Lemma \ref{lemma3}. The moment generating function
$$M_X(s)=L\left[\omega_{\alpha,\beta,\mu}(-s)\right]=\left(-\frac{\alpha\beta}{\mu s}\right)^{\alpha}\psi\left(\alpha,1+\alpha-\beta,-\frac{\alpha\beta}{\mu s}\right)$$
of the $K$-distribution is analytic and zero-free in $\mathbb{C}\setminus[0,\infty)$ and we are going to show that satisfies $\imag\left[\psi'(s)/\psi(s)\right]\geq0$ for $\imag s>0.$ In view of \eqref{ltKdist} and \eqref{pdfome} for $\psi(s)=M_X(s)$ and $\alpha<\beta$ we obtain that
$$\frac{\psi'(s)}{\psi(s)}=-\frac{\alpha}{s}-\frac{\alpha^2\beta}{\mu s^2}\frac{\psi\left(\alpha+1,2+\alpha-\beta,-\frac{\alpha\beta}{\mu s}\right)}{\psi\left(\alpha,1+\alpha-\beta,-\frac{\alpha\beta}{\mu s}\right)}=-\frac{\alpha}{s}+\frac{\alpha}{s}\int_0^{\infty}\frac{-\frac{\alpha\beta}{\mu s}\omega_{\alpha,\beta}(t)}{-\frac{\alpha\beta}{\mu s}+t}dt,$$
that is
$$\frac{\psi'(s)}{\psi(s)}=-\alpha\int_0^{\infty}\frac{\omega_{\alpha,\beta}(t)}{\left(s-\frac{\alpha\beta}{\mu t}\right)}dt$$
and this implies that for $s=x+\mathrm{i}y$ such that $y>0$ we have
$$\imag \left[\frac{\psi'(s)}{\psi(s)}\right]=\int_0^{\infty}\frac{\alpha y\omega_{\alpha,\beta}(t)}{\left[\left(x-\frac{\alpha\beta}{\mu t}\right)^2+y^2\right]}dt>0.$$
Thus, in view of Lemma \ref{lemma3} for $\alpha,\beta,\mu>0$ the $K$-distribution indeed belongs to the class of generalized gamma convolutions in the case when $\alpha<\beta.$ In the case when $\alpha>\beta$ we have a similar proof by repeating the above process for the other version of the Laplace transform \eqref{e1K}.

Finally, we show that the probability density function of $K$-distribution is hyperbolically completely monotone. By using McDonald's integral representation \cite[p. 439]{watson}
\begin{equation}\label{prodK}
K_{\mu}(x)K_{\mu}(y)=\frac{1}{2}\int_0^{\infty}\exp\left(-\frac{t}{2}-\frac{x^2+y^2}{2t}\right)K_{\mu}\left(\frac{xy}{t}\right)\frac{dt}{t},
\end{equation}
where $x,y>0$ and $\mu\in\mathbb{R},$ we obtain that $\omega_{\alpha,\beta,\mu}(uv)\omega_{\alpha,\beta,\mu}\left({u}/{v}\right)$ is equal to
$$\frac{2u^{\alpha+\beta-2}}{\Gamma^2(\alpha)\Gamma^2(\beta)}\left(\frac{\alpha\beta}{\mu}\right)^{\alpha+\beta}
\int_0^{\infty}\exp\left(-\frac{t}{2}-\frac{w}{\mu t}\left(2\alpha\beta u\right)\right)K_{\alpha-\beta}\left(\frac{4\alpha\beta u}{\mu t}\right)\frac{dt}{t},$$
and this is clearly completely monotonic on $(0,\infty)$ as a function of $w=v+1/v$ for all $\alpha,\beta,\mu,u>0.$
\end{proof}

\begin{proof}[\bf Proof of Theorem \ref{Thnewgigd}]
We are going to prove the result in question by applying Lemma \ref{lemma4}. To do this, observe that the Laplace transform of the generalized inverse Gaussian distribution is given by
$$L\left[\pi_{\mu,a,b}(x)\right]=\left(\frac{a}{2x+a}\right)^{\frac{\mu}{2}}\cdot\frac{K_{\mu}\sqrt{b(2x+a)}}{K_{\mu}(\sqrt{ab})}.$$
Now, if denote the above expression by $\phi(x),$ then in view of McDonald's integral representation \eqref{prodK} we clearly have that
$$\phi(uv)\phi\left({u}/{v}\right)=\frac{1}{2}\frac{a^{\mu}\left[K_{\mu}(\sqrt{ab})\right]^{-2}}{(\alpha w+\beta)^{\mu/2}}\int_0^{\infty}e^{-\frac{t}{2}-\frac{b}{t}(uw+a)}K_{\mu}(b\sqrt{\alpha w+\beta})\frac{dt}{t},$$
where $\alpha=2au>0$ and $\beta=4u^2+a^2>0.$ On the other hand, we know that \cite[p. 589]{bariczed} the function $x\mapsto x^{\mu/2}K_{\mu}(\sqrt{x})$ is completely monotonic on $(0,\infty)$ for all $\mu\in\mathbb{R}.$ This implies that the function $w\mapsto (b\sqrt{\alpha w+\beta})^{\mu}K_{\mu}(b\sqrt{\alpha w+\beta})$ is also completely monotonic on $(0,\infty)$ for all $\mu\in\mathbb{R}$ and $b,\alpha,\beta>0.$ Now, since $w\mapsto (\alpha w+\beta)^{-\mu}$ is completely monotonic on $(0,\infty)$ for all $\mu\geq0$ and $\alpha,\beta>0,$ and $w\mapsto e^{-\frac{b}{t}(uw+a)}$ is completely monotonic on $(0,\infty)$ for all $a,b,u,t>0,$ we conclude that indeed the funcion $w\mapsto \phi(uv)\phi\left({u}/{v}\right)$ is completely monotonic on $(0,\infty)$ for all $\mu\geq0$ and $a,b>0.$ For the case when $\mu<0$ we just need to observe that $K_{\mu}(x)=K_{-\mu}(x)$ and in view of the above mentioned results $x\mapsto x^{-\mu/2}K_{\mu}(\sqrt{x})$ is completely monotonic on $(0,\infty)$ for all $\mu<0,$ which implies that the function $w\mapsto (b\sqrt{\alpha w+\beta})^{-\mu}K_{\mu}(b\sqrt{\alpha w+\beta})$ is also completely monotonic on $(0,\infty)$ for all $\mu<0$ and $b,\alpha,\beta>0.$ Thus, we conclude that the funcion $w\mapsto \phi(uv)\phi\left({u}/{v}\right)$ is completely monotonic on $(0,\infty)$ also for all $\mu<0$ and $a,b>0.$ Now, applying Lemma \ref{lemma4}, the proof is complete.
\end{proof}

\begin{proof}[\bf Proof of Theorem \ref{thiskellap1}]
We shall prove the result by verifying the conditions in Lemma \ref{lemma3}. Due to Ismail and Kelker \cite[Theorem 1.9]{kelker} we know that $x\mapsto\rho_{\mu,a}(x)$ is the Laplace transform of a probability distribution, and therefore the moment generating function $s\mapsto\phi_{\mu,a}(s)$ is as follows
$$\phi_{\mu,a}(s)=\rho_{\mu,a}(-s)=\left(\frac{a}{2}\right)^{\mu}\frac{(-s)^{\frac{\mu}{2}}}{\Gamma(\mu+1)}\cdot\frac{1}{I_{\mu}(a\sqrt{-s})}.$$
By using the relation $I_{\mu}(\mathrm{i}x)={\mathrm{i}}^{\mu}J_{\mu}(x)$ we arrive at
\begin{equation}\label{phidef}
\phi_{\mu,a}(s)=\left(\frac{a}{2}\right)^{\mu}\frac{s^{\frac{\mu}{2}}}{\Gamma(\mu+1)}\cdot\frac{1}{J_{\mu}(a\sqrt{s})},
\end{equation}
where $x\mapsto J_{\mu}(x)$ stands for the Bessel function of the first kind of order $\mu$. Therefore $s\mapsto\phi_{\mu,a}(s)$ is analytic and zero-free in $\mathbb{C}\setminus[0, \infty)$. Now, taking the logarithmic derivative of both sides of the equation \eqref{phidef}, we obtain that
$$\frac{\phi_{\mu,a}'(s)}{\phi_{\mu,a}(s)} = \frac{\mu}{2s} - \frac{a}{2\sqrt{s}}\frac{J'_{\mu}(a\sqrt{s})}{J_{\mu}(a\sqrt{s})}.$$
Now, taking the logarithmic derivative of both sides of the well-known infinite product representation
$$J_{\mu}(x) = \frac{\left(\frac{1}{2}x\right)^{\mu}}{\Gamma(\mu+1)} \prod_{n\geq1}\left(1-\frac{x^2}{j^2_{\mu, n}}\right),$$
where $j_{\mu,n}$ stands for the $n$th positive zero of $x\mapsto J_{\mu}(x)$, we obtain the classical Mittag-Leffler expansion
$$\frac{J'_{\mu}(x)}{J_{\mu}(x)} = \frac{\mu}{x} - \sum_{n\geq1}\frac{2x}{j^2_{\mu, n}-x^2},$$
and this in turn implies that
$$\imag\left[\frac{\phi'_{\mu,a}(s)}{\phi_{\mu,a}(s)}\right] = \imag\left[\sum_{n\geq1}\frac{a^2}{j^2_{\mu, n}-a^2s}\right]=\sum_{n\geq1}\dfrac{a^4y}{(j^2_{\mu,n}-a^2x)^2+(a^2y)^2}>0,$$
whenever $x=\real s\in\mathbb{R}$ and $y=\imag s>0.$ Consequently, the conditions in Lemma \ref{lemma3} are satisfied and with this the proof is complete.
\end{proof}

\begin{proof}[\bf Proof of Theorem \ref{theolap1}]
Since $2^{\mu}\Gamma(\mu+1)x^{-\mu}I_{\mu}(x)\to1$ as $x\to0,$ it follows that $\Omega_{\mu,\nu,\sigma,a,b}(x)\to1$ as $x\to0.$ Moreover, by using the well-known recurrence relation \cite[p. 79]{watson}
$$I'_{\mu}(x)=I_{\mu+1}(x)+\frac{\mu}{x}I_{\mu}(x),$$
and the Mittag-Leffler expansion
$$\frac{I_{\mu+1}(x)}{I_{\mu}(x)}= \sum_{n\geq1}\frac{2x}{x^2+j^2_{\mu, n}},$$
we obtain that $-{d\ln\Omega_{\mu,\nu,\sigma,a,b}(x)}/{dx}$ becomes
\begin{align*}
-&\frac{a}{2\sqrt{x}}\frac{I'_{\mu}(a\sqrt{x})}{I_{\mu}(a\sqrt{x})}-\frac{b}{2\sqrt{x}}\frac{I'_{\nu}(b\sqrt{x})}{I_{\nu}(b\sqrt{x})}+
\frac{b}{2\sqrt{x}}\frac{I'_{\mu}(b\sqrt{x})}{I_{\mu}(b\sqrt{x})}+\frac{a}{2\sqrt{x}}\frac{I'_{\nu}(a\sqrt{x})}{I_{\nu}(a\sqrt{x})}+\sum_{n\geq1}\left(\frac{b^2}{b^2x+j^2_{\sigma,n}}\right)\\
&=\sum_{n\geq1} \left[-\frac{1}{x+j^2_{\mu, n}a^{-2}}-\frac{1}{x+j^2_{\nu, n}b^{-2}}+\frac{1}{x+j^2_{\mu, n}b^{-2}}+\frac{1}{x+j^2_{\nu, n}a^{-2}}+\frac{1}{x+j^2_{\sigma, n}b^{-2}}\right]\\
&=\int_{0}^{\infty}e^{-xt}\left[\sum_{n\geq1}\left(-e^{-j^2_{\mu, n}a^{-2}t}-e^{-j^2_{\nu, n}b^{-2}t}+e^{-j^2_{\mu, n}b^{-2}t}+e^{-j^2_{\nu, n}a^{-2}t}+e^{-j^2_{\sigma, n}b^{-2}t}\right)\right]dt\\
&=\int_{0}^{\infty}e^{-xt}\left[\sum_{n\geq1}\left((-e^{-j^2_{\mu, n}a^{-2}t}+e^{-j^2_{\mu, n}b^{-2}t})+(-e^{-j^2_{\nu, n}b^{-2}t}+e^{-j^2_{\sigma, n}b^{-2}t})+e^{-j^2_{\nu, n}a^{-2}t}\right)\right]dt.
\end{align*}
Since $a<b$ and $-1<\sigma < \nu$, we have $e^{-j^2_{\mu, n}a^{-2}t} < e^{-j^2_{\mu, n}b^{-2}t}$ and $e^{-j^2_{\nu, n}b^{-2}t} < e^{-j^2_{\sigma, n}b^{-2}t}$, respectively, for every $n \in \mathbb{N}$, $\mu>-1,$ where we used the well-known fact that $\mu \mapsto j_{\mu,n}$ is increasing for every $n \in \mathbb{N}$ on $(-1, \infty)$. Consequently, the last expression is positive, and this implies that $-{d\ln\Omega_{\mu,\nu,\sigma,a,b}(x)}/{dx}$ is completely monotonic in $x$ on $(0,\infty)$. In view of Lemma \ref{lemma1} we conclude that $\Omega_{\mu,\nu,\sigma,a,b}(x)$ is indeed the Laplace transform of an infinitely divisible distribution.
\end{proof}

\begin{proof}[\bf Proof of Theorem \ref{theolap2}]
The result will be established by verifying the conditions in Lemma \ref{lemma1} and by using the main idea of the proof of \cite[Theorem 1]{ismail}. Clearly $\Omega_{\mu,\nu,a,b}(x)\to 1$ as $x\to 0,$ and by using the same approach as in the proof of Theorem \ref{theolap1}, we obtain that
\begin{align*}
-&\frac{d\left[\ln\Omega_{\mu,\nu,a,b}(x)\right]}{dx}\\&=\frac{b}{2\sqrt{x}}+\sum_{n\geq1}\left(-\frac{1}{x+j^2_{\mu,n}a^{-2}}-\frac{1}{x+j^2_{\nu,n}b^{-2}}+\frac{1}{x+j^2_{\mu,n}b^{-2}}+\frac{1}{x+j^2_{\nu,n}a^{-2}}\right)\\
&= \frac{b}{2\sqrt{x}}-\sum_{n\geq1}\frac{1}{x+j^2_{\nu,n}b^{-2}}+\int_0^\infty e^{-xt}\left[\sum_{n\geq1}\left(-e^{-j^2_{\mu,n}a^{-2}t}+e^{-j^2_{\mu,n}b^{-2}t}+e^{-j^2_{\nu,n}a^{-2}t}\right)\right]dt.
\end{align*}
Note that $e^{-j^2_{\mu,n}a^{-2}t}<e^{-j^2_{\mu,n}b^{-2}t}$ for all $\mu>-1,$ $b>a,$ $t>0$ and $n\in\mathbb{N}.$ Consequently, we have that $$x\mapsto\int_{0}^{\infty}e^{-xt}\left[\sum_{n\geq1}\left(-e^{-j^2_{\mu,n}a^{-2}t}+e^{-j^2_{\mu,n}b^{-2}t}+e^{-j^2_{\nu,n}a^{-2}t}\right)\right]dt$$ is a completely monotonic function on $(0,\infty)$ for all $\mu>-1,$ $\nu>-1,$ and $b>a.$ Now, define $$\eta_{\mu,b}(x)=\frac{b}{2\sqrt{x}}-\sum_{n\geq1}\frac{1}{x+j^2_{\nu,n}b^{-2}}.$$ Observe that $$\frac{b}{2\sqrt{x}}=\frac{1}{\pi}\int_{0}^{\infty}\frac{dt}{x+b^{-2}t^2}$$ and thus we have that
$$\eta_{\mu,b}(x)= \frac{1}{\pi}\int_{0}^{\infty}\frac{dt}{x+b^{-2}t^2}-\sum_{n\geq1}\dfrac{1}{x+j^2_{\nu,n}b^{-2}}.$$
On the other hand for every $s$ such that $j_{\nu,m}\leq s<j_{\nu,m+1},$ we have
$$\int_{0}^{s}\frac{dt}{x+b^{-2}t^{2}}\geq \int_{j_{\nu,0}}^{j_{\nu,m}}\frac{dt}{x+b^{-2}t^{2}}=\sum_{n=1}^m\int_{j_{\nu,n-1}}^{j_{\nu,n}}\frac{dt}{x+b^{-2}t^{2}}\geq \sum_{n=1}^m\int_{j_{\nu,n-1}}^{j_{\nu,n}}\frac{dt}{x+b^{-2}j_{\nu,n}^{2}}=\sum_{n=1}^m\frac{j_{\nu,n}-j_{\nu,n-1}}{x+b^{-2}j^2_{\nu,n}},$$ where we used the fact that $t\mapsto 1\left/\left(x+b^{-2}t^2\right)\right.$ is a decreasing function on $(0,\infty)$ for all $b>0$ and $x>0.$ Consequently, we arrive at
$$\eta_{\mu,b}(x)\geq \frac{1}{\pi}\sum_{n\geq1}\frac{j_{\nu,n}-j_{\nu,n-1}}{x+b^{-2}j^2_{\nu,n}}-\sum_{n\geq1}\dfrac{1}{x+j^2_{\nu,n}b^{-2}}.$$
Since for $\nu>\frac{1}{2}$ we have that $j_{\nu,n}-j_{\nu,n-1}>\pi$ (see for example \cite[Theorem 1.6]{Se00}), we conclude that $\eta_{\mu,b}(x)>0$ for all $x>0$ and $b>0$. Moreover, from the previous results we obtain that
\begin{align*}
\frac{(-1)^k}{k!}\eta_{\mu,b}^{(k)}(x) &= \frac{1}{\pi}\int_{0}^{\infty}\frac{dt}{(x+b^{-2}t^{2})^{k+1}}-\sum_{n\geq1}\dfrac{1}{(x+j^2_{\nu,n}b^{-2})^{k+1}}\\
&\geq \frac{1}{\pi}\sum_{n\geq1}\frac{j_{\nu,n}-j_{\nu,n-1}}{(x+j^2_{\nu,n}b^{-2})^{k+1}}-\sum_{n\geq1}\frac{1}{(x+j^2_{\nu,n}b^{-2})^{k+1}}
\end{align*}
for all $x>0,$ $b>0$ and $k\in\mathbb{N}.$ The above inequality holds true because $t \mapsto 1\left/(x+b^{-2}t^2)^{k+1}\right.$ is a decreasing function on $(0, \infty)$ for every $x>0,$ $b>0$ and $k$ natural number. Consequently, we have that $(-1)^k\eta_{\mu,b}^{(k)}(x)\geq0$ for all $k\in\mathbb{N}$ and $x>0,$ $b>0.$ Therefore $\eta_{\mu,b}$ is also a completely monotonic function on $(0,\infty)$ and this implies that $x \mapsto -{d\left[\ln\Omega_{\mu,\nu,a,b}(x)\right]}/{dx}$ is the sum of two completely monotonic functions on $(0,\infty).$
\end{proof}

\begin{proof}[\bf Proof of Theorem \ref{thIfirst}]
Let $F(z)=e^{-a\sqrt{z}}z^{-\frac{\mu}{2}}I_{\mu}(a\sqrt{z})$. By using the relation \cite[eq. 10.30.1]{nist}, we have that
$$F(z)\sim\frac{e^{-a\sqrt{z}}a^{\mu}}{2^{\mu}\Gamma(\mu+1)}, \quad \text{as} \quad z \to 0.$$
Consequently, we arrive at $F(z)=o(|z|^{-1})$ as $|z| \to 0$. Next, by using \cite[eq. 10.30.4]{nist}, we obtain that
$$F(z)\sim\frac{1}{\sqrt{2\pi a}}z^{-\frac{(\mu+1/2)}{2}} \quad\text{as} \quad |z|\to\infty.$$
The above relations holds uniformly in every sector $|\arg z|\leq\pi-\epsilon$, $\epsilon>0$. Hence conditions {\bf a} and {\bf b} of Lemma \ref{LemStrep1} have been verified and by using the relations in \eqref{recurIKl}, we arrive at
$$F(-t-\mathrm{i}\eta)=F(e^{-\mathrm{i}\pi}(t+\mathrm{i}\eta))=e^{a\mathrm{i}\sqrt{t+\mathrm{i}\eta}}(t+\mathrm{i}\eta)^{-\frac{\mu}{2}}J_{\mu}(a\sqrt{t+\mathrm{i}\eta})$$
and
$$F(-t+\mathrm{i}\eta)=F(e^{\mathrm{i}\pi}(t-\mathrm{i}\eta))=e^{-a\mathrm{i}\sqrt{t-\mathrm{i}\eta}}(t-\mathrm{i}\eta)^{-\frac{\mu}{2}}J_{\mu}(a\sqrt{t-\mathrm{i}\eta}).$$
By using the asymptotic relation (see \cite[eq. 10.7.3]{nist}) as $z\to 0$
\begin{equation}\label{a1}
J_{\nu}(z) \sim \dfrac{\left(z/2\right)^{\nu}}{\Gamma(\nu+1)},\quad \nu\neq-1, -2, \dots,
\end{equation}
we conclude that the function $F(-t-\mathrm{i}\eta)-F(-t+\mathrm{i}\eta)$ is continuous in every rectangle $[t_1, t_2]\times[0, \eta],$ where $t_1,$ $t_2,$ $\eta>0.$ Thus, the limit and the integral in equation \eqref{inv1} can be interchanged and by using the fundamental theorem of calculus, we obtain that
$$\alpha'(t)=\frac{1}{\pi}t^{-\frac{\mu}{2}}J_{\mu}(a\sqrt{t})\sin(a\sqrt{t})$$
and we can observe that $\alpha(t)$ is continuous and $\lim_{t\to 0^{+}}\alpha(t)$ exists. Let $\tilde{\alpha}(t)=\alpha(t)-\alpha(0)$. Thus we get the normalized measure $\tilde{\alpha}(t)$ with $\tilde{\alpha}'(t)=\alpha'(t)$.
\end{proof}

\begin{proof}[\bf Proof of Theorem \ref{thmprod00}]
Let us consider the following function $$F(z)=z^{\frac{\nu-\mu}{2}}I_{\mu}(a\sqrt{z})K_{\nu}(b\sqrt{z}).$$ By using the relations \cite[eq. 10.30.1]{nist} and \cite[eq. 10.30.2]{nist} for $\nu>0$ and $\mu>-1$ we have that as $z\to0$
$$F(z)\sim\left(\frac{1}{2}\right)^{\mu-\nu+1}\frac{\Gamma(\nu)}{\Gamma(\mu+1)}\frac{a^{\mu}}{b^{\nu}}.$$
Moreover, by using the relations \cite[eq. 10.30.1]{nist} and \cite[eq. 10.30.3]{nist}, we obtain that for $\nu=0$ and $\mu>-1$ as $z\to0$
$$F(z)\sim -\frac{a^{\mu}}{2^{\mu}\Gamma(\mu+1)}\ln(b\sqrt{z}).$$
Similarly, by using the relations $K_{\nu}(z)=K_{-\nu}(z)$ (see \cite[eq. 10.27.3]{nist}), \cite[eq. 10.30.2]{nist} and \cite[eq. 10.30.1]{nist}, for $-1<\nu<0$ and $\mu>-1$ we have that as $z\to0$
$$F(z)\sim \frac{a^{\mu}b^{\nu}}{2^{\mu+\nu+1}}\frac{\Gamma(-\nu)}{\Gamma(\mu+1)}z^{\nu+1}.$$
Consequently, we obtain that $F(z)=o\left(|z|^{-1}\right)$ as $z\to0$ for all $\nu>-1$ and $\mu>-1.$ Now, by using the relations \cite[eq. 10.30.4]{nist} and \cite[eq. 10.25.3]{nist}, we have that as $|z|\to \infty$
$$F(z)\sim z^{\frac{\nu-\mu-1}{2}}\dfrac{e^{-(b-a)\sqrt{z}}}{2}$$
holds uniformly in every sector $\left|\arg z\right|\leq\pi-\varepsilon$, $\varepsilon>0$. Consequently, we clearly obtain that $F(z)=o(1)$ as $|z| \to \infty,$ whenever $b\geq a$ and $\nu-\mu<1.$ Thus, the conditions in Lemma \ref{LemStrep1} has been verified because $F(z)$ is analytic in $\left|\arg z\right|<\pi.$ Now, in order to apply Lemma \ref{LemStinv}, we observe that
\begin{align*}
F(-t-\mathrm{i}\eta)=F(e^{-\mathrm{i}\pi}(t+\mathrm{i}\eta))&=\left(e^{-\mathrm{i}\pi}(t+\mathrm{i}\eta)\right)^{\frac{\nu-\mu}{2}}I_{\mu}(ae^{-\frac{\mathrm{i}\pi}{2}}\sqrt{t+\mathrm{i}\eta})
K_{\nu}(be^{-\frac{\mathrm{i}\pi}{2}}\sqrt{t+\mathrm{i}\eta})\\
&=\frac{\mathrm{i}\pi}{2}(t+\mathrm{i}\eta)^{\frac{\nu-\mu}{2}}J_{\mu}(a\sqrt{t+\mathrm{i}\eta})H^{(1)}_{\nu}(b\sqrt{t+\mathrm{i}\eta})
\end{align*}
and
\begin{align*}
F(-t+\mathrm{i}\eta)=F(e^{\mathrm{i}\pi}(t-\mathrm{i}\eta))&=\left(e^{\mathrm{i}\pi}(t-\mathrm{i}\eta)\right)^{\frac{\nu-\mu}{2}}I_{\mu}(ae^{\frac{\mathrm{i}\pi}{2}}\sqrt{t-\mathrm{i}\eta})
K_{\nu}(be^{\frac{\mathrm{i}\pi}{2}}\sqrt{t-\mathrm{i}\eta})\\
&=-\frac{\mathrm{i}\pi}{2}(t-\mathrm{i}\eta)^{\frac{\nu-\mu}{2}}J_{\mu}(a\sqrt{t-\mathrm{i}\eta})H^{(2)}_{\nu}(b\sqrt{t-\mathrm{i}\eta}),
\end{align*}
where we have used the following relations (see \cite[p. 459]{may} and \cite[10.27]{nist})
\begin{equation}\label{recurIKl}
\begin{cases}
H^{(1)}_{\mu}(z) = J_{\mu}(z) + \mathrm{i}Y_{\mu}(z),& H^{(2)}_{\mu}(z) = J_{\mu}(z) - \mathrm{i}Y_{\mu}(z),\\
K_{\mu}(ze^{-\frac{1}{2}\mathrm{i}\pi }) = \frac{1}{2}\mathrm{i}\pi e^{\frac{1}{2}\mathrm{i}\pi\mu}H_{\mu}^{(1)}(z),& K_{\mu}(ze^{\frac{1}{2}\mathrm{i}\pi }) = -\frac{1}{2}\mathrm{i}\pi e^{-\frac{1}{2}\mathrm{i}\pi\mu}H_{\mu}^{(2)}(z),\\
I_{\mu}(ze^{\pm\frac{1}{2}\mathrm{i}\pi}) = e^{\pm\frac{1}{2}\mathrm{i}\mu\pi}J_{\mu}(z),& I_{\mu}(z)=e^{-\frac{1}{2}\mathrm{i}\mu \pi}J_{\mu}(ze^{\frac{1}{2}\mathrm{i}\pi}).
\end{cases}
\end{equation}

Now, we are going to show the justification for interchanging the limit and the integral in equation \eqref{inv1}, in three cases:

\textbf{Case (i) $\nu>0$ and $\mu>-1.$} By using the asymptotic relation (refer \cite[eq. 10.7.7]{nist})
\begin{equation}\label{a2}
H^{(1)}_{\nu}(z) \sim -H^{(2)}_{\nu}(z) \sim -\frac{\mathrm{i}}{\pi}\Gamma(\nu)\left(z/2\right)^{-\nu}, \quad z \to 0,\quad \real \nu>0
\end{equation}
and relation \eqref{a1}, we observe that the function $F(-t-\mathrm{i}\eta)-F(-t+\mathrm{i}\eta)$ is continuous in every rectangle $[t_1, t_2]\times[0, \eta]$, where $t_1$, $t_2$, $\eta>0$. Thus, in this case, we can interchange the limit and the integral in equation \eqref{inv1}.

\textbf{Case (ii) $\nu=0$ and $\mu>-1.$} By using the asymptotic relation (see \cite[eq. 10.7.2]{nist})
\begin{equation}\label{a3}
H^{(1)}_{0}(z)\sim -H^{(2)}_{0}(z) \sim \frac{2\mathrm{i}}{\pi}\ln(z), \quad z \to 0
\end{equation}
and relation \eqref{a1}, we obtain that
$$F(-t-\mathrm{i}\eta) \sim -\frac{a^{\mu}}{2^{\mu}\Gamma(\mu+1)}\ln\left(b\sqrt{t+\mathrm{i}\eta}\right),\quad t\to 0, \eta \to 0$$ and
$$F(-t+\mathrm{i}\eta) \sim -\frac{a^{\mu}}{2^{\mu}\Gamma(\mu+1)}\ln\left(b\sqrt{t-\mathrm{i}\eta}\right),\quad t\to 0, \eta \to 0.$$
By using the fact that $f(z)\sim g(z)$ as $z \to 0,$ implies $f(z)=\mathcal{O}(g(z))$ as $z \to 0,$ in view of the definition of the Landau symbol big-oh, there exist $c_1$, $\alpha_1$, $\beta_1 >0,$ such that $$\left|F(-t-\mathrm{i}\eta)\right|\leq \frac{c_1a^{\mu}}{2^{\mu}\Gamma(\mu+1)}\left|\ln\left(b\sqrt{t+\mathrm{i}\eta}\right)\right|, \quad 0<t<\alpha_1, 0<\eta<\beta_1.$$ Similarly, there exist $c_1$, $\alpha_2$, $\beta_2 >0$, such that $$\left|F(-t+\mathrm{i}\eta)\right|\leq\frac{c_2a^{\mu}}{2^{\mu}\Gamma(\mu+1)}\left|\ln\left(b\sqrt{t-\mathrm{i}\eta}\right)\right|, \quad 0<t<\alpha_2, 0<\eta<\beta_2.$$ We note that the right hand side of both inequalities above are integrable on the interval $[t_1, t_2]$ where $0<t_1<t_2<\alpha$ and $\alpha=\min\{\alpha_1, \alpha_2\}$. Thus, we can interchange again the limit and the integral in equation \eqref{inv1}.

\textbf{Case (iii) $-1<\nu<0$ and $\mu>-1.$} By using the following relations (see \cite[eq. 10.4.6]{nist})
\begin{align}
H^{(1)}_{-\nu}(z)&=e^{\nu\pi \mathrm{i}}H^{(1)}_{\nu}(z)\label{a4}\\
H^{(2)}_{-\nu}(z)&=e^{-\nu\pi \mathrm{i}}H^{(2)}_{\nu}(z)\label{a5}
\end{align}
and the asymptotic relations \eqref{a1} and \eqref{a2}, we have that
$$F(-t-\mathrm{i}\eta) \sim \frac{e^{-\nu\pi \mathrm{i}}\Gamma(-\nu)a^{\mu}b^{\nu}}{2^{\mu+\nu+1}\Gamma(\mu+1)}(t+\mathrm{i}\eta)^{\nu}, \quad t\to 0, \eta \to 0$$ and
$$F(-t+\mathrm{i}\eta) \sim \frac{e^{\nu\pi \mathrm{i}}\Gamma(-\nu)a^{\mu}b^{\nu}}{2^{\mu+\nu+1}\Gamma(\mu+1)}(t-\mathrm{i}\eta)^{\nu}, \quad t\to 0, \eta \to 0.$$
In view of the definition of the Landau symbol big-oh, there exist $c_1>0$, $\alpha_1>0$ and $\beta_1>0$ such that $0<t<\alpha_1$ and $0<\eta<\beta_1$, and we have that
$$\left|F(-t-\mathrm{i}\eta)\right| \leq \frac{c_1\Gamma(-\nu)a^{\mu}b^{\nu}}{2^{\mu+\nu+1}\Gamma(\mu+1)}\left|(t+\mathrm{i}\eta)^{\nu}\right|.$$
Similarly, there exist $c_2>0$, $\alpha_2>0$ and $\beta_2>0$ such that $0<t<\alpha_2$ and $0<\eta<\beta_2$ and we arrive at
$$\left| F(-t+\mathrm{i}\eta)\right| \leq \dfrac{c_2\Gamma(-\nu)a^{\mu}b^{\nu}}{2^{\mu+\nu+1}\Gamma(\mu+1)}\left|(t-\mathrm{i}\eta)^{\nu}\right|.$$
We note that the right hand side of the above inequalities are integrable on the interval $[t_1, t_2],$ where $0<t_1<t_2<\alpha$ and $\alpha=\min\{\alpha_1, \alpha_2\}$. Thus, in this case we can also interchange the limit and the integral in equation \eqref{inv1}.

Consequently, we have that
$$\alpha'(t)=\dfrac{1}{2\pi \mathrm{i}}\lim_{\eta \to 0^{+}}\left[ F(-t-\mathrm{i}\eta)-F(-t+\mathrm{i}\eta)\right]=\frac{t^{\frac{\nu-\mu}{2}}}{2}J_{\mu}(a\sqrt{t})J_{\nu}(b\sqrt{t}).$$
From the above equation we can see that $\alpha(t)$ is continuous and $\lim_{t \to 0}\alpha(t)$ exists. Let $\tilde{\alpha}(t)=\alpha(t)-\alpha(0)$. Note that $\tilde{\alpha}(t)$ is the normalized measure with $\tilde{\alpha}'(t)=\alpha'(t)$.
\end{proof}

\begin{proof}[\bf Proof of Theorem \ref{thprod1}]
We know that (see for example \cite[p. 8]{bondesson}) the function $x\mapsto \phi(x)$ is the Laplace transform of a probability distribution on $(0,\infty)$ if and only if $\phi(0)=1$ and $\phi(x)$ has the form $$\phi(x)=\int_{0}^{\infty}e^{-xt}d\mu(t)$$ and is completely monotonic on $(0,\infty).$ In view of \eqref{eqprod1} this immediately implies that the function $x\mapsto 2\mu I_{\mu}(\sqrt{x})K_{\mu}(\sqrt{x})$ is the Laplace transform of a probability distribution with support $(0,\infty).$ Combining this with the second equality in \eqref{eqprod1} we clearly have that the probability density function of this distribution, that is,
$$\varsigma_{\mu}(x)=\mu\int_{0}^{\infty} e^{-tx}J_{\mu}^2(\sqrt{t})dt$$ is a completely monotonic function on $(0,\infty)$ for all $\mu>0,$ and by using the Goldie-Steutel law we obtain that indeed the function $x\mapsto 2\mu I_{\mu}(\sqrt{x})K_{\mu}(\sqrt{x})$ is the Laplace transform of an infinitely divisible probability distribution with support $(0,\infty).$ Moreover, by using the formula \cite[p. 139]{OB12}
$$\int_{0}^{\infty}e^{-st}J_{\mu}^2(\sqrt{t})dt=\frac{1}{s}e^{-\frac{1}{2s}}I_{\mu}\left(\frac{1}{2s}\right),$$
where $\mu>0,$ we observe that
$$\varsigma_{\mu}(x)=\frac{\mu}{x}e^{-\frac{1}{2x}}I_{\mu}\left(\frac{1}{2x}\right),$$
which completes the proof of the derivation of the probability density function.
\end{proof}

\begin{proof}[\bf Proof of Corollary \ref{corprodcm}]
The assertions in parts {\bf a}, {\bf b} and {\bf c} follow immediately from the integral representation \eqref{eqprod1}. The proof for part {\bf d} is as follows. By using the relation \cite[p. 8]{segura}, we have that
$$[xI_{\mu}(x)K_{\mu}(x)]^{-1}=\frac{I_{\mu-1}(x)}{I_{\mu}(x)}+\frac{K_{\mu-1}(x)}{K_{\mu}(x)}$$
and in view of the three-term recurrence relation $$I_{\mu-1}(x)=I_{\mu+1}(x)+\frac{2\mu}{x}I_{\mu}(x),$$
we arrive at
$$[xI_{\mu}(x)K_{\mu}(x)]^{-1}=\frac{2\mu}{x}+\frac{I_{\mu+1}(x)}{I_{\mu}(x)}+\frac{K_{\mu-1}(x)}{K_{\mu}(x)}.$$
Now, replacing $x$ by $\sqrt{x}$ and multiplying with $\frac{1}{\sqrt{x}}$ both sides of the above equation, we obtain that
$$[xI_{\mu}(\sqrt{x})K_{\mu}(\sqrt{x})]^{-1}=\frac{2\mu}{x}+\frac{1}{\sqrt{x}}\frac{I_{\mu+1}(\sqrt{x})}{I_{\mu}(\sqrt{x})}+\frac{1}{\sqrt{x}}\frac{K_{\mu-1}(\sqrt{x})}{K_{\mu}(\sqrt{x})}.$$
Recall the well-known integral representation (see for example \cite[eq. (1.3)]{ismail00} or \cite[eq. (2.4)]{miller})
\begin{equation}\label{integralKquot}
\frac{K_{\mu-1}(\sqrt{x})}{\sqrt{x}K_{\mu}(\sqrt{x})}=\frac{2}{\pi^2}\int_{0}^{\infty}\frac{t^{-1}}{x+t}\left[J^2_{\mu}(\sqrt{t})+Y^2_{\mu}(\sqrt{t})\right]^{-1}dt,
\end{equation}
where $x>0,$ $\mu\geq0.$ Thus, by using the series and integral representations for the ratios involving modified Bessel functions of the first and second kinds, respectively, we have that
$$[xI_{\mu}(\sqrt{x})K_{\mu}(\sqrt{x})]^{-1}=\dfrac{2\mu}{x}+\sum_{n\geq1}\frac{2}{x+j_{\mu,n}^2}+\frac{2}{\pi^2}\int_{0}^{\infty}\frac{t^{-1}}{x+t}[J^2_{\mu}(\sqrt{t})+Y^2_{\mu}(\sqrt{t})]^{-1} dt,$$
which is the sum of three completely monotonic functions on $(0,\infty)$ for all $\mu>0.$

Finally, for the proof of part {\bf e} we use the upper bound of Landau (see \cite{La00} or \cite{BMPS16})
$$|J_{\mu}(t)|\leq c_{L}|t|^{-\frac{1}{3}}, \quad \mu>0, \quad t\in \mathbb{R}, \quad c_{L}=0.7857468704\dots$$
and the integral representation of $I_{\mu}(\sqrt{x})K_{\mu}(\sqrt{x}).$ We clearly obtain that
$$I_{\mu}(\sqrt{x})K_{\mu}(\sqrt{x})\leq\frac{c^2_{L}}{2}\int_{0}^{\infty}\frac{t^{-\frac{1}{3}}}{x+t} dt=\frac{\pi c^2_{L}}{\sqrt{3}x^{\frac{1}{3}}},$$
which completes the proof.
\end{proof}

\begin{proof}[\bf Proof of Theorem \ref{theprodIKexprepr2}]
Let $$F(z)=z^{\frac{\nu-\mu}{2}}e^{-a\sqrt{z}}I_{\mu}(a\sqrt{z})K_{\nu}(b\sqrt{z}).$$ By using the relations \cite[eq. 10.30.1]{nist} and \cite[eq. 10.30.2]{nist} for the case $\mu>-1$ and $\nu>0$, we have that
$$F(z) \sim {\frac{e^{-a\sqrt{z}}\Gamma(\nu)a^{\mu}}{2^{\mu-\nu+1}\Gamma(\mu+1)b^{\nu}}}\quad \text{as} \quad z \to 0.$$
In the case when $\mu>-1$ and $\nu=0$, by using the relations \cite[eq. 10.30.1]{nist} and \cite[eq. 10.30.3]{nist}, we arrive at
$$F(z) \sim -\frac{e^{-a\sqrt{z}}a^{\mu}\ln(b\sqrt{z})}{2^{\mu}\Gamma(\mu+1)} \quad \text{as} \quad z \to 0,$$
while in the case when $\mu>-1$ and $-1<\nu<0$, by using the relations $K_{\nu}(z)=K_{-\nu}(z)$ (see \cite[eq. 10.27.3]{nist}), \cite[eq. 10.30.1]{nist}, and \cite[eq. 10.30.2]{nist}, we arrive at
$$F(z) \sim {\frac{z^{\nu}e^{-a\sqrt{z}}\Gamma(-\nu)a^{\mu}b^{\nu}}{2^{\mu+\nu+1}\Gamma(\mu+1)}} \quad \text{as} \quad z \to 0.$$
Consequently, we conclude that $F(z)=o(|z|^{-1})$ as $z \to 0$. On the other hand, by using the relations \cite[eq. 10.30.5]{nist} and \cite[eq. 10.25.3]{nist}, we arrive at
$$F(z) \sim \frac{z^{\frac{\nu-\mu-1}{2}}e^{-b\sqrt{z}}}{2\sqrt{ab}} \quad \text{as} \quad |z| \to \infty.$$
The above relation holds uniformly in every sector $|\arg z|\leq\pi-\epsilon$, $\epsilon>0$ and thus we have that $F(z)=o(1)$ as $|z| \to \infty$. Hence the conditions {\bf a} and {\bf b} in Lemma \ref{LemStrep1} have been verified, and this implies that the corresponding Stieltjes transform representation exists. All we need is just to use the inversion theorem to obtain the corresponding measure. In view of Lemma \ref{LemStinv} and equations \eqref{recurIKl} we obtain that
$$F(-t-\mathrm{i}\eta)=F(e^{-\mathrm{i}\pi}(t+\mathrm{i}\eta))=\frac{\mathrm{i}\pi}{2}(t+\mathrm{i}\eta)^{\frac{\nu-\mu}{2}}e^{a\mathrm{i}\sqrt{t+\mathrm{i}\eta}}J_{\mu}(a\sqrt{t+\mathrm{i}\eta})H^{(1)}_{\nu}(b\sqrt{t+\mathrm{i}\eta})$$ and
$$F(-t+\mathrm{i}\eta)=F(e^{\mathrm{i}\pi}(t-\mathrm{i}\eta))=-\frac{\mathrm{i}\pi}{2}(t-\mathrm{i}\eta)^{\frac{\nu-\mu}{2}}e^{-a\mathrm{i}\sqrt{t-\mathrm{i}\eta}}J_{\mu}(a\sqrt{t-\mathrm{i}\eta})H^{(2)}_{\nu}(b\sqrt{t-\mathrm{i}\eta}).$$
The justification for interchanging the limit and the integral in equation \eqref{inv1} goes similarly as in the lines introduced in the proof of Theorem \ref{thmprod00}, and hence it is omitted here. Finally, observe that by using the relations in \eqref{recurIKl}, we obtain that
$$\alpha'(t)=\frac{t^{\frac{\nu-\mu}{2}}}{2}J_{\mu}(a\sqrt{t})\left[J_{\nu}(b\sqrt{t})\cos(a\sqrt{t})-Y_{\nu}(b\sqrt{t})\sin(a\sqrt{t})\right].$$
\end{proof}

\begin{proof}[\bf Proof of Theorem \ref{theprodIKexp}]
Observe that in view of the recurrence relation \cite[p. 79]{watson}
$$K_{\mu}'(x)=-K_{\mu-1}(x)-\frac{\mu}{x}K_{\mu}(x)$$
and the equivalent form of the integral representation \eqref{integralKquot}, that is,
$$\frac{K_{\mu-1}(\sqrt{x})}{\sqrt{x}K_{\mu}(\sqrt{x})}=\frac{4}{\pi^2}\int_{0}^{\infty}\frac{t^{-1}}{x+t^2}\left[J^2_{\mu}(t)+Y^2_{\mu}(t)\right]^{-1}dt$$
where $x>0,$ $\mu\geq0,$ we obtain that
\begin{align*}
-\frac{d\ln\chi_{\mu,\nu,a,b}(x)}{dx}&=\frac{a}{2\sqrt{x}}+\frac{\mu-\nu}{2x}-\frac{a}{2\sqrt{x}}\frac{I_{\mu}'(a\sqrt{x})}{I_{\mu}(a\sqrt{x})}-\frac{b}{2\sqrt{x}}\frac{K_{\nu}'(b\sqrt{x})}{K_{\nu}(b\sqrt{x})}\\
&=\frac{a}{2\sqrt{x}}-\frac{a}{2\sqrt{x}}\dfrac{I_{\mu+1}(a\sqrt{x})}{I_{\mu}(a\sqrt{x})}+\frac{b}{2\sqrt{x}}\frac{K_{\nu-1}(b\sqrt{x})}{K_{\nu}(b\sqrt{x})}\\
&=\frac{1}{\pi}\int_{0}^{\infty}\frac{dt}{x+a^{-2}t^2}-\sum_{n\geq1}\dfrac{1}{x+a^{-2}j^2_{\mu,n}}+\dfrac{2}{\pi^2}\int_{0}^{\infty}\dfrac{t^{-1}}{x+b^{-2}t^2}[J^2_{\nu}(t)+Y^2_{\nu}(t)]^{-1}dt\\
&=\Theta(x)+\dfrac{2}{\pi^2}\int_{0}^{\infty}\dfrac{t^{-1}}{x+b^{-2}t^2}[J^2_{\nu}(t)+Y^2_{\nu}(t)]^{-1}dt,
\end{align*}
where 	
$$\Theta(x)=\frac{1}{\pi}\int_{0}^{\infty}\frac{dt}{x+a^{-2}t^2}-\sum_{n\geq1}\dfrac{1}{x+a^{-2}j^2_{\mu,n}}\geq\frac{1}{\pi}\sum_{n\geq1}
\frac{j_{\mu,n}-j_{\mu,n-1}}{x+a^{-2}j^2_{\mu,n}}-\sum_{n\geq1}\frac{1}{x+a^{-2}j^2_{\mu,n}}\geq 0$$
and
$$\frac{(-1)^{m}}{m!}\frac{d^m\Theta(x)}{dx^m}\geq\frac{1}{\pi}\sum_{n\geq1}\frac{j_{\mu,n}-j_{\mu,n-1}}{(x+a^{-2}j^2_{\mu,n})^{m+1}}-\sum_{n\geq1}\frac{1}{(x+a^{-2}j^2_{\mu,n})^{m+1}}\geq 0$$
since $j_{\mu,n}-j_{\mu,n-1}>\pi$ whenever $\mu>\frac{1}{2}$ and $n\in\mathbb{N}.$ Thus, $x\mapsto -{d\ln\chi_{\mu,\nu,a,b}(x)}/{dx}$ is completely monotonic on $(0,\infty)$ for all $a,b,\nu>0$ and $\mu>\frac{1}{2}$ and consequently in view of Lemma \ref{lemma1} the proof is complete.
\end{proof}

\begin{proof}[\bf Proof of Theorem \ref{thprodK1}]
Let ${F(z)=z^{\frac{\mu+\nu}{2}}K_{\mu}(a\sqrt{z})K_{\nu}(b\sqrt{z})}$. By using the well-known asymptotic relation \cite[eq. 10.30.2]{nist} for $a,b,\mu,\nu>0$, we have that
$$z^{\frac{\mu+\nu}{2}}K_{\mu}(a\sqrt{z})K_{\nu}(b\sqrt{z})\sim\frac{1}{4}\Gamma(\mu)\Gamma(\nu)\frac{2^{\mu+\nu}}{a^{\mu}b^{\nu}}$$
as $z\to 0.$ Consequently, we have that $f(z)=o(|z|^{-1})$ as $|z| \to 0$. Now, by using the asymptotic relation \cite[eq. 10.25.3]{nist}, we have
$$z^{\frac{\mu+\nu}{2}}K_{\mu}(a\sqrt{z})K_{\nu}(b\sqrt{z})\sim\frac{\pi}{2\sqrt{ab}}\sqrt{z}^{\mu+\nu-1}e^{-(a+b)\sqrt{z}}$$
as $|z|\to \infty$ and $\left|\arg z\right|<\frac{3\pi}{2}.$ From the above relation, we obtain that $F(z)=o(1)$ as $|z|\to \infty$, and this holds true uniformly in every sector $\left|\arg z\right|\leq\pi-\varepsilon$, $\varepsilon>0$. Therefore, the conditions of Lemma \ref{LemStrep1} hold true. On the other hand, by using \eqref{recurIKl} we obtain that
$$F(-t-\mathrm{i}\eta)=F(e^{-\mathrm{i}\pi}(t+\mathrm{i}\eta))=-\frac{\pi^2}{4}(t+\mathrm{i}\eta)^{\frac{\mu+\nu}{2}}H^{(1)}_{\mu}(a\sqrt{t+\mathrm{i}\eta})H^{(1)}_{\nu}(b\sqrt{t+\mathrm{i}\eta})$$
and
$$F(-t+\mathrm{i}\eta)=F(e^{\mathrm{i}\pi}(t-\mathrm{i}\eta))=-\frac{\pi^2}{4}(t-\mathrm{i}\eta)^{\frac{\mu+\nu}{2}}H^{(2)}_{\mu}(a\sqrt{t-\mathrm{i}\eta})H^{(2)}_{\nu}(b\sqrt{t-\mathrm{i}\eta}),$$
By using the asymptotic relation \eqref{a2}, we conclude that the function $F(-t-\mathrm{i}\eta)-F(-t+\mathrm{i}\eta)$ is continuous in every rectangle $[t_1, t_2]\times[0, \eta],$ $t_1, t_2, \eta >0.$ Consequently, we can interchange the limit and the integral in equation \eqref{inv1}. In view of Lemma \ref{LemStinv} we arrive at
$$\alpha'(t)=-\frac{\pi}{4}t^{\frac{\mu+\nu}{2}}[J_{\mu}(a\sqrt{t})Y_{\nu}(b\sqrt{t})+J_{\nu}(b\sqrt{t})Y_{\mu}(a\sqrt{t})].$$
From the above equation we can observe that $\alpha(t)$ is continuous and $\lim_{t \to 0^{+}}\alpha(t)$ exists. Let $\tilde{\alpha}(t)=\alpha(t)-\alpha(0)$. Note that $\tilde{\alpha}(t)$ is the normalized measure with $\tilde{\alpha}'(t)=\alpha'(t)$.
\end{proof}

\begin{proof}[\bf Proof of Theorem \ref{thinfdivprodK}]
The infinite divisibility follows naturally by using the result of Ismail and Kelker \cite[Theorem 1.8]{kelker}. Namely, we know that the function $x\mapsto (\sqrt{x})^{\nu}K_{\nu}(\sqrt{x})/\left[2^{\nu-1}\Gamma(\nu)\right]$ is the Laplace transform of an infinitely divisible distribution with support $(0,\infty)$ whenever $\nu>0.$ This implies that the Laplace transform $\vartheta_{\mu,\nu,a,b}(x)$ is in fact the product of two Laplace transforms of two infinitely divisible distributions, and thus it is the Laplace transform of an infinitely divisible distribution, according to \cite[Proposition 2.1]{steutel}. More precisely, in view of \eqref{integralKquot} we obtain that
\begin{align*}
-\frac{d \ln \vartheta_{\mu,\nu,a,b}(x)}{dx}&= -\frac{\mu+\nu}{2x}-\frac{a}{2\sqrt{x}}\frac{K'_{\mu}(a\sqrt{x})}{K_{\mu}(a\sqrt{x})}-\frac{b}{2\sqrt{x}}\frac{K'_{\nu}(b\sqrt{x})}{K_{\nu}(b\sqrt{x})}\\
&=\frac{a}{2\sqrt{x}}\frac{K_{\mu-1}(a\sqrt{x})}{K_{\mu}(a\sqrt{x})}+\frac{b}{2\sqrt{x}}\frac{K_{\nu-1}(b\sqrt{x})}{K_{\nu}(b\sqrt{x})}\\
&=\frac{1}{\pi^2}\left[\int_{0}^{\infty}\frac{t^{-1}}{x+ta^{-2}}\left[J^2_{\mu}(\sqrt{t})+Y^2_{\mu}(\sqrt{t})\right]^{-1}dt + \int_{0}^{\infty} \frac{t^{-1}}{x+tb^{-2}}\left[J^2_{\nu}(\sqrt{t})+Y^2_{\nu}(\sqrt{t})\right]^{-1}dt \right],
\end{align*}
which as a function $x$ is clearly completely monotonic on $(0, \infty)$ for all $a,$ $b,$ $\mu$ and $\nu$ strictly positive real numbers. Applying Lemma \ref{lemma1}, this shows that indeed $x\mapsto \vartheta_{\mu,\nu,a,b}(x)$ is the Laplace transform of an infinitely divisible distribution.

Now, observe that for $\xi(x)=\vartheta_{\mu,\nu,a,b}(x)/\vartheta_{\mu,\nu,a,b}(\alpha x)$ we have that
\begin{align*}
-\frac{d \ln \xi(x)}{dx}&= -\frac{a}{2\sqrt{x}}\frac{K'_{\mu}(a\sqrt{x})}{K_{\mu}(a\sqrt{x})}-\frac{b}{2\sqrt{x}}\frac{K'_{\nu}(b\sqrt{x})}{K_{\nu}(b\sqrt{x})}
+\frac{a\sqrt{\alpha}}{2\sqrt{x}}\frac{K'_{\mu}(a\sqrt{\alpha x})}{K_{\mu}(a\sqrt{\alpha x})}+\frac{b\sqrt{\alpha}}{2\sqrt{x}}\frac{K'_{\nu}(b\sqrt{\alpha x})}{K_{\nu}(b\sqrt{\alpha x})}\\
&=\frac{a}{2\sqrt{x}}\frac{K_{\mu-1}(a\sqrt{x})}{K_{\mu}(a\sqrt{x})}+\frac{b}{2\sqrt{x}}\frac{K_{\nu-1}(b\sqrt{x})}{K_{\nu}(b\sqrt{x})}
-\frac{a\sqrt{\alpha}}{2\sqrt{x}}\frac{K_{\mu-1}(a\sqrt{\alpha x})}{K_{\mu}(a\sqrt{\alpha x})}-\frac{b\sqrt{\alpha}}{2\sqrt{x}}\frac{K_{\nu-1}(b\sqrt{\alpha x})}{K_{\nu}(b\sqrt{\alpha x})}\\
&=\frac{a^2}{\pi^2}\int_{0}^{\infty}\frac{(1-\alpha)}{(a^2\alpha x+t)(a^2x+t)}\left[J^2_{\mu}(\sqrt{t})+Y^2_{\mu}(\sqrt{t})\right]^{-1}dt\\&\quad + \frac{b^2}{\pi^2}\int_{0}^{\infty}\frac{(1-\alpha)}{(b^2\alpha x+t)(b^2x+t)}\left[J^2_{\nu}(\sqrt{t})+Y^2_{\nu}(\sqrt{t})\right]^{-1}dt
\end{align*}
is completely monotonic with respect to $x$ on $(0,\infty)$ for all $\alpha\in(0,1)$ and $a,$ $b,$ $\mu$ and $\nu$ strictly positive real numbers. Thus, applying Lemma \ref{lemma1} and Lemma \ref{lemma2} we conclude that $\xi$ is the Laplace transform of an infinitely divisible distribution, and consequently indeed $x\mapsto \vartheta_{\mu,\nu,a,b}(x)$ is the Laplace transform of a self-decomposable distribution.

Finally, let $w=v+{1}/{v}$. We are going to prove that $w \mapsto \vartheta_{\mu,\nu,a,b}(uv)\vartheta_{\mu,\nu,a,b}\left({u}/{v}\right)$ is completely monotonic on $(0,\infty)$ for all $a, b, \mu, \nu$ strictly positive real numbers. For this observe that
$$\vartheta_{\mu,\nu,a,b}(uv)\vartheta_{\mu,\nu,a,b}\left(\frac{u}{v}\right) =\alpha_{\mu,\nu,a,b}^2\cdot u^{\mu+\nu} K_\mu\left(a\sqrt{uv}\right) K_\mu\left(a\sqrt{\frac{u}{v}}\right) K_\nu\left(b\sqrt{uv}\right) K_\nu\left(b\sqrt{\frac{u}{v}}\right),$$
where
$$\alpha_{\mu,\nu,a,b}=\frac{a^{\mu}b^{\nu}}{2^{\mu+\nu-2}\Gamma(\mu)\Gamma(\nu)},$$
and by using McDonald's integral representation for modified Bessel functions of the second kind \eqref{prodK} we obtain that the expression $\vartheta_{\mu,\nu,a,b}(uv)\vartheta_{\mu,\nu,a,b}\left({u}/{v}\right)$ can be written as
$$\frac{\alpha_{\mu,\nu,a,b}^2}{4}u^{\mu+\nu}\int_{0}^{\infty}\int_{0}^{\infty}\exp\left(-wu\left(\frac{b^2t+a^2s}{2ts}\right)\right)
\exp\left(-\frac{t+s}{2}\right)K_\mu\left(\frac{au^2}{t}\right)K_\nu\left(\frac{bu^2}{s}\right)\frac{dt}{t}\frac{ds}{s}.$$
Now, observe that the factor which contains $w$ makes the above function completely monotonic with respect to $w,$ and consequently in view of Lemma \ref{lemma4} we obtain that $x\mapsto \vartheta_{\mu,\nu,a,b}(x)$ is the Laplace transform of a generalized gamma convolution.
\end{proof}

\begin{proof}[\bf Proof of Theorem \ref{thprodeqI}]
Let $F(z)=e^{-(a+b)\sqrt{z}}z^{-\frac{\mu+\nu}{2}}I_{\mu}(a\sqrt{z})I_{\nu}(b\sqrt{z})$. Note that as $z \to 0$
\begin{equation*}
z^{-\frac{\mu+\nu}{2}}I_{\mu}(a\sqrt{z})I_{\nu}(b\sqrt{z})\sim \left(\frac{1}{2}\right)^{\mu+\nu}\dfrac{a^{\mu}b^{\nu}}{\Gamma(\mu+1)\Gamma(\nu+1)}
\end{equation*}
and consequently $F(z)=o(|z|^{-1})$ as $|z|\to 0$. By using \cite[eq. 10.25.3]{nist}, we obtain that $|z|\to\infty$
\begin{equation*}
F(z)\sim \frac{1}{2\pi\sqrt{ab}}z^{-\frac{\mu+\nu+1}{2}},
\end{equation*}
and the above relation is valid uniformly in every sector $|\arg z|\leq\pi-\varepsilon$, $\varepsilon>0.$ Thus, the conditions {\bf a} and {\bf b} of Lemma \ref{LemStrep1} have been verified. Now, we evaluate the function values $F(-t-\mathrm{i}\eta)$ and $F(-t+\mathrm{i}\eta)$ by using the relations in \eqref{recurIKl}, as follows
$$F(-t-\mathrm{i}\eta)=F(e^{-\mathrm{i}\pi}(t+\mathrm{i}\eta))=e^{(a+b)\mathrm{i}\sqrt{t+\mathrm{i}\eta}}(t+\mathrm{i}\eta)^{-\frac{\mu+\nu}{2}}J_{\mu}(a\sqrt{t+\mathrm{i}\eta})J_{\nu}(b\sqrt{t+\mathrm{i}\eta})$$
and
$$F(-t+\mathrm{i}\eta)=F(e^{\mathrm{i}\pi}(t-\mathrm{i}\eta))=e^{-(a+b)\mathrm{i}\sqrt{t-\mathrm{i}\eta}}(t-\mathrm{i}\eta)^{-\frac{\mu+\nu}{2}}J_{\mu}(a\sqrt{t-\mathrm{i}\eta})J_{\nu}(b\sqrt{t-\mathrm{i}\eta}).$$
By using the asymptotic relation \eqref{a1}, we can validate the interchange of the integral and limit in equation \eqref{inv1}. Now, in view of Lemma \ref{LemStinv} we have that
$$\alpha'(t)=\frac{1}{\pi}t^{-\frac{\mu+\nu}{2}}J_{\mu}(a\sqrt{t})J_{\nu}(b\sqrt{t})\sin((a+b)\sqrt{t})$$
and thus $\alpha(t)$ is continuous and $\lim_{t \to 0^{+}}\alpha(t)$ exists. Let $\tilde{\alpha}(t)=\alpha(t)-\alpha(0)$. Note that $\tilde{\alpha}(t)$ is the normalized measure with $\tilde{\alpha}'(t)=\alpha'(t)$.
\end{proof}

\begin{proof}[\bf Proof of Theorem \ref{thprodeqIinfdiv}]
Due to Ismail \cite[Theorem 2]{ismail} we know that if $\nu>\frac{1}{2},$ then the function $x\mapsto 2^{\nu}\Gamma(\nu+1)x^{-\nu/2}I_{\nu}(\sqrt{x})e^{-\sqrt{x}}$ is the Laplace transform of an infinitely divisible distribution, which is not a generalized gamma convolution. This implies that $\zeta_{\mu,\nu,a,b}(x)$ it is in fact the product of two Laplace transforms of two infinitely divisible distributions, which are not generalized gamma convolutions. Clearly such a product will be also a Laplace transform of an infinitely divisible distribution. For the sake of completeness we briefly state the proof of the infinite divisibility. Moreover, by using the Pick characterization we prove that $\zeta_{\mu,\nu,a,b}(x)$ is not a Laplace transform of a generalized gamma convolution.

By using \cite[eq. 10.30.1]{nist}, we clearly have that $-\ln\zeta_{\mu,\nu,a,b}(x)\to0$ as $x\to0$. Observe that by using the same idea as in the proof of Theorem \ref{theolap2} we obtain that
\begin{align*}
-\frac{\ln \zeta_{\mu,\nu,a,b}(x)}{dx}(x)&=\frac{a+b}{2\sqrt{x}}+\frac{\mu+\nu}{2x}-\frac{a}{2\sqrt{x}}\frac{I_{\mu}'(a\sqrt{x})}{I_{\mu}(a\sqrt{x})}
-\frac{b}{2\sqrt{x}}\frac{I'_{\nu}(b\sqrt{x})}{I_{\nu}(b\sqrt{x})}\\
&=\left[\frac{a}{2\sqrt{x}}-\sum_{n\geq1}\frac{1}{x+j^2_{\mu,n}a^{-2}}\right]+\left[\frac{b}{2\sqrt{x}}-\sum_{n\geq1}\frac{1}{x+j^2_{\nu,n}b^{-2}}\right]\\
&=\left[\frac{1}{\pi}\int_{0}^{\infty}\frac{dt}{x+t^{2}a^{-2}}-\sum_{n\geq1}\frac{1}{x+j^2_{\mu,n}a^{-2}}\right]+
\left[\frac{1}{\pi}\int_{0}^{\infty}\frac{dt}{x+t^{2}b^{-2}}-\sum_{n\geq1}\frac{1}{x+j^2_{\nu,n}b^{-2}}\right]
\end{align*}
is a sum of two completely monotonic functions on $(0,\infty)$ if $a,b>0$ and $\mu,\nu>\frac{1}{2}.$

In view of Lemma \ref{lemma3}, in order to prove that the distribution with the Laplace transform $\zeta_{\mu,\nu,a,b}(x)$ does not belongs to the class of generalized gamma convolutions, it suffices to prove $\imag \left[{\psi'(s)}/{\psi(s)}\right]<0$, whenever $\imag s>0$, where $\psi(s)=\zeta_{\mu,\nu,a,b}(-s)$ is the moment generating function. Since
\begin{align*}
\frac{\psi'(s)}{\psi(s)}&=-\frac{(a+b)\mathrm{i}}{2\sqrt{s}}-\frac{\mu+\nu}{2s}+\frac{a}{2\sqrt{s}}\frac{J'_{\mu}(a\sqrt{s})}{J_{\mu}(a\sqrt{s})}+\frac{b}{2\sqrt{s}}\frac{J'_{\nu}(b\sqrt{s})}{J_{\nu}(b\sqrt{s})}\\
&=-\frac{(a+b)\mathrm{i}}{2\sqrt{s}}+\sum_{n\geq1}\frac{a^2}{s-j^2_{\mu,n}}+\sum_{n\geq1}\frac{b^2}{s-j^2_{\nu,n}}\\
&=\frac{\left({-\sin\frac{\theta}{2}-\mathrm{i}\cos\frac{\theta}{2}}\right)(a+b)}{2\sqrt{r}}+\sum_{n\geq1}\frac{a^2(x-j^2_{\mu,n}-\mathrm{i}y)}{(x-j^2_{\mu,n})^2+y^2}
+\sum_{n\geq1}\frac{b^2(a-j^2_{\nu,n}-\mathrm{i}y)}{(x-j^2_{\nu,n})^2+y^2},
\end{align*}
we arrive at $\imag \left[{\psi'(s)}/{\psi(s)}\right]<0$, whenever $s=re^{\mathrm{i}\theta}=x+\mathrm{i}y$ and $\imag s>0$, which completes the proof.
\end{proof}

\begin{proof}[\bf Proof of Theorem \ref{recprodKrepr}]
Define $F(z)=\frac{e^{-(a+b)\sqrt{z}}}{z^{\frac{\mu+\nu}{2}}K_{\mu}(a\sqrt{z})K_{\nu}(b\sqrt{z})}$ and let $\nu>0$ and $\mu>0$. By using the well-known asymptotic relation \cite[eq. 10.30.2]{nist}, we arrive at
\begin{equation*}
F(z)\sim\frac{a^{\mu}b^{\nu}e^{-(a+b)\sqrt{z}}}{2^{\mu+\nu-2}\Gamma(\mu)\Gamma(\nu)} \quad \text{as}\quad z \to 0.
\end{equation*}
Now, let $\mu>0$ and $\nu=0$. By using the relation \cite[eq. 10.30.3]{nist}, we have that
\begin{equation*}
F(z)\sim -\frac{a^{\mu}e^{-(a+b)\sqrt{z}}}{2^{\mu-1}\Gamma(\mu)\ln(b\sqrt{z})} \quad \text{as}\quad z \to 0.
\end{equation*}
Now, let $\mu>0$ and $\nu<0$. By using the relation $K_{\nu}(z)=K_{-\nu}(z)$ and \cite[eq. 10.30.2]{nist}, we obtain that
\begin{equation*}
F(z) \sim \frac{e^{-(a+b)\sqrt{z}}a^{\mu}b^{-\nu}z^{-\nu}}{2^{\mu-\nu-2}\Gamma(\mu)\Gamma(-\nu)} \quad \text{as} \quad z \to 0.
\end{equation*}
Consequently, we have $F(z)=o(|z|^{-1})$ as $|z| \to 0$ whenever $\mu>0$ and $\nu\in\mathbb{R},$ and by interchanging $\nu$ and $\mu$ we can conclude that in fact $F(z)=o(|z|^{-1})$ as $|z| \to 0$ whenever $\mu\in\mathbb{R}$, $\nu\in\mathbb{R}$ and $\mu+\nu>1$. On the other hand, by using the relation \cite[eq. 10.25.3]{nist}, for $\mu+\nu>1$ we have that
\begin{equation*}
F(z)\sim \frac{2\sqrt{ab}z^{\frac{1-(\mu+\nu)}{2}}}{\sqrt{\pi}} \quad \text{as} \quad |z| \to \infty.
\end{equation*}
Hence, $F(z)=o(1)$ as $|z| \to \infty$ in the region $|\arg z|<\frac{3\pi}{2}$. Hence, the condition {\bf b} in Lemma \ref{LemStrep2} is verified. Note that $K_{\mu}(z)$ has no zeros in $|\arg z|\leq\frac{\pi}{2}$ and it has finitely many zeros in $\mathbb{C}\setminus\Omega$, where $\Omega=\{z: |\arg z|<\frac{\pi}{2}\}$. Consequently, $K_{\mu}(a\sqrt{z})$ has no zeros in $|\arg z|\leq\pi$ and has finitely many zeros in $\pi<|\arg z|<2\pi$. Hence, we can find a $\theta \in (\frac{2}{3}, 1)$ such that the function $F(z)$ is analytic in the region $|\arg |<\frac{\pi}{\theta}$. Thus, the condition {\bf a} of Lemma \ref{LemStrep2} has been also verified. Furthermore, the contour integral in equation \eqref{repSt2} can be evaluated by using the residue theorem
\begin{equation*}
\varrho(t)=\frac{1}{2\pi \mathrm{i}}\int_{C}\frac{ze^{\frac{z}{2}}}{z^2+\pi^2}F(e^zt)dz=\frac{\mathrm{i}}{2}\left[F(te^{\mathrm{i}\pi})-F(te^{-\mathrm{i}\pi})\right],
\end{equation*}
where $C$ is a rectifiable closed curve going around $[-\mathrm{i}\pi, \mathrm{i}\pi]$ in the positive direction and lying in the strip $|\imag z|<\frac{\pi}{\theta}$. On the other hand in view of \eqref{recurIKl}, we have that
$$F(te^{\mathrm{i}\pi})= \frac{-4e^{-(a+b)\mathrm{i}\sqrt{t}}}{\pi^2t^{\frac{\mu+\nu}{2}}H^{(2)}_{\mu}(a\sqrt{t})H^{(2)}_{\nu}(b\sqrt{t})}$$
and
$$F(te^{-\mathrm{i}\pi})= \frac{-4e^{(a+b)\mathrm{i}\sqrt{t}}}{\pi^2t^{\frac{\mu+\nu}{2}}H^{(1)}_{\mu}(a\sqrt{t})H^{(1)}_{\nu}(b\sqrt{t})}.$$
Consequently, we arrive at
$$\varrho(t)=\frac{4}{t^{\frac{\mu+\nu}{2}}\pi^2}\cdot \frac{{}_1T_{\mu,\nu,a,b}(t)\cos((a+b)\sqrt{t})-{}_2T_{\mu,\nu,a,b}(t)\sin((a+b)\sqrt{t})}{\left[J^2_{\mu}(a\sqrt{t})+Y^2_{\mu}(a\sqrt{t})\right]
\left[J^2_{\nu}(b\sqrt{t})+Y^2_{\nu}(b\sqrt{t})\right]}.$$
\end{proof}

\begin{proof}[\bf Proof of Theorem \ref{recprodKinfdiv}]
Due to Ismail \cite[Theorem 1]{ismail} we know that if $\nu>\frac{1}{2},$ then the function $x\mapsto 2^{\nu-1}\Gamma(\nu)x^{-\nu/2}e^{-\sqrt{x}}/K_{\nu}(\sqrt{x})$ is the Laplace transform of an infinitely divisible distribution, which is in fact a generalized gamma convolution. This implies that $\kappa_{\mu,\nu,a,b}(x)$ it is in fact the product of two Laplace transforms of two infinitely divisible distributions, which are generalized gamma convolutions. Clearly such a product will be also a Laplace transform of an infinitely divisible distribution, which is also a generalized gamma convolution (see \cite[Theorem 1]{bonde} or \cite[Proposition 7]{behme}).

For the sake of completeness we briefly state the proof of the infinite divisibility. Observe that
\begin{align*}
-&\frac{d\ln\kappa_{\mu,\nu,a,b}(x)}{dx}\\&=\frac{a+b}{2\sqrt{x}}+\frac{\mu+\nu}{2x}+\frac{a}{2\sqrt{x}}\frac{K'_{\mu}(a\sqrt{x})}{K_{\mu}(a\sqrt{x})}+\frac{b}{2\sqrt{x}}\frac{K'_{\nu}(b\sqrt{x})}{K_{\nu}(b\sqrt{x})}\\
&=\frac{a}{2\sqrt{x}}-\frac{a}{2\sqrt{x}}\frac{K_{\mu-1}(a\sqrt{x})}{K_{\mu}(a\sqrt{x})}+\frac{b}{2\sqrt{x}}-\frac{b}{2\sqrt{x}}\frac{K_{\nu-1}(b\sqrt{x})}{K_{\nu}(b\sqrt{x})}\\
&=\frac{1}{\pi}\int_{0}^{\infty}\frac{1}{x+t^2a^{-2}}\left(1-\frac{2(\pi t)^{-1}}{J^2_{\mu}(t)+Y^2_{\mu}(t)}\right)dt+\frac{1}{\pi}
\int_{0}^{\infty}\frac{1}{x+t^2b^{-2}}\left(1-\frac{2(\pi t)^{-1}}{J^2_{\nu}(t)+Y^2_{\nu}(t)}\right)dt.
\end{align*}
which in view of the inequality $J^2_{\mu}(t)+Y^2_{\mu}(t)>2(\pi t)^{-1}$, whenever $\mu>\frac{1}{2}$ (see \cite[p. 447]{watson}), implies that indeed $x\mapsto -{d\ln\kappa_{\mu,\nu,a,b}(x)}/{dx}$ is completely monotonic on $(0,\infty)$ for all $\mu,\nu>\frac{1}{2}$ and $a,b>0.$ Consequently, in view of Lemma \ref{lemma1} we obtain that $\kappa_{\mu,\nu,a,b}(x)$ is the Laplace transform of an infinitely divisible distribution.
\end{proof}

\begin{proof}[\bf Proof of Theorem \ref{theoquotIK}]
Let $$F(z)=z^{-\frac{\mu+\nu}{2}}e^{-(a+b)\sqrt{z}}\frac{I_{\mu}(a\sqrt{z})}{K_{\nu}(b\sqrt{z})}.$$ In the case when $\nu>0$, by using the relations \cite[eq. 10.30.1]{nist} and \cite[eq. 10.30.2]{nist}, we have that
$$F(z) \sim \frac{a^{\mu}b^{\nu}e^{-(a+b)\sqrt{z}}}{2^{\mu+\nu-1}\Gamma(\mu+1)\Gamma(\nu)} \quad \text{as} \quad z \to 0.$$
Now, in the case when $\nu=0$, by using the relations \cite[eq. 10.30.1]{nist} and \cite[eq. 10.30.3]{nist}, we have that
$$F(z) \sim -\frac{e^{-(a+b)\sqrt{z}}a^{\mu}}{2^{\mu}\Gamma(\mu+1)\ln(b\sqrt{z})} \quad \text{as} \quad z \to 0.$$
Moreover, in the case when $\nu<0$, by using the relations \cite[eq. 10.30.1]{nist}, \cite[eq. 10.27.3]{nist} and \cite[eq. 10.30.2]{nist}, we have that
$$F(z)\sim \frac{a^{\mu}2^{\nu-\mu+1}e^{-(a+b)\sqrt{z}}z^{-\nu}}{b^{\nu}\Gamma(-\nu)\Gamma(\mu+1)} \quad \text{as} \quad z \to 0.$$
On the other hand, by using the relations \cite[eq. 10.30.4]{nist}  and \cite[eq. 10.25.3]{nist}, we obtain that
$$F(z) \sim \sqrt{\frac{b}{a}}\frac{z^{-\frac{\mu+\nu}{2}}}{\pi} \quad \text{as} \quad |z| \to \infty$$
and whenever $a,b>0$ and $\mu+\nu>0.$ Thus, condition {\bf b} of Lemma \ref{LemStrep1} holds true. Note that $K_{\nu}(z)$ has no zeros in $|\arg z|\leq\frac{\pi}{2}$. Consequently, $K_{\nu}(b\sqrt{z})$ has no zeros in $|\arg z|\leq\pi$. Thus, condition {\bf a} of Lemma \ref{LemStrep1} is also verified. By using the corresponding equations in \eqref{recurIKl}, we obtain that
$$F(-t-\mathrm{i}\eta)=F(e^{-\mathrm{i}\pi}(t+\mathrm{i}\eta))=-\frac{2\mathrm{i}(t+\mathrm{i}\eta)^{-\frac{\mu+\nu}{2}}e^{(a+b)\mathrm{i}\sqrt{t+\mathrm{i}\eta}}}{\pi}
\frac{J_{\mu}(a\sqrt{t+\mathrm{i}\eta})}{H_{\nu}^{(1)}(b\sqrt{t+\mathrm{i}\eta})}$$
and
$$F(-t+\mathrm{i}\eta)=F(e^{\mathrm{i}\pi}(t-\mathrm{i}\eta))=\frac{2\mathrm{i}(t-\mathrm{i}\eta)^{-\frac{\mu+\nu}{2}}e^{-(a+b)\mathrm{i}\sqrt{t-\mathrm{i}\eta}}}{\pi}
\frac{J_{\mu}(a\sqrt{t-\mathrm{i}\eta})}{H_{\nu}^{(2)}(b\sqrt{t-\mathrm{i}\eta})}.$$
Here the validation for the interchange of limit and integral in equation \eqref{inv1} is established in three different cases.

\textbf{Case (i) $\nu>0$ and $\mu>-1.$} By using the asymptotic relations \eqref{a1} and \eqref{a2}, we conclude that the function $F(-t-\mathrm{i}\eta)-F(-t+\mathrm{i}\eta)$ is continuous in every rectangle $[t_1, t_2]\times[0, \eta]$.

\textbf{Case (ii) $\nu=0$ and $\mu>-1.$} By using the relations \eqref{a1} and \eqref{a3}, we obtain that
$$F(-t-\mathrm{i}\eta) \sim -\dfrac{a^{\mu}}{2^{\mu}\Gamma(\mu+1)}\dfrac{1}{\ln(b\sqrt{t+\mathrm{i}\eta})},\quad t \to 0, \eta \to 0$$ and
$$F(-t+\mathrm{i}\eta) \sim -\dfrac{a^{\mu}}{2^{\mu}\Gamma(\mu+1)}\dfrac{1}{\ln(b\sqrt{t-\mathrm{i}\eta})},\quad t \to 0, \eta \to 0.$$
Consequently, we note that both functions $F(-t-\mathrm{i}\eta)$ and $F(-t+\mathrm{i}\eta)$ are bounded in every rectangle $[t_1, t_2]\times[0, \eta]$.

\textbf{Case (iii) $\nu<0$ and $\mu>-1.$} By using the asymptotic relations \eqref{a1}, \eqref{a2} and the relation \eqref{a4}, we obtain that
$$F(-t-\mathrm{i}\eta)\sim \dfrac{e^{\nu\pi \mathrm{i}}a^{\mu}}{2^{\mu-\nu-1}b^{\nu}\Gamma(\mu+1)\Gamma(-\nu)}(t+\mathrm{i}\eta)^{-\nu},\quad t \to 0, \eta \to 0$$ and
$$F(-t+\mathrm{i}\eta)\sim \dfrac{e^{-\nu\pi \mathrm{i}}a^{\mu}}{2^{\mu-\nu-1}b^{\nu}\Gamma(\mu+1)\Gamma(-\nu)}(t-\mathrm{i}\eta)^{-\nu},\quad t \to 0, \eta \to 0.$$
Consequently, we observe that the functions $F(-t-\mathrm{i}\eta)$ and $F(-t+\mathrm{i}\eta)$ are bounded in every rectangle $[t_1, t_2]\times[0, \eta]$.

Finally, in view of Lemma \ref{LemStinv}, we have that
$$\frac{d\alpha(t)}{dt}=-\frac{2t^{-\frac{\mu+\nu}{2}}}{\pi^2}J_{\mu}(a\sqrt{t})\cdot\frac{\cos((a+b)\sqrt{t})J_{\nu}(b\sqrt{t})+\sin((a+b)\sqrt{t})Y_{\nu}(b\sqrt{t})}{J^2_{\nu}(b\sqrt{t})+Y^2_{\nu}(b\sqrt{t})}.$$
\end{proof}

\begin{proof}[\bf Proof of Corollary \ref{theoquotIKcoro}]
The proof of this corollary follows naturally by tending with $a$ to zero in Theorem \ref{theoquotIK}, however it is also possible to have a direct proof in view of Lemma \ref{LemStrep2}. In other words, although Corollary \ref{theoquotIKcoro} follows immediately from Theorem \ref{theoquotIK}, for the sake of completeness we present its proof via the Stieltjes transform representation theorem. For this let $$F(z)=(\sqrt{z})^{-\nu}e^{-b\sqrt{z}}\frac{1}{K_{\nu}(b\sqrt{z})}.$$
By using the relation \cite[eq. 10.30.2]{nist}, we obtain that
$$F(z)\sim e^{-b\sqrt{z}}\frac{2^{1-\nu}b^{\nu}}{\Gamma(\nu)} \quad \text{as} \quad z \to 0,$$
and thus, $F(z)\sim o(|z|^{-1})$ as $z \to 0$. On the other hand, by using the relation \cite[eq. 10.25.3]{nist}, we arrive at
$$F(z) \sim (\sqrt{z})^{-\nu+\frac{1}{2}}\sqrt{\frac{2b}{\pi}} \quad \text{as} \quad |z| \to \infty.$$
Hence the condition {\bf b} in Lemma \ref{LemStrep2} holds true. Note that $K_{\nu}(z)$ does not have any zero in the region $|\arg z|\leq\frac{\pi}{2}$ and it has finitely many zeros in the region $\frac{\pi}{2}<\arg z<\pi$ and $-\pi<\arg z<-\frac{\pi}{2}$. Consequently, the function $F(z)$ is analytic in the region $|\arg z|\leq\pi$. Moreover we can find an $\epsilon\in (\frac{2}{3}, 1)$ such that $F(z)$ is analytic in the region $|\arg z|<\frac{\pi}{\epsilon}$. Hence, the condition {\bf a} in Lemma \ref{LemStrep2} also holds true. Now, in view of Lemma \ref{LemStinv} we obtain that
$$\frac{d\alpha(t)}{dt}=\frac{1}{2\pi \mathrm{i}}\left[F(e^{-\mathrm{i}\pi}t)-F(e^{\mathrm{i}\pi}t)\right],$$
where
$$F(e^{-\mathrm{i}\pi}t)=e^{\mathrm{i}\frac{\pi}{2}\nu}(\sqrt{t})^{-\nu}e^{\mathrm{i}b\sqrt{t}}\frac{1}{K_{\nu}(e^{-\mathrm{i}\frac{\pi}{2}}b\sqrt{t})}
=\frac{2(\sqrt{t})^{-\nu}}{\mathrm{i}\pi}e^{\mathrm{i}b\sqrt{t}}\frac{J_{\nu}(b\sqrt{t})-\mathrm{i}Y_{\nu}(b\sqrt{t})}{J^2_{\nu}(b\sqrt{t})+Y^2_{\nu}(b\sqrt{t})}$$
and
$$F(e^{\mathrm{i}\pi}t)=e^{-\mathrm{i}\frac{\pi}{2}\nu}(\sqrt{t})^{-\nu}e^{-\mathrm{i}b\sqrt{t}}\frac{1}{K_{\nu}(e^{\mathrm{i}\frac{\pi}{2}}b\sqrt{t})}=
-\frac{2(\sqrt{t})^{-\nu}}{\mathrm{i}\pi}e^{-\mathrm{i}b\sqrt{t}}\frac{J_{\nu}(b\sqrt{t})+\mathrm{i}Y_{\nu}(b\sqrt{t})}{J^2_{\nu}(b\sqrt{t}+Y^2_{\nu}(b\sqrt{t}))}.$$
Consequently, we arrive at
$$\frac{d\alpha(t)}{dt}=-\frac{2}{\pi^2}(\sqrt{t})^{-\nu}\frac{J_{\nu}(b\sqrt{t})\cos(b\sqrt{t})+Y_{\nu}(b\sqrt{t})\sin(b\sqrt{t})}{J^2_{\nu}(b\sqrt{t})+Y^2_{\nu}(b\sqrt{t})},$$
which completes the proof.
\end{proof}

\begin{proof}[\bf Proof of Theorem \ref{theoquotIKinfdiv}]
Recall that due to Ismail \cite[Theorem 2]{ismail} we know that if $\mu>\frac{1}{2},$ then the function $x\mapsto 2^{\mu}\Gamma(\mu+1)x^{-\mu/2}I_{\mu}(\sqrt{x})e^{-\sqrt{x}}$ is the Laplace transform of an infinitely divisible distribution, which is not a generalized gamma convolution. Moreover, due to Ismail \cite[Theorem 1]{ismail} we also know that if $\nu>\frac{1}{2},$ then the function $x\mapsto 2^{\nu-1}\Gamma(\nu)x^{-\nu/2}e^{-\sqrt{x}}/K_{\nu}(\sqrt{x})$ is the Laplace transform of an infinitely divisible distribution, which is in fact a generalized gamma convolution. This implies that $\varepsilon_{\mu,\nu,a,b}(x)$ it is in fact the product of two Laplace transforms of two infinitely divisible distributions, and clearly such a product will be also a Laplace transform of an infinitely divisible distribution, according to \cite[Proposition 2.1]{steutel}. More precisely, we have that
\begin{align*}
-&\frac{d\ln\varepsilon_{\mu,\nu,a,b}(x)}{dx}(x)\\&=\frac{a+b}{2\sqrt{x}}+\frac{\mu+\nu}{2x}-\frac{a}{2\sqrt{x}}\frac{I'(a\sqrt{x})}{I_{\mu}(a\sqrt{x})}+\frac{b}{2\sqrt{x}}\frac{K'_{\nu}(b\sqrt{x})}{K_{\nu}(b\sqrt{x})}\\
&=\frac{a+b}{2\sqrt{x}}-\frac{a}{2\sqrt{x}}\frac{I_{\mu+1}(a\sqrt{x})}{I_{\mu}(a\sqrt{x})}-\frac{b}{2\sqrt{x}}\frac{K_{\nu-1}(b\sqrt{x})}{K_{\nu}(b\sqrt{x})}\\
&=\left[\frac{1}{\pi}\int_{0}^{\infty}\frac{1}{x+b^{-2}t^2}\left(1-\frac{2(\pi t)^{-1}}{J^2_{\nu}(t)+Y^2_{\nu}(t)}\right)dt\right]+\left[\frac{1}{\pi}\int_{0}^{\infty}\frac{dt}{x+a^{-2}t^2}-\sum_{n\geq1}\dfrac{1}{x+a^{-2}j^2_{\mu,n}}\right]
\end{align*}
is completely monotonic on $(0,\infty)$ for all $\mu,\nu>\frac{1}{2}$ and $a,b>0$ as a sum of two completely monotonic functions. Indeed the expression in the first big brackets is completely monotonic because of the well-known inequality $J^2_{\nu}(t)+Y^2_{\nu}(t)>\frac{2}{\pi t}$ where $\nu>\frac{1}{2}$ (see \cite[p. 447]{watson}), while the expression in the second big brackets satisfies
$$\Theta(x)=\frac{1}{\pi}\int_{0}^{\infty}\frac{dt}{x+a^{-2}t^2}-\sum_{n\geq1}\dfrac{1}{x+a^{-2}j^2_{\mu,n}}\geq\frac{1}{\pi}\sum_{n\geq1}
\frac{j_{\mu,n}-j_{\mu,n-1}}{x+a^{-2}j^2_{\mu,n}}-\sum_{n\geq1}\frac{1}{x+a^{-2}j^2_{\mu,n}}\geq 0$$
and
$$\frac{(-1)^{m}}{m!}\frac{d^m\Theta(x)}{dx^m}\geq\frac{1}{\pi}\sum_{n\geq1}\frac{j_{\mu,n}-j_{\mu,n-1}}{(x+a^{-2}j^2_{\mu,n})^{m+1}}-\sum_{n\geq1}\frac{1}{(x+a^{-2}j^2_{\mu,n})^{m+1}}\geq 0$$
since $j_{\mu,n}-j_{\mu,n-1}>\pi$ whenever $\mu>\frac{1}{2}$ and $n\in\mathbb{N}.$ Thus, in view of Lemma \ref{lemma1} the proof of the first part is complete.

Now, for the second part of this theorem we just need to observe that by using a similar approach as before the corresponding expression
\begin{align*}
-\frac{d}{dx}\left[\ln\left(\frac{e^{-(a+b)\sqrt{x}}}{\varepsilon_{\mu,\nu,a,b}(x)}\right)\right]
&=-\frac{\nu+\mu}{2x}-\frac{b}{2\sqrt{x}}\frac{K'_{\nu}(b\sqrt{x})}{K_{\nu}(b\sqrt{x})}+\frac{a}{2\sqrt{x}}\frac{I'_{\mu}(a\sqrt{x})}{I_{\mu}(a\sqrt{x})}\\
&=\frac{b}{2\sqrt{x}}\frac{K_{\nu-1}(b\sqrt{x})}{K_{\nu}(b\sqrt{x})}+\frac{a}{2\sqrt{x}}\frac{I_{\mu+1}(a\sqrt{x})}{I_{\mu}(a\sqrt{x})}\\
&=\dfrac{2}{\pi^2}\int_{0}^{\infty}\frac{1}{x+b^{-2}t^2}[J^2_{\nu}(t)+Y^2_{\nu}(t)]^{-1}dt + \sum_{n\geq1}\frac{1}{x+a^{-2}j^2_{\mu,n}}
\end{align*}
is completely monotonic on $(0,\infty)$ as a sum of two completely monotonic functions whenever $a,b,\nu>0$ and $\mu>-1.$
\end{proof}

\begin{proof}[\bf Proof of \ref{Thquogamma}]
Let $X$ and $Y$ be independent gamma variables with parameters $(\alpha, \beta)$ and $(\alpha_0, \beta_0)$. Recall that the probability density function of the quotient of these two random variables $Z={X}/{Y}$ is given by (see for example \cite[p. 889]{kelker})
$$f(z)=\frac{\Gamma(\alpha+\alpha_0)}{\Gamma(\alpha)\Gamma(\alpha_0)}\left(\frac{\beta_0}{\beta}\right)^{\alpha}x^{\alpha-1}\left(1+\frac{\beta_0}{\beta}x\right)^{-(\alpha+\alpha_0)},$$
where $x>0.$ Now let $a=\alpha$ and $c=1-\alpha_0.$ The Laplace transform $L(s)$ of the quotient of two gamma random variable is given by \cite[p.889]{kelker}
$$L(s)=\frac{\Gamma(a-c+1)}{\Gamma(1-c)}\psi(a,c,s)$$
and thus the corresponding moment generating function $\phi(s)=L(-s)$ is as follows
\begin{equation}\label{lg}	
\phi(s)=\frac{\Gamma(a-c+1)}{\Gamma(1-c)}\psi(a,c,-s).
\end{equation}	
Since $a>0$, by using \cite[eq. 13.9.13]{nist}, we observe that $\psi(a,c,-s)$ has no zeros in $\mathbb{C}\setminus [0,\infty).$ Thus the first condition of the Pick function characterization theorem, that is, Lemma \ref{lemma3} is verified and by taking the logarithmic derivative of both sides of the equation \eqref{lg}, we obtain that
$$\frac{\phi'(s)}{\phi(s)}=-\frac{\psi'(a,c,-s)}{\psi(a,c,-s)}=\frac{a\psi(a+1,c+1,-s)}{\psi(a,c,-s)}.$$ In the case when $s=x+\mathrm{i}y$, by using the integral representation \eqref{intfor}, which is valid for $|\arg z|<\pi,$ $a>0$ and $c<1,$ we arrive at
\begin{align*}
\frac{\phi'(s)}{\phi(s)}&=\int_{0}^{\infty}\frac{at^{-c}e^{-t}\left|\psi(a,c,te^{\mathrm{i}\pi})\right|^{-2}dt}{(-s+t)\Gamma(a+1)\Gamma(a-c+1)} \\
&=\int_{0}^{\infty}\frac{at^{-c}e^{-t}\left|\psi(a,c,te^{\mathrm{i}\pi})\right|^{-2}dt}{(-x-\mathrm{i}y+t)\Gamma(a+1)\Gamma(a-c+1)}\\
&=\int_{0}^{\infty}\frac{(t-x+\mathrm{i}y)at^{-c}e^{-t}\left| \psi(a,c,te^{\mathrm{i}\pi})\right|^{-2}dt}{((t-x)^2+y^2)\Gamma(a+1)\Gamma(a-c+1)}.
\end{align*}
From the above equation, we clearly conclude that $\imag {\phi'(s)}/{\phi(s)}\geq 0$ whenever $\imag s>0$.
\end{proof}

\begin{proof}[\bf Proof of Theorem \ref{tricrepr}]
By using the relations \cite[eq. (12)-(15), p. 258]{erdelyi1} and $\psi'(a,c,z)=-a\psi(a+1,c+1,z)$, we obtain that
$$\psi(a,c-1,z)=\frac{1-c}{a-c+1}\psi(a,c,z)+\frac{az}{a-c+1}\psi(a+1,c+1,z),$$
$$\psi(a+1,c,z)=\frac{1}{a-c+1}\psi(a,c,z)-\frac{z}{a-c+1}\psi(a+1,c+1,z),$$
$$\psi(a,c+1,z)=\psi(a,c,z)+a\psi(a+1,c+1,z)$$
and
$$\psi(a-1,c,z)=z\psi(a,c,z)-(c-a)\psi(a,c,z)+az\phi(a+1,c+1,z).$$
Dividing both sides of the above equations by $\psi(a,c,z)$ and with the help of the integral representation \eqref{intfor} we arrive to the required results.

Thus, indeed the representations in Theorem \ref{tricrepr} follow naturally from \eqref{intfor}, however for \eqref{triceqnew} we give a more detailed proof via the Stieltjes representation and inversion theorems.
For this let $$F(z)=\frac{\psi(a,c+1,z)}{\psi(a,c,z)}-1.$$ By using \cite[eq. 13.7.3]{nist}, we have that
$$F(z) \sim \left.{\sum_{s\geq1}\frac{(a)_{s}((a-c)_s-(a-c+1)_s)(-z)^{-s}}{s!}}\right/{\sum_{s\geq0}\frac{(a)_s(a-c+1)_s(-z)^{-s}}{s!}} \quad \text{as}\quad z \to \infty,$$
and $|\arg z|<\frac{3}{2}\pi$. Consequently, we obtain that $f(z) \to 0$ as $|z| \to \infty$. Next, for $c<-1$, by using \cite[eq. 13.2.22]{nist}, we have that
$$F(z) \sim \frac{\Gamma(-c)/\Gamma(a-c)}{\Gamma(1-c)/\Gamma(a-c+1)}-1 \quad \text{as} \quad z \to 0,$$
for $c=-1$, by using \cite[eq. 13.2.21]{nist} and \cite[eq. 13.2.22]{nist}, we obtain that
$$F(z) \to \frac{1/\Gamma(a+1)}{\Gamma(1-c)/\Gamma(a-c+1)}-1 \quad \text{as} \quad z \to 0,$$
for $-1<c<0$, by using \cite[eq. 13.2.20]{nist} and \cite[eq. 13.2.22]{nist}, we arrive at
$$F(z) \sim \frac{\Gamma(-c)/\Gamma(a-c)}{\Gamma(1-c)/\Gamma(a-c+1)}-1 \quad \text{as} \quad z \to 0,$$
for $c=0,$ by using \cite[eq. 13.2.19]{nist} and \cite[eq. 13.2.21]{nist}, we have that
$$F(z) \sim \frac{-1/\Gamma(a)(\ln(z)+d)}{1/\Gamma(a+1)}-1 \quad \text{as} \quad z \to 0,$$
where $d$ is the constant in (see \cite[eq. 13.2.19]{nist}), and finally for $0<c<1$, by using \cite[eq. 13.2.18]{nist} and \cite[eq. 13.2.20]{nist}, we have that
$$F(z) \sim \frac{\Gamma(c)/\Gamma(a)z^{-c}+\Gamma(-c)/\Gamma(a-c)}{\Gamma(1-c)/\Gamma(a-c+1)} \quad \text{as} \quad z \to 0.$$
Consequently, we obtain that $F(z)=o(|z|^{-1})$ as $z \to 0,$ for $a>0$ and $c<1$. Hence the conditions in Lemma \ref{LemStrep2} have been verified, because $\psi(a,c,z)$ has no zeros in $|\arg z|<\frac{\pi}{\alpha}$, where $\alpha \in (\frac{2}{3},1).$ By using the residue calculus we obtain that ${d\alpha(t)}/{dt}$ becomes
\begin{align*}
\frac{1}{2\pi \mathrm{i}}&\left[F(te^{-\mathrm{i}\pi})-F(te^{\mathrm{i}\pi})\right]=\frac{1}{2\pi \mathrm{i}}\lim_{\eta \to 0^{+}}\left[F(-t-\mathrm{i}\eta)-F(-t+\mathrm{i}\eta)\right]\\
	&=\frac{1}{2\pi \mathrm{i}}\lim_{\eta \to 0^{+}}\left[\frac{\psi(a,c,e^{\mathrm{i}\pi}(t-\mathrm{i}\eta))\psi(a,c+1,e^{-\mathrm{i}\pi}(t+\mathrm{i}\eta))-\psi(a,c,e^{-\mathrm{i}\pi}(t+\mathrm{i}\eta))\psi(a,c+1,e^{\mathrm{i}\pi}(t-\mathrm{i}\eta))}
{|\psi(a,c,e^{\mathrm{i}\pi}(t-\mathrm{i}\eta))|^2}\right].
\end{align*}
Now, by using \cite[eq.(14), p. 263]{erdelyi1}, for $-a<\min\{0,1-c\}$ we have that
\begin{equation}\label{eqT2}
	\lim_{\eta \to 0^{+}}\psi(a,c,e^{\pm \mathrm{i}\pi}(t\mp \mathrm{i}\eta))=k_1y_1(-t)-e^{\mp \mathrm{i}\pi c}k_2y_2(-t),
\end{equation}
where $$y_1(x)=\Phi(a,c,x),\quad y_2(x)=x^{1-c}\Phi(a-c+1,2-c,x)$$
with $\Phi(a,c,x)$ being the Kummer confluent hypergeometric function and
$$k_1=\frac{\Gamma(1-c)}{\Gamma(a-c+1)},\quad k_2=(-1)^{1-c}\dfrac{\Gamma(c-1)}{\Gamma(a)}.$$
On the other hand, we have that (see \cite[eq. (9), p. 253]{erdelyi1})
\begin{equation}\label{eqT3}
W[y_1(x), y_2(x)]=y_1(x)y'_2(x)-y_2(x)y'_1(x)=(1-c)x^{-c}e^{x}
\end{equation}
and by using the recurrence relation $\psi(a,c+1,x)=-\psi'(a,c,x)+\psi(a,c,x)$ (see \cite[eq.(14), p. 258]{erdelyi1}) and equation \eqref{eqT2}, we obtain that
\begin{equation}\label{eqT4}
	\lim_{\eta \to 0^{+}}\psi(a,c+1,e^{\mp i \pi}(t\pm i\eta))=k_1y'_1(-t)-e^{\pm i\pi c}k_2y'_2(-t)+k_1y_1(-t)-e^{\pm i\pi c}k_2y_2(-t).
\end{equation}
In view of \eqref{eqT4}, \eqref{eqT3} and \eqref{eqT2} it follows that
$$\frac{d\alpha(t)}{dt}=-\dfrac{\sin(\pi c)k_1k_2(y_1(-t)y_2'(-t)-y_2(-t)y_1'(-t))}{\pi |\psi(a,c,e^{i \pi}t)|^{2}}$$
and by using the formula $\sin(\pi c)=\pi \left[\Gamma(c-1)\Gamma(2-c)\right]^{-1}$ (see \cite[p. 890]{kelker}) we arrive at
$$\frac{d\alpha(t)}{dt}=\frac{t^{-c}e^{-t}\left|\psi(a,c,e^{\mathrm{i}\pi}t)\right|^{-2}}{\Gamma(a)\Gamma(a-c+1)}.$$
This completes the proof of \eqref{triceqnew}.
\end{proof}

\begin{proof}[\bf Proof of Theorem \ref{noncentralchihcm}]
First observe that
$$\chi_{\mu,\lambda}(uv)\chi_{\mu,\lambda}\left(\frac{u}{v}\right)=\frac{e^{-\lambda}}{4}\left(\frac{u}{\lambda}\right)^{\frac{\mu}{2}-1}e^{-\frac{u}{2}w}
I_{\frac{\mu}{2}-1}(\sqrt{\lambda uv})I_{\frac{\mu}{2}-1}\left(\sqrt{\frac{\lambda u}{v}}\right).$$
Now, by using the following well-known integral representation for the product of modified Bessel functions of the first kind (see \cite[p. 90]{MOS66})
\begin{equation}\label{intIprod}
I_{\mu}(a)I_{\mu}(b)=\frac{\left(\frac{1}{2}ab\right)^{\mu}}{\sqrt{\pi}\Gamma\left(\frac{1}{2}+\mu\right)}
\int_{0}^{\pi}(a^2+b^2-2ab\cos t)^{-\frac{1}{2}\mu}I_{\mu}((a^2+b^2-2ab\cos t)^{\frac{1}{2}})\sin^{2\mu}t dt,
\end{equation}
where $\mu>-\frac{1}{2},$ in view of the notation $T=\sqrt{\lambda u(w-2\cos t)},$ we obtain that
\begin{equation}\label{chiprodrepr}
\chi_{\mu,\lambda}(uv)\chi_{\mu,\lambda}\left(\frac{u}{v}\right)=c_{\mu,\lambda}(u)
\cdot e^{-\frac{u}{2}w}\int_{0}^{\pi}\frac{I_{\frac{\mu}{2}-1}(T)}{T^{\frac{\mu}{2}-1}}\sin^{\mu-2}t dt
\end{equation}
with
$$c_{\mu,\lambda}(u)=\frac{e^{-\lambda}u^{(\mu-2)}}{2^{\frac{\mu}{2}+1}\sqrt{\pi}\Gamma(\frac{\mu-1}{2})}.$$

Now, observe that ${2dT}/{dw}={\lambda u}/{T}$ and by using the recurrence relation $\left[z^{-\mu}I_{\mu}(z)\right]'=z^{-\mu}I_{\mu+1}(z)$ (see for example \cite[p. 79]{watson}), after differentiating both sides of the equation \eqref{chiprodrepr} with respect to $w=v+{1}/{v}$, we arrive at
\begin{align*}
\frac{d}{dw}\left[\chi_{\mu, \lambda}(uv)\chi_{\mu, \lambda}\left(\frac{u}{v}\right)\right]&=c_{\mu,\lambda}(u)\left[-\frac{u}{2}e^{-\frac{uw}{2}}\int_{0}^{\pi}\dfrac{I_{\frac{\mu}{2}-1}(T)}{T^{\frac{\mu}{2}-1}}\sin^{\mu-2}t dt+e^{-\frac{uw}{2}}\int_{0}^{\pi}\frac{\lambda u}{2}\frac{I_{\frac{\mu}{2}}(T)}{T^{\frac{\mu}{2}}}\sin^{\mu-2}t dt\right]\\
	&=c_{\mu,\lambda}(u)\left[\frac{u}{2}e^{-\frac{uw}{2}}\int_{0}^{\pi}\frac{\sin^{\mu-2}t}{T^{\frac{\mu}{2}-1}}\left(-I_{\frac{\mu}{2}-1}(T)+\lambda\dfrac{I_{\frac{\mu}{2}}(T)}{T}\right)dt\right].
\end{align*}
By using the recurrence relation $\frac{2\mu}{z}I_{\mu}(z)=I_{\mu-1}(z)-I_{\mu+1}(z)$ (see \cite[p. 79]{watson}), we obtain that
$$\frac{\lambda}{T}I_{\frac{\mu}{2}}(T)=\frac{\lambda}{\mu}I_{\frac{\mu}{2}-1}(T)-\frac{\lambda}{\mu}I_{\frac{\mu}{2}+1}(T)$$
and consequently we arrive at
$$\frac{d}{dw}\left[\chi_{\mu,\lambda}(uv)\chi_{\mu,\lambda}\left(\frac{u}{v}\right)\right]
=\frac{\frac{u}{2}\cdot c_{\mu,\lambda}(u)}{e^{\frac{u}{2}w}}\int_{0}^{\pi}\dfrac{\sin^{\mu-2}t}{T^{\frac{\mu}{2}-1}}
\left[\left(\frac{\lambda}{\mu}-1\right)I_{\frac{\mu}{2}-1}(T)-\frac{\lambda}{\mu}I_{\frac{\mu}{2}+1}(T)\right]dt<0$$
whenever $0<\lambda\leq\mu$ and $u>0.$ By using a similar approach we obtain that
\begin{align*}
\frac{d^2}{dw^2}\left[\chi_{\mu, \lambda}(uv)\chi_{\mu, \lambda}\left(\frac{u}{v}\right)\right]
&=c_{\mu,\lambda}(u)\left[-\frac{u^2}{4}e^{-\frac{u}{2}w}\int_{0}^{\pi}\sin^{\mu-2}t\left(-\frac{I_{\frac{\mu}{2}-1}(T)}{T^{\frac{\mu}{2}-1}}+\lambda\dfrac{I_{\frac{\mu}{2}}(T)}{T^{\frac{\mu}{2}}}\right)dt\right.\\
&\left.+\frac{u}{2}e^{-\frac{u}{2}w}\int_{0}^{\pi}\sin^{\mu-2}t\left(-\frac{\lambda u}{2}\frac{I_{\frac{\mu}{2}}(T)}{T^{\frac{\mu}{2}}}+\frac{\lambda^2u}{2}\dfrac{I_{\frac{\mu}{2}+1}(T)}{T^{\frac{\mu}{2}+1}}\right)dt\right]\\
&=c_{\mu,\lambda}(u)\cdot \frac{u^2}{4}e^{-\frac{u}{2}w}\int_{0}^{\pi}\frac{\sin^{\mu-2}t}{T^{\frac{\mu}{2}-1}}\left(I_{\frac{\mu}{2}-1}(T)-\frac{2\lambda}{T}I_{\frac{\mu}{2}}(T)+\frac{\lambda^2}{T^2}I_{\frac{\mu}{2}+1}(T)\right) dt,
\end{align*}
which can be rewritten as
\begin{equation*}
\frac{u^2\cdot c_{\mu,\lambda}(u)}{4e^{\frac{u}{2}w}}\int_{0}^{\pi}\dfrac{\sin^{\mu-2}t}{T^{\frac{\mu}{2}-1}}
\left[\left(1-\frac{2\lambda}{\mu}\right)I_{\frac{\mu}{2}-1}(T)+\frac{2\lambda}{\mu}I_{\frac{\mu}{2}+1}(T)+\frac{\lambda^2}{\mu^2}I_{\frac{\mu}{2}+1}(T)\right] dt
\end{equation*}
and this is strictly positive whenever $0<\lambda\leq{\mu}/{2}$ and $u>0.$
\end{proof}

\begin{proof}[\bf Proof of Theorem \ref{thprodIabsmon}]
Observe that if $t \in \left(0,\frac{\pi}{2}\right]$ and $a,b>0,$ then we have that $a^2+b^2-2ab\cos t>(a-b)^2>0$ and if $t \in \left(\frac{\pi}{2}, \pi\right)$ and $a,b >0$ then we also have that $a^2+b^2-2ab\cos t>0.$ Replacing $a=uv$ and $b={u}/{v}$ in the above integral representation \eqref{intIprod}, we obtain that
$$I_{\mu}(uv)I_{\mu}\left(\frac{u}{v}\right)=\frac{\left(\frac{1}{2}u^2\right)^{\mu}}{\sqrt{\pi}\Gamma(\frac{1}{2}+\mu)}\int_{0}^{\pi}f_{\mu}(S) \sin^{2\mu}t dt,$$
where $f_{\mu}(S)=S^{-\mu}I_{\mu}(S)$ and $S=\sqrt{u^2[(w^2-2)-2\cos t]}>0.$ In view of the recurrence relation $\left[z^{-\mu}I_{\mu}(z)\right]'=z^{-\mu}I_{\mu+1}(z)$ (see \cite[p. 79]{watson})
we obtain that each of the expressions
$$\frac{d}{dw}f_{\mu}(S)=u^2wf_{\mu+1}(S), \quad \frac{d^2}{dw^2}f_{\mu}(S)=u^2f_{\mu+1}(S)+(u^2w)^2f_{\mu+2}(S),$$
$$\frac{d^3}{dw^3}f_{\mu}(S)=3wu^4f_{\mu+2}(S)+(u^2w)^3f_{\mu+3}(S), \quad \frac{d^4}{dw^4}f_{\mu}(S)=3u^4f_{\mu+2}(S)+6u^6w^2f_{\mu+3}(S)+(u^2w)^4f_{\mu+4}(S)$$
are positive for all $a,b>0,$ $w>2$ and $\mu>-\frac{1}{2}.$ In view of the above relations we make the induction hypothesis that the $(2n+1)$th order derivative is of the form
$$\frac{d^{2n+1}}{dw^{2n+1}}f_{\mu}(S)=(u^2w)^{2n+1}f_{\mu+2n+1}(S)+\alpha_{2n}(u)w^{2n-1}f_{\mu+2n}(S)+\ldots+\alpha_2(u)wf_{\mu+k}(S)$$
and this expression is positive, where $k\leq 2n$ and the constants $\alpha_{2n}(u),$ $\alpha_{2n-2}(u),$ $\ldots,$ $\alpha_2(u)$ are non-negative. Moreover, we also make the
induction hypothesis that the $2n$th order derivative is of the form
$$\frac{d^{2n}}{dw^{2n}}f_{\mu}(S)=(u^2w)^{2n}f_{\mu+2n}(S)+\beta_{2n}(u)w^{2n-1}f_{\mu+2n-1}(S)+\ldots+\beta_2(u)f_{\mu+k}(S)$$
and this expression is positive, where $k\leq 2n-1$ and the constants $\beta_{2n}(u),$ $\beta_{2n-2}(u),$ $\ldots,$ $\beta_2(u)$ are non-negative. By using the recurrence relation ${d}f_{\mu}(S)/{dw}=u^2wf_{\mu+1}(S)$ repeatedly we can see that the derivatives ${d^{2n+3}f_{\mu}(S)}/{dw^{2n+3}}$ and ${d^{2n+2}f_{\mu}(S)}/{dw^{2n+2}}$ will have a similar form as the derivatives ${d^{2n+1}f_{\mu}(S)}/{dw^{2n+1}}$ and ${d^{2n}f_{\mu}(S)}/{dw^{2n}},$ respectively, moreover each of these expressions are positive. Consequently, for all $\mu>-\frac{1}{2},$ $u,v>0,$ $w>2$ and $n\in\mathbb{N}$ we have that
$$\frac{d^n}{dw^n}\left[I_{\mu}(uv)I_{\mu}\left(\frac{u}{v}\right)\right]=\frac{\left(\frac{1}{2}u^2\right)^{\mu}}{\sqrt{\pi}\Gamma(\frac{1}{2}+\mu)}\int_{0}^{\pi}\frac{d^n}{dw^n}f_{\mu}(S) \sin^{2\mu}t dt \geq 0.$$
\end{proof}

\subsection*{Acknowledgments} The research of Dhivya Prabhu K was supported in part by the Council of Scientific and Industrial Research India (Grant No. 09/1022(11054)/2021-EMR-I) and in part by Short-term Fellowship of the Indian Institute of Technology Indore, and was completed during his visit between March and May 2024 to Babe\c{s}-Bolyai University of Cluj-Napoca. This author is grateful to
the Department of Business Administration, Extension Sf\^{a}ntu-Gheorghe, Faculty of Economics and Business Administration, Babe\c{s}-Bolyai University of Cluj-Napoca for hospitality.

\end{document}